\magnification 1200 
\baselineskip 5mm
\hsize=160truemm
\vsize=236truemm
\hoffset=5truemm
\voffset=5truemm
\hfuzz=10pt

\expandafter\ifx\csname amssym.def\endcsname\relax \else\endinput\fi
%
\expandafter\edef\csname amssym.def\endcsname{%
       \catcode`\noexpand\@=\the\catcode`\@\space}
\catcode`\@=11
%

\def\undefine#1{\let#1\undefined}
\def\newsymbol#1#2#3#4#5{\let\next@\relax
 \ifnum#2=\@ne\let\next@\msafam@\else
 \ifnum#2=\tw@\let\next@\msbfam@\fi\fi
 \mathchardef#1="#3\next@#4#5}
\def\mathhexbox@#1#2#3{\relax
 \ifmmode\mathpalette{}{\m@th\mathchar"#1#2#3}%
 \else\leavevmode\hbox{$\m@th\mathchar"#1#2#3$}\fi}
\def\hexnumber@#1{\ifcase#1 0\or 1\or 2\or 3\or 4\or 5\or 6\or 7\or 8\or
 9\or A\or B\or C\or D\or E\or F\fi}

\font\tenmsa=msam10
\font\sevenmsa=msam7
\font\fivemsa=msam5
\newfam\msafam
\textfont\msafam=\tenmsa
\scriptfont\msafam=\sevenmsa
\scriptscriptfont\msafam=\fivemsa
\edef\msafam@{\hexnumber@\msafam}
\mathchardef\dabar@"0\msafam@39
\def\dashrightarrow{\mathrel{\dabar@\dabar@\mathchar"0\msafam@4B}}
\def\dashleftarrow{\mathrel{\mathchar"0\msafam@4C\dabar@\dabar@}}

\def\ulcorner{\delimiter"4\msafam@70\msafam@70 }
\def\urcorner{\delimiter"5\msafam@71\msafam@71 }
\def\llcorner{\delimiter"4\msafam@78\msafam@78 }
\def\lrcorner{\delimiter"5\msafam@79\msafam@79 }
\def\yen{{\mathhexbox@\msafam@55}}
\def\checkmark{{\mathhexbox@\msafam@58}}
\def\circledR{{\mathhexbox@\msafam@72}}
\def\maltese{{\mathhexbox@\msafam@7A}}

\font\tenmsb=msbm10
\font\sevenmsb=msbm7
\font\fivemsb=msbm5
\newfam\msbfam
\textfont\msbfam=\tenmsb
\scriptfont\msbfam=\sevenmsb
\scriptscriptfont\msbfam=\fivemsb
\edef\msbfam@{\hexnumber@\msbfam}

\def\widehat#1{\setbox\z@\hbox{$\m@th#1$}%
 \ifdim\wd\z@>\tw@ em\mathaccent"0\msbfam@5B{#1}%
 \else\mathaccent"0362{#1}\fi}
\def\widetilde#1{\setbox\z@\hbox{$\m@th#1$}%
 \ifdim\wd\z@>\tw@ em\mathaccent"0\msbfam@5D{#1}%
 \else\mathaccent"0365{#1}\fi}
\font\teneufm=eufm10
\font\seveneufm=eufm7
\font\fiveeufm=eufm5
\newfam\eufmfam
\textfont\eufmfam=\teneufm
\scriptfont\eufmfam=\seveneufm
\scriptscriptfont\eufmfam=\fiveeufm

\let\goth\frak

\csname amssym.def\endcsname

\catcode`\@=11


\pretolerance=500 \tolerance=1000 \brokenpenalty=5000

\catcode`\;=\active
\def;{\relax\ifhmode\ifdim\lastskip>\z@
\unskip\fi\kern.2em\fi\string;}

\catcode`\:=\active
\def:{\relax\ifhmode\ifdim\lastskip>\z@
\unskip\fi\penalty\@M\ \fi\string:}

\catcode`\!=\active
\def!{\relax\ifhmode\ifdim\lastskip>\z@
\unskip\fi\kern.2em\fi\string!}

\catcode`\?=\active
\def?{\relax\ifhmode\ifdim\lastskip>\z@
\unskip\fi\kern.2em\fi\string?}

\frenchspacing



\newsymbol\ltimes 226E
\newsymbol\rtimes 226F

\newtoks\auteurcourant
\newtoks\titrecourant

\font\bb= msbm10 
\font\sbb=msbm10 at 7pt
\font\ssbb=msbm10 at 5pt
\newfam\cmfam\textfont\cmfam=\bb\scriptfont\cmfam=
\sbb\scriptscriptfont\cmfam=\ssbb

\font\got=eufm10 
\font\sgot=eufm10 at 7pt
\font\ssgot=eufm10 at 5pt
\newfam\gfam\textfont\gfam=\got\scriptfont\gfam=\sgot
\scriptscriptfont\gfam=\ssgot
\def\goth{\fam\gfam\got}

\font\tbf=cmbx10 at 12pt

\font\bftwelve=cmbx10 at 12pt

\font\petcap=cmcsc10

\def\truc{\unskip\kern 3pt\penalty 500
\hbox{\vrule\vbox to 5pt{\hrule width 4pt\vfill\hrule}\vrule}\kern
3pt}

\def\qed{\hfill $\truc$\par\goodbreak\bigskip}

\def\r#1{\hbox{$\cal #1$}}

\def\g#1{\hbox{\goth #1}}
\def\sg#1{\hbox{\sgot #1}}

\def\tch#1#2{{#1\check{#2}}}

\def\bc{\buildrel\circ\over}\null

\def\parni{\par\noindent}

\def\N{\hbox{\bb N}}
\def\sN{\hbox{\sbb N}}

\def\Z{\hbox{\bb Z}}

\def\R{\hbox{\bb R}}
\def\sR{\hbox{\sbb R}}

\def\C{\hbox{\bb C}}

\def\A{\hbox{\bb A}}
\def\sA{\hbox{\sbb A}}

\def\qa{\alpha}     	\def\qb{\beta}	    	\def\qd{\delta}
\def\qe{\varepsilon} \def\qf{\varphi}    \def\qg{\gamma}
       \def\qk{\kappa}    	\def\ql{\lambda} 
\def\qm{\mu}         \def\qn{\nu}        \def\qo{\omega}  
\def\qp{\pi}									\def\qr{\rho}      	 	
 	            
\def\qx{\xi}         
\def\qD{\Delta}	     \def\qF{\Phi}       \def\qG{\Gamma}	    
     \def\qO{\Omega} 	           
\def\qS{\Sigma}

\def\cf{{\it cf.\/}\ }
\def\ie{{\it i.e.\/}\ }
\def\eg{{\it e.g.\/}\ }

\def\lc{{\it l.c.\/}\ }
\def\LC{{\it loc. cit.\/}\ }

\def\ref#1{\par [\hbox to 0,95 cm{#1] \hfil}}

\def\^#1{\if#1i{\accent"5E\i}\else{\accent"5E #1}\fi}
\def\"#1{\if#1i{\accent"7F\i}\else{\accent"7F #1}\fi}

\catcode`\@=12



\headline={\ifnum \pageno=1 {\hfill}\else{\ifodd \pageno {\hfill\tenrm\the
\titrecourant\hfill\tenbf\folio}\else {\tenbf\folio\hfill\tenrm\the
\auteurcourant\hfill}\fi}\fi}

\footline={\hfil}



\auteurcourant= {Guy Rousseau}
\titrecourant= {Masures affines}

\centerline{\tbf Masures affines}

\bigskip
\centerline{par}
\medskip
\centerline{\petcap   Guy Rousseau}
\bigskip
\vskip 0.5cm   \noindent Octobre 2008 

\vskip 0.5cm   \noindent
{\bf Abstract.} {\sl We give an abstract definition of affine hovels which
generalizes the definition of affine buildings (eventually non simplicial) given
by Jacques Tits and includes the hovels built by St\'ephane Gaussent and the author
for some Kac-Moody groups over ultrametric fields. We prove that, in such an affine
hovel \r I, there exist retractions with center a sector germ and that we can add
at the infinity of \r I a pair of twin buildings or two microaffine buildings. For
some affine hovels \r I, we prove that the residue at a point of \r I has a natural
structure of pair of twin buildings and that there exists on
\r I a preorder which induces on each apartment the preorder associated to the
Tits cone.}

\vskip 0.5cm   \noindent
{\bf R\'esum\'e.} {\sl On donne une d\'efinition abstraite de masure affine qui
g\'en\'eralise celle d'immeuble affine (\'eventuellement non simplicial) donn\'ee par
Jacques Tits et qui inclut les masures (hovels) construits par St\'ephane Gaussent
et l'auteur pour certains groupes de Kac-Moody sur des corps ultram\'etriques.
On montre que, dans une telle masure affine, il existe des r\'etractions de centre
un germe de quartier et qu'\`a l'infini on peut construire une paire d'immeubles
jumel\'es et deux immeubles microaffines. Pour certaines masures affines \r I, on
montre que le r\'esidu en chaque point de \r I a une structure naturelle de paire
d'immeubles jumel\'es et qu'il existe sur \r I un pr\'eordre qui induit sur chaque
appartement le pr\'eordre associ\'e au c\^one de Tits.}

\vskip 0.6cm\noindent

{\tbf Introduction.}
\bigskip

\par Si $G$ est un groupe de Kac-Moody d\'eploy\'e sur un corps ultram\'etrique,
St\'ephane Gaussent et l'auteur ont introduit un espace \r I sur lequel le groupe
$G(K)$ agit [Gaussent-Rousseau-08]; il y a en fait (pour l'instant?) des conditions
techniques restrictives sur $G$ et $K$, voir \S\kern 2pt 6 ci-dessous. La construction est
calqu\'ee, d'aussi pr\`es que possible, sur celle des immeubles de Bruhat-Tits pour
les groupes r\'eductifs sur les corps locaux [Bruhat-Tits-72]. Cette g\'en\'eralisation
au cas Kac-Moody pose un certain nombre de probl\`emes techniques et r\'ev\`ele une
difficult\'e majeure: cet espace \r I est r\'eunion d'appartements, mais l'axiome
fondamental des syst\`emes d'appartements n'est pas v\'erifi\'e: il peut exister deux
points de \r I qui ne sont pas dans un m\^eme appartement. \`A cause de cette
mauvaise propri\'et\'e, \r I s'est vu attribuer le nom de masure (hovel).

\par Par ailleurs Jacques Tits a donn\'e dans [Tits-86] une d\'efinition abstraite des
immeubles affines (\'eventuellement non discrets) sans supposer cet axiome
fondamental des syst\`emes d'appartements (ce qui est par contre le cas de la
d\'efinition de l'auteur dans [Rousseau-08]). L'un de ses axiomes essentiels est
l'existence d'un appartement contenant n'importe quelle paire de germes de
quartier. Comme cette propri\'et\'e et d'autres analogues sont v\'erifi\'ees par les
masures de [Gaussent-Rousseau-08], il est naturel d'envisager la g\'en\'eralisation
de la d\'efinition de Tits aux masures affines [\lc; remark 4.6].

\par On va donc, dans cet article, essayer de g\'en\'eraliser l'essentiel des
r\'esultats de [Tits-86] dans un cadre incluant les masures d\'ej\`a construites. On
introduit donc une d\'efinition abstraite de masure affine en imposant comme
axiomes un certain nombre de propri\'et\'es de ces masures concr\`etes, choisies comme
de bonnes g\'en\'eralisations des axiomes de [Tits-86] \cf 2.5.

\par On introduit d'abord au \S\kern 2pt 1 l'appartement affine t\'emoin $\A$ de la masure
affine. On part, essentiellement, de la repr\'esentation g\'eom\'etrique classique d'un
groupe de Coxeter (\'eventuellement infini) $W^v$ dans un espace vectoriel r\'eel $V$
de dimension finie (mais pas forc\'ement euclidien) [Bourbaki-68; V \S\kern 2pt 4]. Les
g\'en\'erateurs fondamentaux de $W^v$ sont des r\'eflexions lin\'eaires par rapport aux
murs d'un c\^one convexe, ouvert et simplicial $C^v_f$ (la chambre vectorielle
fondamentale positive). Les conjugu\'es par $W^v$ de ces murs sont les murs
vectoriels, leur ensemble est not\'e $\r M^v$; \`a chaque $M^v\in\r M^v$ est associ\'ee
une unique r\'eflexion $r_{M^v}\in W^v$ d'hyperplan $M^v$. Les conjugu\'es par $W^v$
de $C^v_f$ (ou de ses facettes) sont les chambres (ou les facettes) vectorielles
positives. La r\'eunion (disjointe) de toutes ces facettes est le c\^one de Tits \r T
qui est convexe. Les chambres ou facettes n\'egatives sont les images par $-Id_V$
des chambres ou facettes positives. De plus on est souvent oblig\'e de consid\'erer des
murs vectoriels "imaginaires" dont l'ensemble est stable par $W^v$, \cf 1.1.5.

\par L'appartement $\A$ est un espace affine sous $V$, il est muni d'un ensemble
\r M d'hyperplans affines appel\'es murs. Si $M\in \r M$, sa direction $M^v$ est
dans $\r M^v$ et on note $r_M$ la r\'eflexion d'hyperplan $M$ dont l'application
lin\'eaire associ\'ee est $r_{M^v}$. Le groupe de Weyl $W$, engendr\'e par les $r_M$
pour $M\in\r M$, doit stabiliser \r M. Les murs imaginaires sont les hyperplans
de direction un mur vectoriel imaginaire. Comme l'ensemble des murs n'est pas
localement fini, on ne peut pas d\'efinir des facettes dans $\A$ comme des
sous-ensembles de $\A$: ce sont des filtres de parties de $\A$ comme dans
[Bruhat-Tits-72], \cf 1.7. Un quartier (resp. une face de quartier) de
$\A$ est une partie de $\A$ de la forme $\g f=x+F^v$ pour $x\in\A$ et $F^v$ une
chambre (resp. une facette) vectorielle. On a besoin d'\'epaissir ces faces de
quartier: une chemin\'ee est un ensemble de la forme $\g r=F+{\overline{F^v}}$ (ou
plus pr\'ecis\'ement son enclos \cf 1.10) o\`u $F$ est une facette et $F^v$ une facette
vectorielle. La face de quartier $x+F^v$ (resp. la chemin\'ee
$cl(F+{\overline{F^v}})$) est dite sph\'erique (resp. \'evas\'ee) si $F^v$ est sph\'erique
\ie contenue dans l'int\'erieur de $\pm \r T$. Le germe d'une face de quartier $x+F^v$
(resp. d'une chemin\'ee $cl(F+{\overline{F^v}})$) est le filtre des parties de $\A$
contenant $x+F^v+\qx$ (resp.  $cl(F+{\overline{F^v}}+\qx)$) pour un certain $\qx\in
{\overline{F^v}}$.

\par Apr\`es ces d\'efinitions techniques, on d\'efinit au \S\kern 2pt 2 une masure affine comme
un ensemble \r I muni d'un recouvrement par un ensemble \r A d'appartements, tous
isomorphes \`a $\A$. On peut r\'esumer l'essentiel des axiomes en disant que l'on
demande de plus qu'une facette, une facette et un germe de chemin\'ee \'evas\'ee ou
deux germes de chemin\'ees \'evas\'ees sont contenus dans un m\^eme appartement, unique \`a
isomorphisme fixant ces objets pr\`es (2.1). Les immeubles de Bruhat-Tits, les
immeubles affines simpliciaux (au sens classique, voir par exemple [Brown-89]) et
vraisemblablement tous les immeubles affines de [Tits-86] sont des masures affines
(2.4). On d\'eduit aussit\^ot de la d\'efinition l'existence dans une masure affine de
r\'etractions de centre un germe de quartier (2.6).

\par Au \S\kern 2pt 3 on commence l'\'etude des propri\'et\'es \`a l'infini d'une masure affine.
Deux germes de faces de quartier sph\'eriques sont dits parall\`eles si, dans un
appartement les contenant tous deux, ils correspondent \`a la m\^eme facette
vectorielle. On montre que ceci d\'efinit une relation d'\'equivalence (3.2) et que
les classes d'\'equivalences constituent les facettes sph\'eriques de deux
immeubles jumel\'es (th\'eor\`emes 3.4 et 3.7). On obtient ainsi, \`a l'infini d'une
masure affine, une paire d'immeubles jumel\'es qui g\'en\'eralise donc l'immeuble
sph\'erique \`a l'infini d'un immeuble affine construit dans [Tits-86].

\par Au \S\kern 2pt 4 on examine la structure de l'ensemble des germes de faces de quartier
sph\'eriques positives (resp. n\'egatives). On montre que cela constitue un immeuble
microaffine au sens de [Rousseau-06]; les facettes de cet immeuble microaffine
correspondent aux germes de chemin\'ees \'evas\'ees positives (resp. n\'egatives) \cf4.4.
Les germes de faces de quartier sph\'eriques dans une m\^eme classe de parall\'elisme
constituent un immeuble affine (4.3). Dans le cas des cloisons de quartier, on
obtient ainsi des arbres qui g\'en\'eralisent ceux utilis\'es par J. Tits dans
[Tits-86] pour d\'efinir une valuation sur la donn\'ee radicielle associ\'ee (en
g\'en\'eral) \`a l'immeuble sph\'erique \`a l'infini d'un immeuble affine. On n'arrive
cependant pas ici \`a une classification des masures affines (qui g\'en\'eraliserait
[Tits-86]) pour (au moins) deux raisons: la classification des immeubles jumel\'es
n'est pas aussi g\'en\'erale que celle des immeubles sph\'eriques et surtout on ne sait
(encore) pas construire une masure affine associ\'ee \`a une donn\'ee radicielle valu\'ee
correspondant \`a un syst\`eme de racines infini (4.12).

\par Au \S\kern 2pt 5 on s'int\'eresse aux propri\'et\'es \`a distance finie d'une masure \r I. On
montre (5.3) qu'en chaque point de \r I il existe un r\'esidu sous la forme de deux
immeubles. Si la masure \r I est "ordonn\'ee" (2.3), ces deux immeubles sont en
fait jumel\'es (5.6). On a donc encore une paire d'immeubles jumel\'es qui g\'en\'eralise
l'immeuble sph\'erique qu'est le r\'esidu pour un immeuble affine. Dans l'appartement
$\A$ il y a un pr\'eordre naturel donn\'e par: $x\le y$ si et seulement si $y-x\in\r T$.
On d\'efinit une relation dans \r I par $x\le y$ si et seulement si il existe un
appartement $A$ contenant $x$ et $y$ et si $x\le y$ dans cet appartement. On montre
que cette relation est bien d\'efinie et que c'est un pr\'eordre si la masure est
ordonn\'ee (th\'eor\`eme 5.9). Dans le cas des immeubles affines cette relation est
triviale (car $\r T=V$); ce n'est pas le cas en g\'en\'eral et on a donc ainsi une
structure d'un type nouveau sur les masures affines.

\par Enfin au \S\kern 2pt 6 on exhibe des exemples de masures affines qui ne sont pas des
immeubles affines, en montrant que (presque toutes) les masures de
[Gaussent-Rousseau-08] sont des masures affines, ordonn\'ees, \'epaisses et
semi-discr\`etes (th\'eor\`eme 6.11).

\bigskip
\parni {\S$\,${\bf 1.}$\quad${\bftwelve Appartements affines }}
\bigskip
\par Dans ce premier paragraphe on laisse beaucoup de d\'emonstrations au lecteur.
Sauf indication particuli\`ere, ce sont de simples cons\'equences de l'alg\`ebre
lin\'eaire et de la th\'eorie des in\'equations lin\'eaires dans $\R^n$.

\bigskip \noindent
{\bf 1.1. Le groupe de Weyl vectoriel et les racines} 

\par Dans tout cet article on se donne un quadruplet
$(V,W^v,(\alpha_i)_{i\in I},(\tch{\alpha}{_i})_{i\in I})$ form\'e d'un espace
vectoriel r\'eel $V$ de dimension finie, d'un sous-groupe $W^v$ de $GL(V)$,
d'une famille $(\tch{\alpha}{_i})_{i\in I}$ dans $V$ et d'une famille libre
$(\alpha_i)_{i\in I}$ dans le dual $V^*$ de $V$ (toutes deux index\'ees par le m\^eme
ensemble fini $I$). On suppose v\'erifi\'ees les propri\'et\'es suivantes:

\par 1) Pour $i\in I$, $\alpha_i(\tch{\alpha}{_i})=2$ et donc la formule 
$r_i(v) = v - \alpha_i(v)\tch{\alpha}{_i} $ d\'efinit une involution de $V$ , 
plus pr\'ecis\'ement une r\'eflexion d'hyperplan $Ker(\alpha_i)$.

\par 2) Le groupe de Weyl $W^v$ est le sous-groupe de $GL(V)$ de syst\`eme
g\'en\'erateur $S=\{r_i\mid i\in I\}$. C'est un groupe de Coxeter, appel\'e {\it groupe
de Weyl vectoriel}. Le coefficient $m(i,j)$ de la matrice de Coxeter est l'ordre
de $r_ir_j$.

\par 3) On note $\Phi$ l'ensemble des {\it racines} r\'eelles c'est \`a dire
 des formes lin\'eaires sur $V$ de la forme $\alpha = w(\alpha_i)$ avec 
$w \in W^v$ et $ i \in I$. Si $\alpha \in \Phi$, alors $r_{\alpha} = 
w.r_i.w^{-1} $ est bien d\'etermin\'e par $\alpha$, ind\'ependamment du choix 
de $w$ et de $i$ tels que $\alpha = w(\alpha_i)$. Pour $v\in V$ on a 
$r_{\alpha}(v) = v - \alpha(v)\tch{\alpha}{ } $ pour un $\alpha\check{ }  \in V $
avec $\alpha(\alpha\check{ }) = 2 $ ; ainsi $r_{\alpha}$ est une r\'eflexion par 
rapport \`a l'hyperplan $M^v(\alpha) = Ker(\alpha)$ que l'on appelle {\it mur
(vectoriel)}  de $\alpha$ . On note ${\cal M}^v$ l'ensemble de ces murs
vectoriels (r\'eels). Le {\it demi-appartement vectoriel} associ\'e
\`a $\alpha$ est $D(\alpha) = \{ v \in V \mid \alpha(v) \geq 0 \} $.

\par 4) Si $\Phi^+ = \Phi\cap(\bigoplus_{i\in I}\; \R^+\alpha_i)$ et $\Phi^- = -
\Phi^+$ , on a $\Phi = \Phi^+ \bigsqcup \Phi^-$. Le groupe $W^v$ permute $\qF$,
mais seul l'\'el\'ement neutre stabilise $\qF^+$. Plus pr\'ecis\'ement pour $w\in W^v$, 
$w(\qa_i)\in\qF^+\Leftrightarrow\ell(wr_i)>\ell(w)$ o\`u $\ell$ d\'esigne la longueur
dans $W^v$ relative \`a $S$. Les racines de $\qF^+$ (resp. $\qF^-$) sont dites {\it
positives} (resp. {\it n\'egatives}). Pour $\qa\in\qF$, on a
$\qF\cap\R^+\qa=\{\qa\}$; ainsi $\qF^+$ (resp. $\qF$) est en bijection avec
${\cal M}^v$ (resp. les demi-appartements de $V$).

\par5) On consid\`ere un ensemble $\qD_{im}\subset V^*$ de racines imaginaires
disjoint de $\R\qF$ et stable par $\pm W^v$; on note $\qD_{re}=\qF$ et
$\qD=\qF\cup\qD_{im}$. Trois cas particuliers sont particuli\`erement int\'eressants:

\par $\bullet$ Cas sans imaginaire: $\qD_{im}=\emptyset$; c'est un cas particulier
du cas suivant.

\par $\bullet$ Cas mod\'er\'ement imaginaire: Si
$\qD_{im}^+=\qD_{im}\cap(\bigoplus_{i\in I}\; \R^+\qa_i)$ et
$\qD_{im}^-=-\qD_{im}^+$, on a
$\qD_{im}=\qD_{im}^+\cup\qD_{im}^-$ et $\qD_{im}^+$ comme $\qD_{im}^-$ sont stables
par $W^v$. On peut par exemple prendre $\qD_{im}^+=\cap_{w\in W^v}\;(\bigoplus_{i\in
I}\; \R^+\qa_i)\setminus \R\qF$.

\par $\bullet$ Cas tr\`es imaginaire: $\R^*\qD_{im}=V^*\setminus\R\qF$, par exemple
$\qD_{im}=V^*\setminus\R\qF$.

\par On note $\r M^{vi}$ l'ensemble des {\it murs vectoriels imaginaires}
$M^v(\qa)=$Ker$(\qa)$, pour $\qa\in\qD_{im}$.

\bigskip \noindent
{\bf 1.2. Exemples}
\par La situation d\'ecrite ci-dessus se rencontre dans (au moins) deux cas (non
disjoints):

\par 1) {\it Cas Kac-Moody}: \cf [Kac-90] ou [Rousseau-06]

\par On consid\`ere une r\'ealisation $(V,(\alpha_i)_{i\in I},
(\tch{\alpha}{_i})_{i\in I})$ d'une matrice de Kac-Moody $A$ (cela signifie que 
$\alpha_j(\tch{\alpha}{_i})$ est le coefficient $a_{i,j}$ de la matrice $A$
[Kac-90]) avec $(\alpha_i)_{i\in I}$ libre dans $V^*$. Alors $W^v$ en est le
groupe de Weyl, $\qF$ l'ensemble des racines r\'eelles et $\qD_{im}$ l'ensemble des
racines imaginaires; on est dans le cas mod\'er\'ement imaginaire. De mani\`ere
\'equivalente on peut consid\'erer un syst\`eme g\'en\'erateur de racines
$(A,X,Y,(\alpha_i)_{i\in I}, (\tch{\alpha}{_i})_{i\in I})$ au sens de [Bardy-96]
avec $V=Y\otimes\R$ comme dans [Rousseau-06]. Ce cas existe si et seulement si les
coefficients de la matrice de Coxeter de $W^v$ sont $1$, $2$, $3$,
$4$, $6$ ou $\infty$. De plus on a $\qD\subset \bigoplus_{i\in I}\;\Z\alpha_i$.

\par Si $W^v$ est fini, la matrice de Kac-Moody est une matrice de Cartan et
$(V,\qF,W^v)$ est associ\'e \`a une alg\`ebre de Lie r\'eductive complexe, ce cas sera
dit {\it classique}. Si $A$ est une matrice de Kac-Moody de type affine
[Kac-90; IV], on parlera de cas {\it affine}.

\par On peut modifier ce cas Kac-Moody en se pla{\c c}ant dans les cas sans imaginaire
ou tr\`es imaginaire.

\par 2) {\it Cas Coxeter g\'en\'eral}:
\parni \cf [Bourbaki-68; V \S\kern 2pt 4], [Deodhar-89], [Humphreys-90; V] ou [Tits-88].

\par On consid\`ere un groupe de Coxeter quelconque $W^v$ de syst\`eme fondamental de
g\'en\'erateurs $S=\{s_i\mid i\in I\}$ fini et matrice de Coxeter $(m(i,j))_{i,j\in
I}$. On lui associe un espace vectoriel r\'eel $V^*$ de base $(\alpha_i)_{i\in I}$, 
que l'on munit de la forme bilin\'eaire sym\'etrique $B$, telle que $B(\qa_i,\qa_i)=1$ 
et $B(\qa_i,\qa_j)=-cos({\qp\over m(i,j)})$. On d\'efinit alors 
$\tch{\alpha}{_i}\in V$ par $2B(f,\qa_i)=f(\tch{\alpha}{_i})$ pour tout  $f\in
V^*$. D'apr\`es les r\'ef\'erences ci-dessus les propri\'et\'es de 1.1 sont v\'erifi\'ees. Pour
les racines imaginaires, on peut envisager les 3 cas de 1.1.5.

\par Bien sur on peut aussi multiplier $V$ par un espace vectoriel sur lequel
$W^v$ n'agit pas.

\par 3) Remarques: On pourrait sans doute abandonner dans les hypoth\`eses de 1.1 les
conditions de libert\'e de la famille $(\qa_i)_{i\in I}$ et de finitude de $I$ ou de
la dimension de $V$. Il faudrait prendre quelques pr\'ecautions suppl\'ementaires, par
exemple imposer que $C_f^v$ (d\'efini ci-dessous en 1.3) engendre l'espace  vectoriel
$V$. Des exemples de cette situation et leur int\'er\^et sont expliqu\'es dans
[Moody-Pianzola-95] et [Bardy-96].

 \par M\^eme sans abandonner ces deux conditions, les syst\`emes g\'en\'erateurs de racines
plus g\'en\'eraux de [Bardy-96] fournissent des exemples diff\'erents de ceux de 1),
ni mod\'er\'ement imaginaires, ni sans imaginaires, ni tr\`es imaginaires.

\par Dans la suite on s'int\'eresse essentiellement au cas mod\'er\'ement imaginaire (sans
l'imposer). Le cas sans imaginaire est id\'eal mais souvent semblable, \cf 1.8.3.

\bigskip \noindent
{\bf 1.3. Le c\^one de Tits} 

\par Les r\'esultats de ce num\'ero se d\'emontrent comme dans [Bourbaki-68; V \S\kern 2pt 4] ou
[Humphreys-90; V] \`a partir des hypoth\`eses de 1.1.

\par La {\it chambre fondamentale positive} $C^v_f = \{ v \in V \mid \alpha_i(v) >
0\quad\forall i \in I \} $ est un c\^one convexe ouvert non vide. Son adh\'erence
 est r\'eunion disjointe des {\it facettes vectorielles} $ F^v(J) = \{\; v \in V
\mid  \alpha_i(v)  = 0 \; \forall i \in J ;
 \; \alpha_i(v) > 0 \;	\forall i \notin J \} $ pour $J\subset I$. On a
$C^v_f=F^v(\emptyset)$ et $V_0=F^v(I)$ est un sous-espace vectoriel. Si
$V_0=\{0\}$, on dit que $W^v$ est {\it essentiel}.

\par On dit que la facette $ F^v(J)$ ou que $J$ est {\it sph\'erique} si 
$W^v(J) = \langle r_iÊ \mid i \in J \rangle$ est fini; c'est le cas de la chambre 
$C^v_f$ ou des {\it cloisons} $F^v(\{i\})$ , $\forall i \in I$. Plus g\'en\'eralement
les {\it chambres} (resp. {\it facettes, facettes sph\'eriques, cloisons}) positives
sont les transform\'ees par $W^v$  de celles d\'ej\`a d\'efinies. Le {\it type} de
$w.F^v(J)$ est $J$.

\par On dit que l'on est dans le cas {\it fini} si (de mani\`ere \'equivalente)
$\qF$ est fini, $W^v$ est fini ou $I$ (et donc toute facette) est sph\'erique.

\par Le {\it c\^one de Tits}  est la r\'eunion disjointe $\r T$ des 
facettes $w.F^v(J)$ pour $w$ dans $W^v$ et $J\subset I$. C'est 
un c\^ one convexe stable par $W^v$ . L'action de $W^v$ sur les chambres
 est simplement transitive. Le fixateur ou le stabilisateur de $F^v(J)$ est 
$W^v(J)$. Pour $x$ dans \r T , on note $F^v(x)$ la facette qui le contient.
Celle-ci est sph\'erique si et seulement si $x$ est int\'erieur \`a \r T: la r\'eunion
des facettes sph\'eriques est l'int\'erieur $\r T^\circ$ de \r T, c'est un c\^one
convexe ouvert non vide, appel\'e {\it c\^one de Tits ouvert}. De m\^eme l'adh\'erence
$\overline{\r T}$ de \r T est un c\^one convexe.

\par On consid\'erera aussi le {\it c\^one de Tits n\'egatif} $-\r T$ et toutes les
facettes, chambres, cloisons n\'egatives obtenues par l'op\'eration
$\qO\mapsto-\qO$.  

\par Il peut arriver qu'une facette soit positive et n\'egative,
c'est le cas par exemple de $V_0$. Dans le cas  
fini, on a $\r T=-\r T=V$ et toutes les facettes sont sph\'eriques, 
positives et n\'egatives. Dans le cas non fini on a $\r T^\circ\cap -\r T^\circ
=\emptyset$: une facette sph\'erique ne peut \^etre positive et n\'egative. Dans le cas
affine $\r T^\circ=\{ v \in V \mid \delta(v) > 0 \}$ avec $\qd\in\qD^+_{im}$ et $\r
T=V_0\cup\r T^\circ$.

\par Dans le cas mod\'er\'ement imaginaire toute racine imaginaire positive est
positive sur \r T. En particulier le cas mod\'er\'ement imaginaire fini est un cas sans
imaginaire.

\bigskip \noindent
{\bf 1.4. L'appartement affine t\'emoin $\A$} 

\par On choisit un espace affine $\A$ sous l'espace vectoriel $V$ muni d'une
collection \r M d'hyperplans affines (appel\'es {\it murs} ou {\it murs r\'eels}) telle
que:

\par a) pour tout $M\in\r M$, il existe une racine $\qa_M\in\qF$ telle que la
direction $M^v$ de $M$ soit le mur vectoriel $M^v(\qa_M)=Ker(\qa_M)$,

\par b) toute racine $\qa\in\qF$ peut s'\'ecrire $\qa=\qa_M$ pour une infinit\'e de
murs $M$,

\par c)  pour $M\in\r M$, on note $r_M$ la r\'eflexion d'hyperplan $M$ dont
l'application lin\'eaire associ\'ee est $r_{\qa_M}$ et on suppose que le groupe $W$
engendr\'e par les $r_M$ stabilise \r M.

\par d) On consid\`ere aussi l'ensemble $\r M^i$ des {\it murs imaginaires} qui sont
tous les hyperplans affines dont la direction est de la forme Ker$(\qa)$ pour
$\qa\in\qD_{im}$. Cet ensemble $\r M^i$ est stable par $W$ (et m\^eme par
le groupe $W^e_{\sR}$ d\'efini en 1.6.1).

\par Le groupe $W$ est le {\it groupe de Weyl} (affine) de $\A$.

\par La {\it dimension} (resp. le {\it rang}  (vectoriel)) de $\A$
est $dim(V)$ (resp. $\vert I\vert$).

\par On dit que $\A$ est {\it semi-discret} (resp. {\it dense}) si, pour toute
racine r\'eelle $\qa$, l'ensemble des murs de direction $Ker(\qa)$ est discret (resp.
dense). L'ensemble \r M n'est un ensemble discret d'hyperplans que si, de plus,
on est dans le cas fini; on dit alors que $\A$ est {\it discret}.

\par Un point $x\in\A$ est {\it sp\'ecial} si tout mur r\'eel admet un mur parall\`ele
passant par $x$ (il existe toujours de tels points). 

\par Si on choisit un point sp\'ecial $x_0$ comme origine, les murs r\'eels ou
imaginaires peuvent \^etre d\'ecrits comme les $M(\qa,k)=\{v\in\A\mid\qa(v)+k=0\}$ pour
$\qa\in\qD$ et $k\in\qG_{\qa}$, o\`u $\qG_{\qa}$ est un sous-ensemble de $\R$; si
$\qa$ est r\'eelle, $\qG_\qa$ est infini et contient $0$, si $\qa$ est imaginaire,
$\qG_\qa=\R$.  On note $D(\qa,k)=\{v\in\A\mid\qa(v)+k\ge 0\}$ (resp.
$D^\circ(\qa,k)= \{v\in\A\mid\qa(v)+k>0\}$); si $\qa$ est r\'eelle on l'appelle un
{\it demi-appartement} (resp. {\it demi-appartement-ouvert}). On note
$D(\qa,\infty)=D^\circ(\qa,\infty)=\A$ et, pour $\qa$ r\'eelle, 
$r_{\qa,k}=r_{M(\qa,k)}$.

\parni{\bf N.B.} L'ensemble $\r M^i$ des murs imaginaires ne servira en fait qu'\`a
d\'efinir l'enclos d'une partie de $\A$, \cf 1.7.2.

\bigskip \noindent
{\bf 1.5. Le groupe de Weyl $W$ et les relations de pr\'eordre sur $\A$} 

\par 1) Pour $\qa\in\qF$, il est clair que $W$ contient les translations de
vecteurs dans ${\tilde \qG}_\qa.\qa\,\check{ }$, o\`u ${\tilde \qG}_\qa$ est le
sous-groupe de $\R$ engendr\'e par $\qG_\qa$ et que $\qG_\qa=\qG_\qa+2{\tilde
\qG}_\qa$ (condition c), en particulier $\qG_\qa=-\qG_\qa$. Par contre on peut
avoir $\qG_\qa\ne {\tilde \qG}_\qa$ si
$\qG_\qa$ n'est pas discret [Bruhat-Tits-72; 6.2.17].

\par On note $Q \,\check{ }=\sum_{\qa\in\qF}\;{\tilde \qG}_\qa.\qa\,\check{ }$,
c'est un sous-groupe de $V$ stable par $W^v$ (car $\qG_{w\qa}=\qG_\qa$). On
l'identifie \`a un groupe de translations de $\A$ (contenu dans $W$). Pour
$\qa\in\qF$, 
$M=M(\qa,k)$ et $N=M(\qa,0)$, $r_M$ est le compos\'e de $r_N$ et de la
translation de vecteur ${k\qa\,\check{ }}$. Donc, si on identifie $W^v$ au
sous-groupe de $W$ fixant l'origine (sp\'eciale) $x_0$, on a une d\'ecomposition en
produit semi-direct:\qquad $W=W^v\ltimes Q \,\check{ }$.

\par 2) {\it Remarque} : Si $\A$ est discret, quitte \`a quotienter par
$V_0$ on est dans le cas essentiel et l'ensemble des points sp\'eciaux est discret
dans $\A$. Ainsi $Q \,\check{ }$ (qui les stabilise) est un sous-groupe discret de
$V$. Comme $W^v$ est fini, un produit scalaire $W^v-$invariant identifie $V$ \`a
$V^*$ de fa{\c c}on que $\qa_i$ soit colin\'eaire \`a $\qa\,\check{_i}$, donc $Q \,\check{
}$ est un r\'eseau de $V$. On sait alors que $W^v$ est cristallographique
[Bourbaki-68; VI 2.5]: les \'el\'ements de la matrice de Coxeter sont $1$, $2$, $3$,
$4$ ou $6$. On peut donc se ramener au cas classique de 1.2.1.

\par 3) On note $P\, \check{ }$  le sous-groupe de $V$ stabilisant \r M, il
contient $Q\, \check{ }$ et $V_0$. Alors $W_P=W^v\ltimes P \,\check{ }$ est le
plus grand  sous-groupe de $W_{\sR}^e=W^v\ltimes V$ stabilisant \r M.

\par 4) {\it Pr\'eordres}: Le c\^one de Tits \r T, son adh\'erence $\overline{\r T}$
et son int\'erieur ${\r T}^{\circ}$ sont des c\^ones convexes et
$W^v$ invariants, on obtient donc trois pr\'eordres $W_{\sR}^e-$invariants sur $\A$
en posant:
\par\quad
$x\le y\;\Leftrightarrow\; y-x\in\r T$\quad ,\quad
$x\overline{\le }y\;\Leftrightarrow\; y-x\in\overline{\r T}$\quad ,\quad
$x{\bc\le }y\;\Leftrightarrow\; y-x\in\r T^\circ\cup\{0\}$ .

\par On a souvent $\r T\cap-\r T\ne \{0\}$, ce qui implique que $\le $ ou
$\overline{\le }$ n'est pas une relation d'ordre. Par contre en dehors du cas fini
on a $\r T^\circ\cap-\r T^\circ=\emptyset$, donc ${\bc\le }$ est une relation
d'ordre. Dans le cas fini, on a $\r T^\circ=V$, donc $x\le y$ et $x{\bc\le }y$ pour tous
$x,y$ dans $\A$. Dans le cas affine essentiel  $\le $ et $\bc\le $ co\"{\i}ncident; de
plus pour tous $x,y$ dans $\A$ on a $x\overline{\le }y$ ou $x\ge y$ (et les deux \`a la fois
seulement pour $x=y$).


\bigskip \noindent
{\bf 1.6. Exemples} 

\par 1) On peut consid\'erer l'ensemble $\r M_{\sR}$ de tous les hyperplans de
direction dans $\r M^v$. On a donc 4 appartements affines (tr\`es li\'es) $\A$,
$\A_{\sR}$, $\A^{si}$ et $\A_{\sR}^{si}$ associ\'es aux syst\`emes de murs 
$\r M\cup\r M^i$, $\r M_{\sR}\cup\r M^i$, $\r M$ et $\r M_{\sR}$. L'appartement 
$\A_{\sR}$ ou $\A_{\sR}^{si}$ est dense et tous ses points sont sp\'eciaux. Le groupe
de Weyl (resp. le groupe $W_P$) correspondant est
$W_{\sR}=W^v\ltimes Q \,\check{ }{_{\sR}}$ (resp. $W_{\sR}^e=W^v\ltimes V$) o\`u $Q
\,\check{ }{_{\sR}}=\sum_{\qa\in\qF}\;{\R}.\qa\,\check{}$.

\par Un mur ou
demi-appartement  relatif \`a $\r M_{\sR}$ mais pas \`a \r M sera dit {\it
fant\^ome}; par opposition un mur ou demi-appartement  relatif \`a
$\r M$ sera parfois dit {\it vrai}.

\par 2) Dans le cas Kac-Moody sans imaginaire ou mod\'er\'ement imaginaire, on peut
choisir un sous-groupe fixe non trivial
$\qG$ de $\R$ et poser $\qG_\qa=\qG$, $\forall\qa\in\qF$. Alors $Q \,\check{ }
=\oplus_{i\in\ I}\;{\qG}.\qa\,\check{_i }$ et donc $W$ stabilise bien \r M.
L'appartement $\A$ est semi-discret (resp. dense) si et seulement si $\qG$
est discret (resp. dense) dans $\R$ et alors $Q \,\check{ }$ est discret
(resp. dense) dans $Q\,\check{ }{_{\sR}}$. Le cas $\qG=\Z$ est d\'evelopp\'e dans
[Gaussent-Rousseau-08].

\par Dans le cas mod\'er\'ement imaginaire, il est en fait naturel de poser, comme dans
\LC, $\qG_\qa=\qG$, $\forall\qa\in\qD_{im}$, ce qui est contraire \`a la convention de
1.4.d. Mais comme tout multiple entier d'une racine imaginaire $\qa$ est une racine,
l'ensemble des murs de direction Ker$\qa$ est alors dense dans $\A$. Comme la
d\'efinition ci-dessous de l'enclos permet de faire intervenir
$\bigcap_{n\in\sN^*}\;D(n\qa,k_{n\qa})$ pour des $k_\qb\in\qG$, cette convention et
celle de 1.4.d donnent les m\^emes enclos.

\par Ceci explique la d\'efinition simplificatrice de $\r M^i$ adopt\'ee en 1.4.d. On
notera cependant que les 2 conventions peuvent donner (en 1.7.3 ci-dessous) des
facettes l\'eg\`erement diff\'erentes (car une intersection d\'ecroissante de demi-espaces
ouverts est un demi-espace ferm\'e), m\^eme si les facettes ferm\'ees sont les m\^emes.

\par 3) Dans le cas Kac-Moody sans imaginaire, si $(\qf_\qa)_{\qa\in\qF}$ est une
valuation r\'eelle d'une donn\'ee radicielle $(G,(U_\qa)_{\qa\in\qF},T)$ (\cf
[Rousseau-06]), on note
$\qG_\qa=\qf_\qa(U_\qa\setminus\{1\})$. Si $\qF$ est classique, on d\'efinit
ainsi un appartement affine $\A$ [Bruhat-Tits-72]. Il est vraisemblable que ceci
reste vrai dans le cas Kac-Moody g\'en\'eral.

\bigskip \noindent
{\bf 1.7. Facettes et autres filtres} 

\par Les facettes de $\A$ sont associ\'ees aux murs et demi-appartements; mais,
comme dans [Bruhat-Tits-72], ce ne sont a priori des sous-ensembles de $\A$ que dans
le cas discret. En g\'en\'eral on doit les d\'efinir comme des filtres de parties de $\A$.
On adopte pour les filtres les conventions de [Rousseau-08] ou
[Gaussent-Rousseau-08].

\par 1) Pour $x\in\A$, on note $\r V_x$ le filtre des voisinages dans $\A$ du
point $x$. Si $\qO$ est un sous-ensemble de $\A$ contenant $x$ dans son adh\'erence,
le {\it germe de $\qO$ en $x$} est le filtre $germ_x(\qO)=\qO\cap\r V_x$ form\'e
des sous-ensembles de $\A$ contenant un voisinage de $x$ dans $\qO$. 

\par En particulier, pour $x\ne y\in\A$, $[x,y)=germ_x([x,y])$ (resp.
$]x,y)=germ_x(]x,y])$) est le germe en $x$  du segment $[x,y]$ (resp. de
l'intervalle $]x,y]=[x,y]\setminus\{x\}$). Ce germe est dit {\it positif} (resp. 
{\it n\'egatif}) si $x\le y$ (resp. $y\le x$). Si $x\le y$ ou $y\le x$, on dit que $[x,y]$, 
$[x,y)$ ou $]x,y)$ et la demi-droite d'origine $x$ contenant $y$ sont {\it
pr\'eordonn\'es}; ils sont {\it g\'en\'eriques} si la facette
$F^v(\pm (y-x))$ est sph\'erique \ie si $x\bc\le y$ ou $y\bc\le x$.

\par Si $F$ est un filtre de parties de $\A$ et $x\in\A$, on dit que $x\in F$ si
$x\in S$, $\forall S\in F$ \ie si $\{x\}\subset F$

\par Le {\it support} $supp(F)$ d'un filtre $F$ de parties de $\A$ est le plus
petit sous-espace affine $E$ de $\A$ contenant $F$; c'est aussi le support de
l'adh\'erence $\overline F$ de $F$. La {\it dimension} de $F$ est la dimension de
$E$. 


\par 2) {\it L'enclos} $cl(F)$ d'un filtre $F$ de parties de $\A$ est le filtre
form\'e des sous-ensembles de $\A$ contenant un \'el\'ement de $F$ de la forme
$\cap_{\qa\in\qD}\;D(\qa,k_\qa)$, avec pour chaque $\qa$,
$k_\qa\in\qG_\qa\cup\{\infty\}$ (en particulier chaque $D(\qa,k_\qa)$ contient
$F$). Un filtre est dit {\it clos} s'il est \'egal \`a son enclos. Tout enclos est clos.

\parni{\bf N.B.} a) Dans le cas fini non discret, l'enclos tel que d\'efini ici peut
\^etre l\'eg\`erement plus gros que celui de [Bruhat-Tits-72; 7.1.2]. Par exemple
l'enclos d'une partie $\qO$ de $\A$  peut n'\^etre qu'un filtre et non une partie de
\A ; mais l'intersection des \'el\'ements de $cl(\qO)$ est l'enclos tel que d\'efini
dans \lc

\par b) Cet enclos est plus petit que l'enclos d\'efini \`a partir de \r M uniquement
(cas sans imaginaire) que l'on notera $cl^{si}(F)$.

\par Le {\it $\R-$enclos} $cl_{\sR}(F)$ du filtre $F$ est d\'efini de la m\^eme
mani\`ere que l'enclos avec $\r M_{\sR}$ \`a la place de \r M. Dans le cas tr\`es
imaginaire le $\R-$enclos n'est rien d'autre que l'enveloppe convexe $conv(F)$ de
$F$. On note $cl_{\sR}^{si}(F)$ l'enclos de $F$ associ\'e \`a $\r M_{\sR}$ (sans $\r
M^i$). Bien sur les enclos $cl$, $cl_{\sR}$, $cl^{si}$ et $cl_{\sR}^{si}$ sont en
fait associ\'es \`a $\A$, $\A_{\sR}$, $\A^{si}$ et $\A_{\sR}^{si}$. Pour tout filtre
$F$ on a : $conv(F)\subset cl_{\sR}(F)\subset cl_{\sR}^{si}(F)\subset cl^{si}(F)$ et
 $cl_{\sR}(F)\subset cl(F)\subset cl^{si}(F)$.

\par 3) Une {\it facette-locale} de $\A$ est associ\'ee \`a un point $x\in\A$ et une
facette vectorielle $F^v$ dans $V$; c'est le filtre $F^\ell(x,F^v)
=germ_x(x+F^v)$. La {\it facette} associ\'ee \`a $F^\ell(x,F^v)$ est le filtre
$F(x,F^v)$ form\'e des ensembles contenant une intersection de demi-espaces
$D(\qa,k_\qa)$ ou $D^\circ(\qa,k_\qa)$ (un seul
$k_\qa\in\qG_\qa\cup\{\infty\}$ pour chaque $\qa\in\qD$), cette intersection devant
contenir $F^\ell(x,F^v)$. La facette-locale $F^\ell$ et la facette $F$ ont le  m\^eme
enclos qui est aussi la {\it facette-ferm\'ee} $\overline F$, adh\'erence de $F$.

\par Cette d\'efinition donne des facettes plus grandes (au sens large) que celles de
[Bruhat-Tits-72; 7.2.4] dans le cas fini. Ce sont les m\^emes facettes si l'on
remplace intersection par intersection finie en [\lc; 7.2.4 ligne 4]; cette
modification ne change pas le groupe associ\'e et donc pas l'ensemble des
appartements contenant la facette. Dans le cas discret, la facette est bien le
sous-ensemble de $\A$ consid\'er\'e par Bruhat et Tits.

\par Si toute racine de $\qD$ est nulle sur $V_0$, la facette associ\'ee \`a $\r
M_{\sR}$, $\r M^i$, $F^v$ et $x$ est 
$F_{\sR}(x,F^v)= F^\ell(x,F^v)+V_0$. Dans le cas dense on a
$F(x,F^v)=F_{\sR}(x,F^v)$ pour $x$ sp\'ecial.

\par 4) Une facette-locale ou une facette est dite {\it positive} ou {\it n\'egative}
si on peut l'\'ecrire  $F=F^\ell(x,F^v)$ ou $F=F(x,F^v)$ avec $F^v$ de ce signe; elle
peut \^etre positive et n\'egative.

\par Cette facette-locale ou facette $F$ est dite {\it sph\'erique} si la direction
de son support rencontre le c\^one de Tits ouvert (et donc a un fixateur fini dans
$W^v$), alors son fixateur $W_F$ dans $W$ est fini. Si $F^v$ est sph\'erique,  alors
$F^\ell(x,F^v)$ ou $F(x,F^v)$ l'est aussi et c'est \'equivalent pour une
facette-locale.

\par Si $F=F^\ell(x,F^v)$ est une facette-locale, on a $dim(F)=dim(F^v)$. Ce
n'est en g\'en\'eral pas vrai pour une facette ou une facette-ferm\'ee.

\par 5) {\bf Lemme.} {\sl Si $F$ est une facette, une facette-ferm\'ee, une facette
locale ou un filtre clos, alors tout $\qO\in F$ contient un ouvert non vide du
support de
$F$.}

\parni{\bf N.B.} Par contre le support d'un filtre clos n'est pas forc\'ement clos.

\parni{\bf D\'emonstration.} Une base du filtre $F$ est form\'ee d'ensembles convexes
et il est clair qu'un ensemble convexe contient un ouvert non vide de l'espace
affine qu'il engendre.\qed

\par 6) Dans le cas mod\'er\'ement imaginaire une facette-ferm\'ee est l'enclos d'un germe
de segment pr\'eordonn\'e.

\bigskip\noindent
{\bf 1.8. Chambres, cloisons...} 

\par 1) Il y a un ordre sur les facettes: les assertions "$F$ est une face de $F'$
" , "$F'$ domine $F$ " et "$F\le F'$ " sont par d\'efinition \'equivalentes \`a \quad
$F\subset{\overline{F'}}$.

\par Tout point $x\in\A$ est contenu dans (l'adh\'erence d')une unique facette
minimale 
$F(x,V_0)$ ; le point $x$ est un {\it sommet} si
et seulement si $F(x,V_0)=\{x\}$.

\par 2) Une {\it chambre} (ou alc\^ove) est une facette maximale ou, de mani\`ere
\'equivalente, une facette que l'on peut \'ecrire sous la forme $C=F(x,C^v)$ o\`u $C^v$
est une chambre vectorielle. Sa dimension est $n$ et elle est sph\'erique.

\par Une cloison est une facette maximale parmi les facettes contenues dans au
moins un mur (r\'eel), ou de mani\`ere \'equivalente une facette que l'on peut \'ecrire
sous la forme $F=F(x,F^v)$ o\`u $F^v$ est une cloison vectorielle et $x$ appartient \`a
un mur de direction $supp(F^v)$. Sa dimension est $n-1$, son support est un mur et
elle est sph\'erique. 


\par Un {\it mur d'une chambre} $C$ est le support $M$ d'une cloison $F$ domin\'ee
par $C$. Deux chambres sont dites {\it adjacentes} (le long de $M$ ou $F$) si
elles dominent une m\^eme cloison $F$ de support $M$. Mais il peut exister une
chambre ne dominant aucune cloison, et donc n'ayant aucun mur. Donc, $\A$ n'est
pas forc\'ement "connexe par galeries".

\par 3) {\bf Exemple}: On se place dans le cas (assez g\'en\'eral) de 1.2.1: le cas
Kac-Moody (sans imaginaire ou mod\'er\'ement imaginaire) avec des $\qG_\qa$ comme en
1.4 (suppos\'es seulement infinis). On va mieux d\'ecrire ci-dessous la facette
$F=F(x,F^v(J))$ m\^eme quand $x$ n'est pas sp\'ecial. On note
$\qF_x=\{\qa\in\qF\mid\qa(x)\in\qG_\qa\}$, $\qD(J)=\qD\cap(\oplus_{j\in
J}\;\R\qa_j)$ et $\qF(J)=\qF\cap\qD(J)$.

\par Si $J$ n'est pas sph\'erique, alors (que ce soit avec ou sans imaginaires) on a
dim$(F)\le n-1$. Plus pr\'ecis\'ement $supp(F)$ est contenu dans l'intersection $x+V'(J)$
des hyperplans affines $x+$Ker$\qb$ pour les $\qb\in V^*\setminus\{0\}$ tels que
$\R\qb$ soit point d'accumulation (dans ${\hbox{\bb P}}(V^*)$ ) des $\R\qg$ pour
$\qg\in\qF(J)$. L'espace $V'(J)$ est stable par $W^v(J)$. De plus toute racine
$\qd\in\qD_{im}\cap\qD(J)$ est nulle sur $V'(J)$ [Kac-90; ex. 5.12], donc $F$ n'est
pas sph\'erique.

\par Si $J$ est sph\'erique, alors $F$ est \'egalement sph\'erique. On a dim$(F)=n$ si et
seulement si $\qF_x\cap\qF(J)=\emptyset$ et alors $F=F(x,C^v_f)$ est une chambre.
On a dim$(F)=n-1$ si et seulement si $\vert\qF_x\cap\qF(J)\vert=1$; on peut supposer
que $\qF_x^+\cap\qF(J)=\{\qa_{j_0}\}$ pour $j_0\in J$ et alors $F=F(x,F^v(\{j_0\}))$
est une cloison.

\par Ainsi une facette de dimension $n$ est une chambre et une facette de
dimension $n-1$ est soit une cloison soit non sph\'erique. Donc l'ensemble \r F des
facettes de $\A$ et le c\^one de Tits d\'eterminent l'ensemble \r M des murs. 

\bigskip \noindent
{\bf 1.9. Quartiers} 

\par Un {\it quartier} dans $\A$ est un translat\'e $\g q=x+{{C^v}}$
d'une chambre vectorielle; $x$ est son {\it sommet} et ${{C^v}}$
sa {\it direction}. Deux quartiers ont la m\^eme direction si et seulement si ils
sont translat\'es l'un de l'autre et si et seulement si leur intersection contient
un autre quartier. 

\par Le {\it germe de quartier} du quartier $\g q=x+{{C^v}}$ est le
filtre $germ_\infty(\g q)=\g Q$ form\'e des sous-ensembles de $\A$ contenant un
translat\'e de
\g q; il est bien d\'etermin\'e par la direction ${{C^v}}$. Ainsi l'ensemble des
classes de quartiers de $\A$ \`a translation pr\`es, l'ensemble des chambres
vectorielles de $V$ et l'ensemble des germes de quartiers de $\A$ sont en
bijection canonique.

\par Une {\it face} (resp. {\it cloison}){\it de quartier} dans $\A$ est un
translat\'e $\g f=x+{{F^v}}$ d'une facette (resp. cloison) vectorielle. Le {\it
germe de face de quartier} de la face de quartier $\g f=x+{{F^v}}$ est le
filtre $germ_\infty(\g f)=\g F$ form\'e des sous-ensembles de $\A$ contenant un
{\it raccourci} de \g f \ie un translat\'e quelconque $\g f'$ de \g f par un
\'el\'ement de $\overline{F^v}$ (donc $\g f'\subset\g f$). Si $F^v$ est sph\'erique,
alors \g f et \g F sont  aussi dits {\it sph\'eriques} (ou {\it \'evas\'es}). Le {\it
signe} de \g f ou \g F est le signe de leur {\it direction} $F^v$.

\bigskip \noindent
{\bf 1.10. Chemin\'ees} 

\par Une {\it chemin\'ee} est associ\'ee \`a une facette $F=F(x,F^v_0)$ (sa {\it base})
et une facette vectorielle $F^v$ (sa {\it direction}), c'est le filtre

\par\quad $\g
r(F,F^v)=cl(F+F^v)=cl(germ_x(x+F_0^v)+F^v)\supset cl(F)+{\overline{F^v}} =
{\overline{F}}+{\overline{F^v}}$ .

\par La chemin\'ee $\g r(F,F^v)$ est dite {\it \'evas\'ee} si $F^v$ est sph\'erique, son
{\it signe} est celui de $F^v$. Cette chemin\'ee est dite {\it solide} (resp.  {\it
pleine}) si la direction de son support  a un fixateur fini dans $W^v$ (resp. si
son support est \A).

\bigskip \noindent
{\bf Remarques 1.11.} 

\par 1) Si $F^v$ ou $F$ est sph\'erique (resp. une chambre), alors $\g r(F,F^v)$
est solide (resp. pleine). En particulier une chemin\'ee \'evas\'ee est solide.

\par 2) La pl\'enitude de la chemin\'ee \g r \'equivaut au fait que tout $S\in\g r$
contient un ouvert non vide de $\A$ (lemme 1.7.5), alors le fixateur (point par
point)
$W_{\g r}$ de \g r dans $W$ est r\'eduit \`a l'\'el\'ement neutre.

\par 3) Dans le cas fini toute chemin\'ee est \'evas\'ee et toute face de quartier ou
facette est sph\'erique.

\par 4) Dans le cas fini ou si $F^v$ n'est pas sph\'erique, il peut arriver que la
chemin\'ee soit \`a la fois positive et n\'egative.

\par 5) Si $F^v=V_0$ est la facette vectorielle minimale, alors $\g
r(F,F^v)={\overline{F}}+V_0\subset cl^{si}({\overline{F}})$ et est \'egal \`a la facette
ferm\'ee ${\overline{F}}$ si toute $\qa\in\qD$ est nulle sur $V_0$ (\eg dans le cas
mod\'er\'ement imaginaire). Toute facette-ferm\'ee engendre donc une
chemin\'ee (d\'eg\'en\'er\'ee) qui n'est \'evas\'ee que dans le cas fini. 

\par 6) L'enclos d'une face de quartier $\g f=x+F^v$ est une chemin\'ee (avec
$F_0^v=V_0$ ou $F_0^v=F^v$) de m\^eme direction que \g f. La face de quartier est
sph\'erique si et seulement si la chemin\'ee correspondante est \'evas\'ee.

\par 7) Dans le cas mod\'er\'ement imaginaire et sauf si la direction est r\'eduite \`a
$\{0\}$, on peut d\'ecrire l'adh\'erence d'une face de quartier comme le $\R-$enclos
d'une demi-droite pr\'eordonn\'ee et une chemin\'ee comme l'enclos d'une demi-droite
pr\'eordonn\'ee et d'un germe de segment pr\'eordonn\'e de m\^eme origine. La face de
quartier est sph\'erique (ou la chemin\'ee est \'evas\'ee) si et seulement si la
demi-droite est g\'en\'erique.

\bigskip \noindent
{\bf 1.12. Germes de chemin\'ees ou demi-droites} 

\par Un {\it raccourci} d'une demi-droite $\qd$ d'origine $x$ est une demi-droite
$\qd_y$ d'origine $y\in\qd$ d\'efinie par $\qd_y=\qd+(y-x)=\qd\setminus[x,y[$. Le
{\it germe} de cette demi-droite est le filtre $germ_\infty(\qd)$ des parties de
$\A$ contenant un de ses raccourcis.

\par Un {\it raccourci} d'une chemin\'ee $\g r(F,F^v)$ est d\'efini par un \'el\'ement
$\qx\in{\overline{F^v}}$, c'est la chemin\'ee $cl(F+\qx+{F^v})$. Le
{\it germe} de cette chemin\'ee est le filtre $\g R(F,F^v)=germ_\infty(\g
r(F,F^v))$ des parties de
$\A$ contenant un de ses raccourcis.

\parni{\bf Remarques:} 1) La direction $F^v$ est bien d\'etermin\'ee par \g r ou \g R
(ce n'est pas vrai pour la base): ${\overline{F^v}}$ est la r\'eunion des
demi-droites vectorielles telles que tout \'el\'ement du filtre \g R contienne une
demi-droite de cette direction.

\par 2) Un raccourci d'une chemin\'ee a la m\^eme direction, le m\^eme support, le
m\^eme signe et donc les m\^emes propri\'et\'es (\'evas\'ee, solide ou pleine) que la chemin\'ee
originelle; ces objets et ces propri\'et\'es sont donc attach\'es au germe.

\par 3) Dans le cas d\'eg\'en\'er\'e o\`u $F^v=V_0$, tout raccourci de \g r est \'egal \`a \g
r, donc $\g r=\g R$.

\par 4) Un germe de quartier est un germe de chemin\'ee \'evas\'ee pleine.

\bigskip \noindent
{\bf 1.13. Appartements de type $\A$} 

\par Un {\it appartement de type} $\A$ est un ensemble $A$ muni d'un ensemble
$Isom(\A,A)$ de bijections $f:\A\rightarrow A$ (appel\'ees {\it isomorphismes}) tel
que, si $f_0\in Isom(\A,A)$, alors $f\in Isom(\A,A)$ si et seulement si il existe
$w\in W$  v\'erifiant $f=f_0\circ w$.

\par Un {\it isomorphisme} entre deux appartements $A$ et $A'$ est une bijection
$\qf:A\rightarrow A'$ telle qu'il existe $f_0\in Isom(\A,A))$ avec $\qf\circ
f_0\in Isom(\A,A')$ (\ie $f\in Isom(\A,A))\Leftrightarrow\qf\circ f\in
Isom(\A,A')$ ).

\par En particulier un appartement $A$ poss\`ede une structure naturelle d'espace
affine sous un espace vectoriel r\'eel $V(A)$, avec des ensembles 
$\qF(A)\subset\qD(A)\subset V(A)^*$ de racines, des murs ou murs imaginaires, des
c\^ones de Tits
$\pm \r T(A)$ dans $V(A)$ et des relations de pr\'eordre $\le _A$, ${\overline\le }_A$,
${\bc\le }_A$. Il contient des facettes, des demi-appartements, des quartiers, des
chemin\'ees ... et plus g\'en\'eralement tout ce que l'on vient de d\'efinir, puisque ces
notions sont invariantes par $W$. De plus un isomorphisme d'appartements \'echange
ces structures, objets et relations.

\bigskip
\parni {\S$\,${\bf 2.}$\quad${\bftwelve Masures affines }}
\bigskip

\par Pour simplifier on se place dor\'enavant dans le cas mod\'er\'ement imaginaire. Il
faudrait sinon rajouter \`a la d\'efinition 2.1 ci-dessous au moins certaines des
propri\'et\'es d\'emontr\'ees en 2.2. Le cas tr\`es imaginaire avec enclos r\'eel (\'egal \`a
l'enveloppe convexe) pourrait quand m\^eme \^etre int\'eressant.

\bigskip \noindent
{\bf D\'efinition 2.1.} Une {\it masure affine} de type $\A$ est un ensemble \r I
muni d'un recouvrement par un ensemble \r A de sous-ensembles appel\'es {\it
appartements} tel que:

\par (MA1) Tout $A\in \r A$ est muni d'une structure d'appartement de type $\A$.

\par (MA2) Si $F$ est un point, un germe d'intervalle pr\'eordonn\'e, une demi-droite
g\'en\'erique ou une chemin\'ee solide d'un appartement $A$ et si $A'$ est un autre
appartement contenant $F$, alors $A\cap A'$ contient l'enclos $cl_A(F)$ de $F$
dans $A$ et il existe un isomorphisme de $A$ sur $A'$ fixant cet enclos.

\par (MA3) Si \g R est un germe de chemin\'ee \'evas\'ee, si $F$  est une facette  ou
un germe de chemin\'ee solide, alors \g R et $F$ sont toujours contenus dans un m\^eme
appartement.

\par (MA4) Si deux appartements $A$, $A'$ contiennent \g R et $F$ comme en (MA3),
alors leur intersection $A\cap A'$ contient l'enclos $cl_A(\g R\cup F)$ de $\g
R\cup F$ dans $A$ et il existe un isomorphisme de $A$ sur $A'$ fixant cet enclos.

\medskip\par Par r\'ef\'erence aux d\'ecompositions correspondantes dans les groupes de
Kac-Moody (\cf \S\kern 2pt 6), si \g R  est un germe de quartier et $F$ une facette (resp.
un germe de quartier de m\^eme signe, de signe oppos\'e) on dit que (MA3) ou (MA4)
est la partie existence ou unicit\'e de la {\it propri\'et\'e d'Iwasawa} (resp. {\it
de Bruhat, de Birkhoff}); on ajoute le qualificatif {\it mixte} si on remplace
germe de quartier par germe de chemin\'ee \'evas\'ee non d\'eg\'en\'er\'ee ou solide non
d\'eg\'en\'er\'ee. La {\it propri\'et\'e de Bruhat-Iwahori} correspondrait au cas o\`u \g R et
$F$ sont des facettes (\cf (I1) et (I2) en 2.2.6 ci-dessous); elle est rarement
v\'erifi\'ee par les masures affines: seulement par les immeubles affines (\cf 2.4,
2.5).

\parni{\bf Remarque.} Les d\'efinitions des masures affines de type $\A$, $\A_{\sR}$,
$\A^{si}$ et $\A_{\sR}^{si}$ ne diff\`erent donc que par les contraintes plus ou
moins grandes impos\'ees par le fait de contenir les enclos correspondants. Le
type $\A^{si}$ est le plus contraignant, on esp\`ere pouvoir toujours s'y ramener. Le
type $\A_{\sR}$ est le moins contraignant.

\bigskip \noindent
{\bf 2.2. Premi\`eres cons\'equences}

\par 1) Si la conclusion de (MA2) est v\'erifi\'ee pour un filtre $F$, elle l'est
aussi pour tout filtre $F'$ tel que $F\subset F'\subset cl_A(F)$. En particulier
si (MA2) est v\'erifi\'e, il s'applique aussi bien quand $F$ est un germe de
segment pr\'eordonn\'e, une facette-locale, une facette, une facette-ferm\'ee ou une
face de quartier sph\'erique. De m\^eme (MA2) s'applique alors aussi quand $F$ est un
germe de demi-droite g\'en\'erique, un germe de face de quartier sph\'erique ou un
germe de chemin\'ee solide.

\par 2) L'axiome (MA2) et la remarque finale de 1.13 montrent que les notions de
germe de segment (ou d'intervalle) pr\'eordonn\'e, facette, facette-locale,
facette-ferm\'ee, demi-droite g\'en\'erique, face de quartier sph\'erique, chemin\'ee
solide et les germes associ\'es sont bien d\'efinis dans la masure ind\'ependamment de
l'appartement qui les contient. De m\^eme les qualificatifs qui s'appliquent \`a ces
objets (positif, n\'egatif, sph\'erique, g\'en\'erique, \'evas\'e, pleine, dimension, type,
chambre, cloison) ne d\'ependent pas de l'appartement. On utilise ceci dans les
axiomes (MA3) et (MA4).

\par 3) D'apr\`es (MA2) les axiomes (MA3) et (MA4) s'appliquent aussi dans une
masure quand $F$ est un point, un germe de segment pr\'eordonn\'e, une facette-locale
ou une facette-ferm\'ee et quand \g R ou $F$ est un germe de demi-droite g\'en\'erique
ou un germe de face de quartier sph\'erique.

\par Bien sur (MA3) (mais pas forc\'ement (MA4)) s'applique aussi si \g R ou $F$
est contenu dans l'un des objets indiqu\'es, par exemple si \g R est un
germe de demi-droite g\'en\'erique et si $F$ est un germe de demi-droite pr\'eordonn\'ee.

\par 4) Les enclos $cl_A(-)$ intervenant dans les axiomes (MA2) et (MA4) ne
d\'ependent finalement pas de l'appartement $A$ choisi dans la masure et sont donc
not\'es simplement $cl(-)$.

\par 5) Un mur (dans un appartement) est l'enclos de deux germes de cloisons de 
quartier; ainsi les appartements contenant un mur sont isomorphes et la notion de
mur est intrins\`eque dans une masure. De m\^eme pour un demi-appartement 
qui est l'enclos d'un germe de quartier et d'un germe
de cloison de quartier.

\par 6) Dans le cas fini, tout appartement $A$ est muni d'une m\'etrique
euclidienne  $d_A$, $W_A-$invariante. De plus on a vu (1.11.5) qu'une
facette-ferm\'ee est un germe de chemin\'ee \'evas\'ee, une masure affine v\'erifie donc
dans ce cas les trois axiomes habituels des immeubles (\cf [Tits-74] ou
[Brown-89] dans le cas discret et [Rousseau-08] dans le cas affine g\'en\'eral):

\par (I0) Tout $A\in\r A$ est un appartement euclidien.

\par (I1) Deux facettes sont toujours contenues dans un m\^eme appartement.

\par (I2) Si $A$ et $A'$ sont deux appartements, leur intersection est une
r\'eunion de facettes et, pour toutes facettes $F$, $F'$ dans $A\cap A'$, il existe
un isomorphisme de $A$ sur $A'$ fixant les facettes-ferm\'ees $\overline{F}$ et 
$\overline{F'}$.

\par En particulier la propri\'et\'e de Bruhat-Iwahori est v\'erifi\'ee.

\bigskip \noindent
{\bf D\'efinitions 2.3.} La masure affine \r I est dite {\it ordonn\'ee} si elle
v\'erifie l'axiome suppl\'ementaire suivant:

\par (MAO) Soient $x,y$ deux points de \r I et $A,A'$ deux appartements les
contenant; si $x\le _A y$, alors les segments $[x,y]_A$ et $[x,y]_{A'}$ d\'efinis par
$x$ et $y$ dans $A$ et $A'$ sont \'egaux.

\par La masure affine \r I est dite {\it \'epaisse} si toute cloison de \r I est
domin\'ee par au moins trois chambres.

\par Tous les qualificatifs qui s'appliquent \`a $\A$ s'appliquent aussi \`a \r I,
\eg semi-discret, discret, dense, sans imaginaire, dimension, rang (vectoriel), type
fini, type classique, ...

\parni{\bf Remarques:} 1) La seconde d\'efinition est habituelle. La premi\`ere est
une condition de "convexit\'e pr\'eordonn\'ee" des intersections d'appartements. On
verra au \S\kern 2pt 5 qu'elle permet de d\'efinir des pr\'eordres globaux sur \r I.

\par 2) Evidemment une masure \r I de type $\A$ est munie d'une structure
naturelle $\r I_{\sR}$ de masure de type $\A_{\sR}$ (\cf 1.6.1), il suffit de
rajouter tous les murs fant\^omes \`a la structure. Mais $\r I_{\sR}$ ne peut \^etre
\'epaisse (sauf si $\r M=\r M_{\sR}$). Il n'est pas sur que cette masure de type $\A$
puisse \^etre munie (comme on l'esp\`ere) d'une structure de masure de type $\A^{si}$
ou $\A^{si}_{\sR}$ (sans imaginaire).

\bigskip \noindent
{\bf 2.4. Exemples} 

\par 1) Soit \r I l'immeuble de Bruhat-Tits d'une donn\'ee radicielle valu\'ee (au
sens de [Bruhat-Tits-72]) que l'on munit de son syst\`eme canonique \r A
d'appartements. L'appartement t\'emoin $\A$ est un appartement affine relevant du
cas classique (donc fini) de 1.2.1 (sans imaginaire). D'apr\`es [Bruhat-Tits-72;
6.2.11] on peut supposer que la donn\'ee radicielle est g\'en\'eratrice et que l'image du
groupe
$N$ de cette donn\'ee dans le groupe affine de $\A$ est notre groupe $W$. De plus pour
$\emptyset\ne \qO\subset\A$, le groupe $U_\qO$ d\'efini dans \LC fixe le filtre
$cl(\qO)$, m\^eme avec notre d\'efinition plus restrictive d'enclos. Ainsi [\lc;
7.4.8 et 7.4.9] montre que l'intersection $A\cap A'$ de deux appartements
$A$, $A'$ est close (m\^eme pour notre d\'efinition plus restrictive) et que
$A$ et $A'$ sont isomorphes par un isomorphisme (au sens de 1.13) fixant cette
intersection. On en d\'eduit les axiomes (MA1), (MA2), (MA4) et (MAO). Le th\'eor\`eme
7.4.18 de \LC montre l'axiome (MA3) dans les cas impliquant des facettes ou des
germes de quartier. Cependant les chemin\'ees n'interviennent pas dans \LC Les
r\'esultats manquants pour l'axiome (MA3) sont d\'emontr\'es dans [Rousseau-77], voir
aussi [Rousseau-01]. 

\par Un immeuble affine irr\'eductible, de dimension $\ge 3$ et muni d'un syst\`eme
d'appartements v\'erifiant les axiomes de [Tits-86], est associ\'e \`a une donn\'ee
radicielle valu\'ee, mais en un sens qui peut \^etre plus g\'en\'eral que celui de
[Bruhat-Tits-72]. On n'est donc sur d'obtenir une masure affine que dans un
nombre plus restreint de cas. Mais la d\'emonstration de [Rousseau-77] devrait
pouvoir \^etre adapt\'ee \`a ce cadre, comme elle l'a \'et\'e ici pour les groupes de
Kac-Moody, \cf 6.7.

\par Un immeuble affine simplicial (\ie discret) muni de son syst\`eme complet
d'apparte-ments est une masure affine ordonn\'ee. L'essentiel des r\'esultats se
trouve dans les r\'ef\'erences classiques: [Brown-89], [Garrett-97], [Ronan-89],
[Scharlau-95] ou [Weiss-08]. Les r\'esultats manquants sur les chemin\'ees sont dans
[Charignon].

\par Les immeubles affines de [Tits-86] ont \'et\'e revisit\'es par Anne Parreau
[2000], voir aussi [Kleiner-Leeb-97]. Si on les munit de leur syst\`eme complet
d'appartements, il est vraisemblable que l'on puisse d\'emontrer que deux
demi-droites sont (quitte \`a les raccourcir) contenues dans un m\^eme appartement.
Joint aux r\'esultats connus, ceci montrerait qu'un tel immeuble affine est une
masure affine ordonn\'ee.

\par On obtient ainsi de nombreux exemples de masures affines ordonn\'ees \'epaisses
de type classique (qui sont aussi des immeubles). Mais ceci ne
donne pas de nouveaux objets.

\par 2) Soit $G$ un groupe de Kac-Moody sur un corps $K$, au sens minimal
("alg\'ebrique") de [Tits-87]. Ce groupe est muni d'une donn\'ee radicielle
g\'en\'eratrice de type un syst\`eme de racines r\'eelles $\qF$ de Kac-Moody (\ie comme
en 1.2.1) [Rousseau-06; 1.4]. Si $K$ est muni d'une valuation r\'eelle $\qo$ non
triviale, la donn\'ee radicielle de $G$ est naturellement munie d'une valuation
[Rousseau-06; 2.2] et on peut d\'efinir un appartement affine $\A$ correspondant
(voir [Gaussent-Rousseau-08; 2.2] pour $\qo$ discr\`ete, mais le cas g\'en\'eral
est semblable); cet appartement est mod\'er\'ement imaginaire.

\par Si $G$ est sym\'etrisable, si la valuation $\qo$ de $K$ est discr\`ete et si le
corps r\'esiduel contient $\C$, on construit dans [Gaussent-Rousseau-08]
une masure (hovel) sur laquelle le groupe $G$ agit. Cette masure n'est en g\'en\'eral
pas un immeuble, faute de satisfaire l'axiome (I1), \ie la propri\'et\'e de
Bruhat-Iwahori. C'est en fait (au moins sous la condition (SC) de 6.1) une masure
affine ordonn\'ee \'epaisse, mais certaines des propri\'et\'es exig\'ees dans la d\'efinition
2.1 ne sont pas prouv\'ees dans \LC Les d\'emonstrations manquantes sont regroup\'ees \`a
la fin de cet article (\S\kern 2pt 6).

\par Malheureusement il n'y a pas encore de construction de masure affine
associ\'ee \`a une donn\'ee radicielle valu\'ee quelconque. On l'esp\`ere au moins dans le
cas pr\'esent des groupes de Kac-Moody pour des corps r\'eellement valu\'es quelconques
[Rousseau-?].

\bigskip \noindent
{\bf 2.5. Commentaires sur les d\'efinitions} 

\par La d\'efinition de masure affine adopt\'ee ici est largement inspir\'ee de celle
d'immeuble de type affine donn\'ee par J. Tits [Tits-86]. Toutes les deux d\'ependent
d'un syst\`eme d'appartements choisi et font intervenir des \'el\'ements \`a l'infini: les
quartiers. On a en fait ici rajout\'e l'exigence de propri\'et\'es d'autres \'el\'ements \`a
l'infini, les chemin\'ees d\'efinies dans [Rousseau-77]. Cela permet de se passer de
l'axiome (A5) introduit par J. Tits (ou de ses variantes, \cf [Ronan-89;
appendice 2]). Par ailleurs on raisonne ici avec $\A$ plut\^ot qu'avec $\A_{\sR}$
comme dans [Tits-86].

\par Les axiomes pr\'esent\'es ici sont compliqu\'es par le fait que beaucoup de cas
diff\'erents sont consid\'er\'es. Je n'ai pas trouv\'e de formulation plus simple qui
n'exclurait pas l'exemple fondamental de Kac-Moody ci-dessus (2.4.2) et donc
ne risquerait pas de ne s'appliquer qu'\`a l'ensemble vide (en dehors du cas fini).
Il fallait en particulier \'eviter d'exiger l'axiome (I1) (ou propri\'et\'e de
Bruhat-Iwahori). Cet axiome me semblant fondamental pour les immeubles, j'ai
adopt\'e ici le nom de masure.

\par J'ai simplement demand\'e dans les axiomes ce que je savais sur les masures de
[Gaussent-Rousseau-08]. Les r\'esultats qui vont suivre semblent l\'egitimer ce choix
et il ne semble pas y avoir d'hypoth\`ese inutile. J'ai cependant choisi de s\'eparer
l'axiome (MAO): il ne fait intervenir que des \'el\'ements \`a distance finie et
pourrait d'ailleurs peut-\^etre s'av\'erer cons\'equence des autres axiomes. Les
propri\'et\'es sp\'ecifiques des masures affines ordonn\'ees sont abord\'ees dans le \S\kern 2pt 5.

\par La d\'efinition de masure affine donn\'ee ici est pr\'ecise. La notion g\'en\'erale de
masure ("immeuble sans l'axiome I1") est beaucoup plus vague. Il existe cependant
au moins un autre exemple int\'eressant: [Broussous-04].

\par Les masures, comme les immeubles, sont des exemples d'{\it espaces
couverts}, c'est \`a dire d'ensembles munis d'un recouvrement par des
sous-ensembles (appartements ou plats) tous isomorphes par "suffisamment
d'isomorphismes". C'est aussi le cas des espaces riemanniens sym\'etriques munis
des plats (maximaux); dans ce cas (de mani\`ere analogue \`a l'axiome (I1) des
immeubles) deux points sont toujours contenus dans un m\^eme plat. Des exemples ne
v\'erifiant pas cette propri\'et\'e, donc plus proches des masures, semblent \`a chercher
parmi les espaces pseudo-riemanniens sym\'etriques. 

\bigskip \noindent
{\bf 2.6. R\'etractions} 

\par Soit \g R un germe de chemin\'ee \'evas\'ee pleine dans un appartement $A$ d'une
masure affine \r I.

\par Pour $x\in\r I$, choisissons un appartement $A_x$ contenant \g R et $x$
(MA3). Il existe alors un isomorphisme $\qf:A_x\rightarrow A$ fixant \g R. Si
$\qf$ et $\qf'$ sont deux tels isomorphismes, $\qf^{-1}\circ\qf'$ est un
automorphisme de $A_x$ fixant la chemin\'ee pleine \g R, c'est donc l'identit\'e
(1.11.2) et $\qf$ est unique. De m\^eme d'apr\`es (MA4) $\qf(x)$ ne d\'epend pas du
choix de l'appartement $A_x$.\qquad\qquad On peut donc d\'efinir \qquad
$\qr_{A,\sg R}(x)=\qf(x)$\qquad.

\par L'application $\qr_{A,\sg R}:\r I\rightarrow A$ est une r\'etraction de \r I
sur $A$. Elle ne d\'epend que de $A$ et \g R; on l'appelle la {\it r\'etraction de \r
I sur $A$ de centre} \g R.

\par Si $A'$ est un autre appartement contenant \g R, on a \'evidemment $\qr_{A',\sg
R}=\psi\circ\qr_{A,\sg R}$ o\`u $\psi$ est l'unique isomorphisme de $A$ sur $A'$
fixant \g R. Ainsi, \`a isomorphisme unique pr\`es, $\qr_{A,\sg R}$ ne d\'epend que de
\g R.

\parni{\bf Remarque:} Dans le cas fini, on retrouve ainsi les r\'etractions bien
connues (d'immeubles affines) de centre des chambres ou des germes de
quartier. Mais on a aussi des r\'etractions de centre des germes de chemin\'ees
pleines. En fait celles-ci sont des chambres des diff\'erents immeubles affines
intervenant dans la compactification de Satake de l'immeuble (\cf \S\kern 2pt 4).

\bigskip
\parni {\bf Proposition 2.7.} {\sl  
\par 1) Soient \g R un germe de chemin\'ee \'evas\'ee et $x$, $y$ deux points dans un
appartement $A$ d'une masure affine v\'erifiant $x\le _Ay$. Il existe une subdivision
$z_0=x,z_1,\dots,z_n=y$ du segment $[x,y]_A$ et des appartements $A_1,\dots,A_n$
tels que, pour $1\le i\le n$, l'appartement $A_i$ contienne \g R et le segment
$[z_{i-1},z_i]_A$.

\par 2) Soient $F$ une facette ou un germe de chemin\'ee solide et \g F, \g F' deux
germes de faces de quartier sph\'eriques dans une masure affine \r I. On suppose
qu'il existe un appartement $A$ de \r I contenant \g F et \g F' et que, dans cet
appartement, \g F et \g F' ont m\^eme direction. Alors il existe dans $A$ des
germes de faces de quartier (sph\'eriques) de m\^eme direction $\g F_0=\g F,\g
F_1,\dots,\g F_n=\g F'$ et des appartements $A_1,\dots,A_n$ tels que, pour
$1\le i\le n$, l'appartement $A_i$  contienne $\g F_{i-1}$, $\g F_i$ et $F$.
}
\parni {\bf  N.B.} a) Bien sur, dans l'\'enonc\'e ci-dessus ont peut remplacer
facette par facette-locale, facette-ferm\'ee, point ou germe de segment pr\'eordonn\'e;
 germe de chemin\'ee solide ou \'evas\'ee par germe de face de quartier sph\'erique ou
demi-droite g\'en\'erique et face de quartier sph\'erique par demi-droite g\'en\'erique.

\par b) On sait de plus (MA4) que le fait d'avoir m\^eme direction ne d\'epend pas de
l'appartement contenant les deux faces de quartier sph\'eriques ou demi-droites
g\'en\'eriques.

\parni {\bf  D\'emonstration.} Les preuves des deux r\'esultats sont tr\`es proches. De
plus la premi\`ere preuve est essentiellement celle de [Gaussent-Rousseau-08;
sect. 4.4]. On ne fait donc que la seconde preuve.

\par Notons $\g f=x_-+F^v$ et $\g f'=x_++F^v$ des repr\'esentants de \g F et $\g
F'$ dans $A$. Le cas o\`u $F^v=V_0$ ne correspond \`a des faces de quartier
sph\'eriques que dans le cas fini et rel\`eve de la partie 1).
 On peut donc supposer que $F^v$ contient une demi-droite g\'en\'erique
$\qd$ d'origine $0$ dans $V$. 

\par Quitte \`a remplacer $x_+$ par $x_++\qx$ (et donc \g
f' par une raccourcie) pour $\qx$ grand dans $\qd$, on peut supposer $x_-\le _Ax_+$
ou $x_-\ge _Ax_+$. Pour $z\in[x,y]_A$, l'ensemble $\g r_z^\pm =cl_A([z,x_\pm )+F^v)$ (si
$z\ne x_\pm $) est une chemin\'ee \'evas\'ee. D'apr\`es (MA3) il existe un appartement $A^\pm _z$
contenant $F$ et un raccourci $\g r_z^\pm +\qx^\pm _z=cl_A([z,x_\pm )+\qx^\pm _z+F^v)$ pour
$\qx^\pm _z$ assez grand dans $\qd$; cet appartement contient en fait une bande 
$[z,z+\qe^\pm _z(x_\pm -z)]+\qx^\pm _z+F^v$ pour $\qe^\pm _z>0$ petit. Mais
$[z+\qe^-_z(x_--z),z+\qe^+_z(x_+-z)]$ est un voisinage de $z$ dans $[x_-,x_+]$ qui
est compact. On peut donc recouvrir $[x_-,x_+]$ par un nombre fini de tels
intervalles: il existe une subdivision $z_0=x_-,z_1,\dots,z_n=x_+$ de $[x_-,x_+]$,
des appartements $A_1,\dots,A_n$ et des \'el\'ements $\qx_1,\dots,\qx_n\in\qd$ tels
que, pour $1\le i\le n$, l'appartement $A_i$ contient $F$ et $[x_{i-1},x_i]+\qx_i+F^v$.
Il suffit alors de poser $\g f_i=x_i+\qx_i+F^v$ (on peut m\^eme supposer tous les
$\qx_i$ \'egaux).
\qed

\parni {\bf Corollaire 2.8.} {\sl  Dans la situation de 2.6 la r\'etraction
$\qr=\qr_{A,\sg R}$ est "croissante". Plus pr\'ecis\'ement, si $x$ et $y$ sont deux
points d'un appartement $A'$ tels que $x\le _{A'}y$, l'image par $\qr$ du segment
$[x,y]_{A'}$ est une ligne bris\'ee $[\qr z_0,\qr z_1]\cup[\qr z_1,\qr z_2]
\cup\dots\cup[\qr z_{n-1},\qr z_n]$ avec $\qr z_0=\qr x$, $\qr z_n=\qr y$, $\qr
z_{i-1}\le _A\qr z_i$ (pour $1\le i\le n$) et donc $\qr x\le _A\qr y$. Le m\^eme r\'esultat est
valable pour la relation $\bc\le $. }
\parni{\bf N.B.} Pour les masures affines pr\'eordonn\'ees on verra  au \S\kern 2pt 5 que les
diff\'erentes relations $\le _A$ ou $\bc\le _A$ se recollent en une relation de pr\'eordre
sur \r I; ceci donne un sens plus habituel au qualificatif "croissante".

\parni {\bf  D\'emonstration.} C'est une cons\'equence imm\'ediate de la proposition et
du fait que les isomorphismes d'appartements \'echangent les relations de pr\'eordre
(1.13).
\qed

\parni {\bf Proposition 2.9.} {\sl  1) Soit $D$ un demi-appartement de la masure
affine \r I et $M=\partial D$ son mur (\ie son bord). On consid\`ere une cloison
$F$ dans $M$ et une chambre $C$ de \r I dominant $F$. Il existe alors un
appartement $A$ de \r I contenant $D$ et $C$.
\par 2) Soient $D$ et $D_1$ deux demi-appartements de \r I dont l'intersection
est r\'eduite au mur bordant chacun d'eux: $M=D\cap D_1=\partial D=\partial D_1$.
Alors $D\cup D_1$ est un appartement de \r I.
 }
\parni {\bf  D\'emonstration.} 1) Soit $\g f=x+F^v$  une cloison de quartier
contenue dans $M$ de sommet $x\in\overline F$. Il existe un appartement $A_1$
contenant $\overline C$ et le germe \g F de \g f, donc aussi $cl(F\cup\g
F)\supset\overline F+F^v$. On note $\g r=\g r(C,F^v)$ et \g R son germe (dans
$A_1$). Il est clair que $cl(F\cup\g R)\supset\g r\supset\overline C$. Soit $\g
f'=x-F^v$ la cloison de quartier dans $M$ de sommet $x$ et de direction oppos\'ee \`a
\g f. On note \g q' le quartier de $D$ de sommet $x$ tel que $\g
f'\subset\overline{\g q'}$. Il existe un appartement $A_2$ contenant son germe \g
Q' et $\g R\supset\g F$. Ainsi
$D=cl(\g F\cup\g Q')$ est dans $A_2$ qui contient aussi $cl(F\cup\g
R)\supset\overline C$.

\par 2) Soit $x\in M$, on choisit \g f, \g f' ($\subset M$) et \g q' ($\subset
D$) comme ci-dessus. 

\par Soit $\g f_1\subset D_1\setminus M$ une cloison de quartier
de m\^eme direction que \g f et $A_1$ un appartement contenant les germes \g Q' et
$\g F_1$; montrons que $A_1$ contient \g F (et donc $D=cl(\g Q'\cup\g F)$). On se
ram\`ene par r\'ecurrence, gr\^ace \`a 2.7.2, au cas o\`u \g Q', \g F et $\g F_1$ sont dans
un m\^eme appartement $A_2$. Si, dans $A_2$, $\g F\subset cl(\g Q'\cup\g F_1)$ le
r\'esultat est \'evident. Sinon $\g F_1\subset cl(\g Q'\cup\g F)=D$, c'est absurde
d'apr\`es l'hypoth\`ese $D\cap D_1=M$.

\par Soit maintenant $\g q_1$ le quartier de $D_1$ de sommet $x$ tel que $\g
f\subset\overline{\g q_1}$. Il existe un appartement $A_3$ contenant les germes
$\g Q'$ et $\g Q_1$; il contient donc une cloison de quartier $\g f_1\subset
D_1\setminus M$ de m\^eme direction que \g f. D'apr\`es l'alin\'ea pr\'ec\'edent $A_3$
contient \g F et m\^eme $D$. Ainsi $A_3$ contient $D_1=cl(\g F'\cup\g Q_1)$.
Finalement $D\cup D_1\subset A_3$ et on a \'egalit\'e car, dans un appartement, un mur
 ne borde que deux demi-appartements.
\qed

\parni {\bf Corollaire 2.10.} {\sl  Si la masure \r I est \'epaisse, un
demi-appartement est l'intersection de deux appartements.}

\parni{\bf Cons\'equence.} Ainsi dans une masure \'epaisse un sous-ensemble clos
d'appartement est une intersection d'appartements. Par contre on n'est pas
sur que l'intersection de deux appartements soit close.

\parni{\bf Remarque.} On dit qu'un mur est {\it \'epais} si tout demi-appartement
qu'il borde est intersection de deux appartements. On va montrer qu'un mur est
\'epais si (et seulement si) une (resp. toute) cloison qu'il contient est \'epaisse
\ie domin\'ee par au moins trois chambres.

\parni {\bf  D\'emonstration.} Soit $D$ un demi-appartement dans un appartement
$A_1$ de \r I. On choisit une cloison $F=F(x,F^v)$ (avec $F^v$ cloison
vectorielle) dans le mur $M$ bordant $D$. On note $C$ (resp. $C_1$) la chambre de
$A_1$ dominant $F$ contenue dans $D$ (resp. $A_1\setminus D$). Par hypoth\`ese il
existe une chambre $C_2$ dominant $F$ diff\'erente de $C$ et $C_1$. D'apr\`es la
proposition il existe un appartement $A_2$ contenant $D$ et $C_2$; montrons que
$D=A_1\cap A_2$. Sinon il existe $y$ dans $(A_1\cap A_2)\setminus D$. Calcul\'e dans
$A_1$, l'enclos de $\{y\}\cup germ_\infty(x\pm F^v)$ contient une chemin\'ee raccourcie
de $\g r^\pm =\g r(C_1,\pm F^v)$. Ainsi $A_1\cap A_2$ contient les germes $\g R^+$,
$\g R^-$ et donc leur enclos: $A_1\cap A_2$ contient l'espace de $A_1$ compris
entre $M$ et un mur de $A_1\setminus D$ parall\`ele \`a $M$. On en d\'eduit aussit\^ot que
$A_1\cap A_2$ contient $C_1$, ainsi $A_2$ contient trois chambres distinctes $C$,
$C_1$ et $C_2$ dominant $F$, c'est absurde.
\qed

\bigskip
\parni {\S$\,${\bf 3.}$\quad${\bftwelve Parall\'elisme et immeubles jumel\'es \`a
l'infini }}
\bigskip
\par On consid\`ere une masure affine \r I

\bigskip \noindent
{\bf D\'efinitions 3.1.} 

\par Soient \g f et \g f' deux faces de quartier sph\'eriques, il existe un
appartement $A$ contenant les germes correspondants \g F et \g F'.

\par On dit que \g f est {\it parall\`ele} \`a \g f' (not\'e $\g f\parallel\g f'$) si \g
F et \g F' ont m\^eme direction dans $A$, \ie on peut \'ecrire 
$\g F=germ_\infty(x+F^v)$ et $\g F'=germ_\infty(y+F^v)$ dans $A$.

\par On dit que \g f et \g f' sont {\it de directions oppos\'ees} (ou {\it 
oppos\'ees} pour abr\'eger) (not\'e $\g f \;opp\;\g f'$) si on peut \'ecrire 
$\g F=germ_\infty(x+F^v)$ et $\g F'=germ_\infty(y-F^v)$ dans $A$.

\par On dit que \g f {\it domine} \g f' (not\'e $\g f\ge \g f'$) si on peut \'ecrire 
$\g F=germ_\infty(x+F^v)$ et $\g F'=germ_\infty(y+F'^v)$ dans $A$ avec
$F'^v\subset\overline{F^v}$.

\medskip
\parni{\bf Remarques.} a) Comme $A$ est unique \`a isomorphisme fixant \g F
et
\g F' pr\`es, ces d\'efinitions ne d\'ependent pas du choix de $A$ et bien sur elles ne
d\'ependent que des germes.

\par b) Si \g f est parall\`ele \`a \g f' ou domine \g f' (resp. est oppos\'ee \`a \g f'),
les deux faces de quartier ont m\^eme signe (resp. des signes oppos\'es).

\par c) Pour les quartiers, \g q est parall\`ele \`a \g q' si et seulement si $\g Q=\g
Q'$.

\par d) Ces d\'efinitions et la proposition suivante sont aussi valables pour les
demi-droites g\'en\'eriques. On pourrait aussi parler de parall\'elisme et opposition
entre un germe de segment pr\'eordonn\'e et un germe de demi-droite g\'en\'erique, \cf
[Rousseau-01].

\bigskip \noindent
{\bf Proposition 3.2.} {\sl  

\par 1) Le parall\'elisme est une relation d'\'equivalence sur
les faces de quartier sph\'eriques.

\par 2) La domination est un pr\'eordre (partiel) sur les faces de quartier
sph\'eriques qui est compatible avec le parall\'elisme. Plus pr\'ecis\'ement le
parall\'elisme est la relation d'\'equivalence associ\'ee \`a ce pr\'eordre. 

\par On a :\qquad
$\g f\ge \g f'\quad\Leftrightarrow\quad\exists\;\g f'_1\;,\;\g f'_1\parallel\g
f'\;,\;\g f'_1\subset \overline{\g f}\quad\Leftrightarrow\quad\exists\;\g
f_1\;,\;\g f_1\parallel\g f\;,\;\g f'\subset \overline{\g f_1}$

\par 3) L'opposition est sym\'etrique et compatible avec le parall\'elisme et
la domination:

\par on a:\qquad $\g f_1\parallel\g f\;opp\;\g f'\parallel\g f'_1\Rightarrow\g
f_1\;opp\;\g f'_1$\qquad ,
\qquad $\g f\;opp\;\g f'\ge \g f''\Leftrightarrow\exists\; \g f'_1,\; \g f\ge \g
f'_1\; opp\;\g f''$
\parni et \qquad\qquad $\g f\;opp\;\g f'\le \g f''\Leftrightarrow\exists\; \g
f'_1,\; \g f\le \g f'_1\; opp\;\g f''$\qquad.
}

\parni {\bf  D\'emonstration.} 

\par Soient \g f, \g f' et $\g f''$ trois faces de quartier sph\'eriques telles que
$\g f'\parallel\g f''$ et $\g f\parallel\g f'$ (resp. $\g f\ge \g f'$, $\g f\le \g f'$,
$\g f\;opp\;\g f'$), il est \'evident que $\g f\parallel\g f''$ (resp. $\g f\ge \g
f''$, $\g f\le \g f''$, $\g f\;opp\;\g f''$) si les trois germes sont situ\'es dans un
m\^eme appartement. Mais la proposition 2.7 permet de se ramener \`a ce cas par
r\'ecurrence; on a donc la conclusion en g\'en\'eral.
\par On vient donc de prouver 1) et la compatibilit\'e de la domination ou de
l'opposition avec le parall\'elisme. La sym\'etrie de l'opposition est claire.

\par Pour la derni\`ere assertion de 2), on peut donc remplacer \g f (ou \g f') par
une face de quartier parall\`ele $\g f_1$ (ou $\g f'_1$) de m\^eme origine et dans un
m\^eme appartement que \g f' (ou \g f) et alors les \'equivalences annonc\'ees sont
claires. La preuve de la transitivit\'e de la domination se ram\`ene de m\^eme au cas
facile de trois faces de quartier de m\^eme sommet. On a \'egalement $\g f\parallel\g
f'\Leftrightarrow\g f\ge \g f'$ et $\g f'\ge \g f$; d'o\`u la preuve compl\`ete de 2).

\par Pour les deux derni\`eres assertions de 3), si $\g f\;opp\;\g f'\ge \g f''$ (resp.
$\g f\;opp\;\g f''\le \g f'$), on peut supposer \g f et \g f' dans un m\^eme
appartement $A$ puis, quitte \`a raccourcir \g f', remplacer $\g f''$ par une face 
de quartier parall\`ele, de m\^eme sommet et dans un m\^eme appartement que \g f'. On a
alors $\g f''\subset\overline{\g f'}$  et les trois faces de quartier sont dans 
un m\^eme appartement. Les r\'esultats annonc\'es sont alors \'evidents.
\qed

\parni
{\bf 3.3. Cons\'equences}

\par 1) La classe d'\'equivalence d'une face de quartier sph\'erique \g f de germe \g
F est appel\'ee sa {\it direction} $\g f^\infty=\g F^\infty$. La {\it direction}
d'une chemin\'ee (ou d'un germe de chemin\'ee) \'evas\'ee est la direction d'un germe de
face de quartier de dimension maximale qu'il contient. Ces notions co\"{\i}ncident avec
les notions de direction dans un appartement

\par Les directions de faces de quartier
sph\'eriques, aussi appel\'ees {\it facettes \`a l'infini}, forment un ensemble $\r
I^\infty$ r\'eunion de deux sous-ensembles: $\r I^{+\infty}$ (resp. $\r
I^{-\infty}$) est l'ensemble des directions de faces de quartier sph\'eriques
positives (resp. n\'egatives). Ces deux sous-ensembles sont confondus dans le cas
fini et disjoints dans les autres cas.

\par On notera que la masure $\r I_{\sR}$ (\cf 2.3.2) fournit les m\^emes ensembles
de facettes \`a l'infini: $\r I_{\sR}^\infty=\r I^\infty$ et $\r
I_{\sR}^{\pm \infty}=\r I^{\pm \infty}$.

\par 2) La domination induit un ordre (partiel) sur $\r I^\infty$, qui ne m\'elange
pas $\r I^{+\infty}$ et $\r I^{-\infty}$ (s'ils sont disjoints). La face \g f est
un quartier (resp. une cloison de quartier) si et seulement si $\g f^\infty$ est
maximale (resp. maximale parmi les facettes \`a l'infini non maximales); on dit
alors que $\g f^\infty$ est une {\it chambre \`a l'infini} (resp. {\it cloison \`a
l'infini}).

\par 3) \`A un appartement $A$ on associe {\it l'appartement \`a l'infini}
$A^\infty=A^{+\infty}\cup A^{-\infty}$, o\`u $A^{+\infty}$ (resp. $A^{-\infty}$) 
est l'ensemble des classes de faces de quartier sph\'eriques positives (resp.
n\'egatives) dont un repr\'esentant est dans $A$. Par d\'efinition des relations,
$A^{+\infty}$ (resp. $A^{-\infty}$) est isomorphe \`a l'ensemble ordonn\'e des
facettes vectorielles sph\'eriques positives (resp. n\'egatives) de $V(A)$ \ie de $\r
T^\circ(A)$ (resp. $-\r T^\circ(A)$). Par ces isomorphismes l'opposition
correspond \`a $F^v\mapsto-F^v$ et induit donc une bijection croissante de
$A^{+\infty}$ sur $A^{-\infty}$.

\par Il n'est pas clair (pour l'instant, \cf 3.9) que l'application $A\mapsto
A^{+\infty}$ ou $A\mapsto A^{-\infty}$ soit injective; par contre $A\mapsto
A^{\infty}$ l'est (par convexit\'e).

\par 4) Si $M$ est un mur (resp. $D$ est un demi-appartement) de l'appartement
$A$ de \r I, l'ensemble $M^\infty$ (resp. $D^\infty$) des classes de faces de
quartier sph\'eriques dont un repr\'esentant est dans $M$ (resp. $D$) est, par
d\'efinition, un {\it mur} (resp. {\it demi-appartement}) {\it \`a l'infini} (dans
$A^\infty$) appel\'e {\it direction} de $M$ (resp. $D$). De m\^eme 
$M^{\pm \infty}=M^\infty\cap\r I^{\pm \infty}$ (resp.
$D^{\pm \infty}=D^\infty\cap\r I^{\pm \infty}$) est un {\it mur} (resp. un {\it
demi-appartement}) dans
$A^{\pm \infty}$. Deux murs parall\`eles (resp. deux demi-appartements inclus l'un
dans l'autre) donnent le m\^eme mur (resp. demi-appartement) \`a l'infini. On peut
parler de $M^\infty$ et
$D^\infty$ m\^eme pour $M$ et $D$ fant\^omes (1.6.1).

\par Si $M=M(\qa,k)$ (resp. $D=D(\qa,k)$), $M^{\pm \infty}$ (resp. $D^{\pm \infty}$)
s'identifie au "mur" $\pm \r T^\circ(A)\cap Ker(\qa)$ (resp. au "demi-appartement"
$\pm \r T^\circ(A)\cap D(\qa)$) dans $\pm \r T^\circ(A)\subset V(A)$.

\bigskip
\parni {\bf Th\'eor\`eme 3.4.} {\sl  L'ensemble $\r I^{+\infty}$ muni de sa relation
d'ordre et de l'ensemble $\r A^{+\infty}$ de ses appartements \`a l'infini
$A^{+\infty}$  (pour $A\in\r A$) est un immeuble combinatoire de type $W^v$. Si
la masure \r I est \'epaisse, alors l'immeuble $\r I^{+\infty}$ est \'epais. }

\parni{\bf Remarques.} a) On a le m\^eme r\'esultat avec $-\infty$.

\par b) A proprement parler $\r I^{+\infty}$ n'est pas un complexe simplicial
puisque dans l'adh\'erence d'un quartier  on ne garde que les faces sph\'eriques. On
a cependant bien toutes les chambres et cloisons du complexe. Il para\^{\i}t donc
opportun d'adopter (ci-dessous) le langage des syst\`emes de chambres.

\par c) D'apr\`es la derni\`ere assertion de 3.3.1 et la remarque 2.3.2, l'immeuble
$\r I^{+\infty}$ peut \^etre \'epais sans que \r I le soit. De plus, comme il y a une
infinit\'e de murs parall\`eles, la d\'emonstration ci-dessous prouve que, quand \r I
est \'epaisse, une cloison \`a l'infini est domin\'ee par une infinit\'e de chambres \`a
l'infini.

\par d) D'apr\`es 2.7.2 une r\'etraction $\qr$ de centre un germe de chemin\'ee \'evas\'ee
pleine positive \g R est compatible avec le parall\'elisme et est donc d\'efinie sur
$\r I^{\pm \infty}$. Si \g R est un germe de quartier, cette r\'etraction co\"{\i}ncide sur
$\r I^{+\infty}$ avec celle d\'eduite de sa structure d'immeuble et sur $\r
I^{-\infty}$ avec celle d\'eduite de la structure d'immeuble jumel\'e sur $\r
I^{\infty}$, \cf 3.7.

\par e) L'ensemble des classes de parall\'elisme de demi-droites positives g\'en\'eriques
est une r\'ealisation g\'eom\'etrique de $\r I^{+\infty}$: \`a $(x+F^v)^\infty$ on associe
l'ensemble des directions des demi-droites $x+\R_+\qx$ pour $\qx\in F^v$.

\parni {\bf  D\'emonstration.}  
\par V\'erifions les axiomes des immeubles comme \'enonc\'es en 2.2.6 ou, de mani\`ere
plus appropri\'ee aux immeubles combinatoires, dans [Brown-89; 4.1] ou [R\'emy-02;
2.4.1]. Chaque appartement $A^{+\infty}$ est le complexe de Coxeter associ\'e \`a
$W^v$, puisqu'il est isomorphe \`a l'ensemble ordonn\'e des facettes (sph\'eriques) de
$\r T^\circ(A)$, d'o\`u (I0). Les axiomes (I1) et (I2) sont des cons\'equences
imm\'ediates des axiomes (MA3) et (MA4).

\par Soit $\g f^\infty$ une cloison de $\r I^{+\infty}$ et $\g f=x+F^v$ une
cloison de quartier de \r I correspondante. On peut supposer que \g f a pour
support un (vrai) mur $M$ d'un appartement $A_1$. Si la masure \r I est \'epaisse,
on a construit en 2.10 un second appartement $A_2$ tel que $D=A_1\cap A_2$ est
l'un des deux demi-appartements de $A_1$ bord\'e par $M$. Notons \g q (resp $\g
q_1$, $\g q_2$) le quartier de $D$ (resp. $A_1\setminus D$, $A_2\setminus D$) de
sommet $x$ contenant \g f dans son adh\'erence; $\g f^\infty$ est bien domin\'ee par
les chambres \`a l'infini $\g q^\infty$, $\g q_1^\infty$ et $\g q_2^\infty$. Si $\g
f'=x-F^v$, on a $cl(\g F'\cup\g Q)=D$, $cl(\g F'\cup\g
Q_1)=\overline{A_1\setminus D}$ et $cl(\g F'\cup\g Q_2)=\overline{A_2\setminus
D}$; donc les trois chambres \`a l'infini $\g q^\infty$, $\g q_1^\infty$ et $\g
q_2^\infty$ sont bien distinctes.
\qed

\parni {\bf 3.5. $W^v-$distance}

\par $W^v$ est le groupe de Weyl vectoriel de $\A$; son syst\`eme de g\'en\'erateurs
canonique $S$ est associ\'e \`a la chambre fondamentale $C^v_f$ de $\A$. La
$W^v-$distance sur les chambres vectorielles de $\A$ est donc donn\'ee par
$d_+^{\sA}(wC^v_f,w'C^v_f)=w^{-1}w'$, elle est invariante par l'action de $W^v$.
On a $d_+^{\sA}(C^v_1,C^v_2)\in S\cup\{1\}$ si et seulement si $C_1^v$ et $C_2^v$
sont adjacentes \ie dominent une m\^eme cloison vectorielle.

\par Soient $\g Q_1$ et $\g Q_2$ deux germes de quartiers positifs et $A$ un
appartement les contenant; on note $\g Q_i=germ_\infty(x+C_i^v)$ dans $A$. Si
$f\in Isom(\A,A)$, on d\'efinit $d_+^{A}(\g Q_1,\g
Q_2)=d_+^{\sA}(f^{-1}(C^v_1),f^{-1}(C^v_2))$; un autre choix de $f$ est de la
forme $f\circ w$ pour $w\in W^v$, donc 
$d_+^{\sA}(w^{-1}f^{-1}(C^v_1),w^{-1}f^{-1}(C^v_2))
=d_+^{\sA}(f^{-1}(C^v_1),f^{-1}(C^v_2))$, ainsi  $d_+^{A}$ ne d\'epend pas du
choix de $f$. Par ailleurs deux choix de $A$ diff\`erent par un isomorphisme fixant
$\g Q_1$ et $\g Q_2$, donc $d_+^{A}(\g Q_1,\g Q_2)$ ne d\'epend pas de $A$; c'est
la {\it $W^v-$distance} $d_+(\g Q_1,\g Q_2)$ de $\g Q_1$ et $\g Q_2$ dans $\r
I^{+\infty}$.

\par Comme cette $W^v-$distance est associ\'ee \`a un immeuble, elle v\'erifie
les axiomes des $W^v-$distances, \cf [R\'emy-02; 2.3]. La v\'erification directe est
facile et laiss\'ee au lecteur: le premier axiome est clair; pour le troisi\`eme il
suffit de choisir les germes de quartier $x$, $y$ et $z$ dans le m\^eme appartement
; quant au second on peut s'inspirer de la d\'emonstration de 3.7 ci-dessous.

\par On a bien sur les m\^emes r\'esultats dans $\r I^{-\infty}$: pour des germes de
quartiers n\'egatifs et avec les notations ci-dessus on pose: $d_-(\g Q_1,\g
Q_2)=d_+^{\sA}(-f^{-1}(C_1^v),-f^{-1}(C_2^v))$.

\bigskip
\parni {\bf 3.6. Codistance}

\par Soient $\g Q_1$ et $\g Q_2$ deux germes de quartier de signes oppos\'es
($\qe_1$ et $\qe_2=-\qe_1$), $A$ un appartement les contenant et $f\in
Isom(\A,A)$; on note $\g Q_i=germ_\infty(x+C_i^v)$ dans $A$ et
$f^{-1}(C_i^v)=\qe_i w_iC^v_f$ pour $w_i\in W^v$.

\par On d\'efinit alors \qquad$d_*(\g Q_1,\g Q_2)=d_{\qe_1}(C_1^v,-C^v_2) =
d_{\qe_2}(-C_1^v,C^v_2)=w_1^{-1}w_2$.

\parni Pour les m\^emes raisons qu'en 3.5, cette d\'efinition ne d\'epend pas des choix
effectu\'es; on dit que $d_*$ est la {\it codistance} dans $\r I^\infty$.

\par Deux germes de quartiers (de signes diff\'erents) sont oppos\'es dans \r I si et
seulement si leur codistance est $1$. Deux germes de faces de quartier sph\'eriques
$\g F$ et $\g F'$ sont oppos\'es si et seulement si, pour tout germe de quartier \g
Q dominant \g F, il existe un germe de quartier $\g Q'$ oppos\'e \`a \g Q dominant
$\g F'$ et inversement en \'echangeant \g F et $\g F'$. On peut donc traduire
l'opposition par des codistances.
\medskip
\parni{\bf Cas fini:} Supposons $W^v$ fini et notons $w_0$ son \'el\'ement de plus
grande longueur. Alors $\r I^{+\infty}=\r I^{-\infty}=\r I^{\infty}$ est un
immeuble sph\'erique et on retrouve donc bien le classique immeuble sph\'erique \`a
l'infini de l'immeuble affine \r I.

\par On a d\'efini trois "distances" sur l'ensemble des germes de quartiers. Pour
tous germes de quartiers $\g Q_1$, $\g Q_2$ on a: $d_-(\g Q_1,\g Q_2)=w_0d_+(\g
Q_1,\g Q_2)w_0$ et $d_*(\g Q_1,\g Q_2)=w_0d_+(\g Q_1,\g Q_2)$ (resp. $d_+(\g
Q_1,\g Q_2)w_0$) si on consid\`ere $\g Q_1$ dans $\r I^{-\infty}$ et $\g Q_2$ dans
$\r I^{+\infty}$ (resp. l'inverse). On retrouve les formules de [Tits-92; prop.1]
 \`a condition d'\'echanger $w_0d_+$ et $d_+w_0$ dans ces formules (ce qui semble
n\'ecessaire dans \lc); ainsi dans ce cas
$d_*$ est le jumelage naturel de $\r I^{\infty}$ avec lui-m\^eme.

\par Le th\'eor\`eme qui va suivre est donc une g\'en\'eralisation naturelle de ce cas
connu. \`A l'infini de la masure affine on trouve un immeuble jumel\'e \`a la place de
l'immeuble sph\'erique \`a l'infini d'un immeuble affine.

\bigskip
\parni {\bf Th\'eor\`eme 3.7.} {\sl  La codistance $d_*$ d\'efinit un jumelage des
immeubles $(\r I^{+\infty},d_+)$ et $(\r I^{-\infty},d_-)$.    }

\parni {\bf  D\'emonstration.} 

\par Le premier axiome des jumelages (\cf [Tits-92; 2.2] ou [R\'emy-02; 2.5.1])
est clairement v\'erifi\'e. Le troisi\`eme est aussi facile: il suffit de choisir tous
les germes de quartier dans un m\^eme appartement. Pour le second, soient $\g
Q_1\in\r I^{-\infty}$ et $\g Q_2,\;\g Q_3\in\r I^{+\infty}$ des germes de
quartier tels que: $d_*(\g Q_1,\g Q_2)=w\in W^v$ et $d_+(\g Q_2,\g Q_3)=s\in S$,
avec $\ell(ws)=\ell(w)-1$. Soit $A$ un appartement contenant $\g Q_1$ et $\g
Q_2$; on note $\g Q_1=germ_\infty(x+C_1^v)$ et $\g Q_2=germ_\infty(x+C_2^v)$ dans
$A$. On choisit une isom\'etrie de $\A$ sur $A$ qui identifie $C^v_f$ avec $-C_1^v$
et donc $C_2^v$ avec $wC^v_f$. Comme $\ell(ws)=\ell(w)-1$, le mur $M^v$ support
de la cloison $F^v$ de $C_2^v$ de type $\{s\}$ s\'epare $C_2^v$ de $C^v_f=-C_1^v$.
Ainsi $C_1^v$ et $C_2^v$ sont du m\^eme c\^ot\'e de $M^v$ dans $A$ et $C_2^v\subset
cl(C_1^v\cup F^v)$: tout appartement contenant $\g Q_1$ et $\g
F=germ_\infty(x+F^v)$ contient $\g Q_2$. Gr\^ace \`a la proposition 2.7 on voit
facilement que c'est encore vrai si on remplace $\g F$ par un germe de cloison de
quartier parall\`ele. Mais $d_+(\g Q_2,\g Q_3)=s$, donc les cloisons de type $\{s\}$
de deux quartiers quelconques $\g q_2$ et $\g q_3$, de germes $\g Q_2$ et $\g
Q_3$, sont parall\`eles. Ainsi tout appartement $A'$ contenant $\g Q_1$ et $\g Q_3$
contient $\g Q_2$; dans cet appartement le calcul des codistances et distances 
est facile et donne le r\'esultat attendu: $d_*(\g Q_1,\g Q_3)=ws$.
\qed
 
\parni {\bf 3.8. Cons\'equences}\qquad On suppose la masure \r I \'epaisse.

\par 1) La construction d'un immeuble (combinatoire) \'epais $\r I^{\pm \infty}$ de
type $W^v$ implique des limitations sur $W^v$. Si $J\subset I$, on sait
construire \`a partir de $\r I^{\pm \infty}$ un r\'esidu qui est un immeuble \'epais de
type $W^v(J)$ [Ronan-89; III 3.5]. Les limitations connues sur les groupes de
Coxeter des immeubles \'epais sph\'eriques (\cf [Tits-74; addenda], [Tits-Weiss-03;
40.3 et 40.22] ou [Weiss-03; sect. 12.3]) montrent donc que le diagramme de Coxeter
de
$W^v$ ne peut contenir de sous-diagramme de type $H_3$ ou $H_4$. Par contre il n'y
a, a priori, pas de limitation sur les coefficients de la matrice de Coxeter de
$W^v$,
\cf [Ronan-89; exercice 3.21].

\par 2) D'apr\`es [Tits-92; 5.6 cor. 3] et [M\H uhlherr-Ronan-95; (4) p 73], le
jumelage entre $\r I^{+\infty}$ et $\r I^{-\infty}$ montre que, si tous les
coefficients de la matrice de Coxeter sont finis (cas $2-$sph\'erique), tous les
$2-$r\'esidus (sauf peut-\^etre ceux facteurs directs) sont de Moufang. On obtient
alors plus de limitations sur $W^v$: d'apr\`es [Tits-74; addenda] et
[Tits-Weiss-02; 17.1] les coefficients de la matrice de Coxeter de $W^v$ ne
peuvent \^etre que  $1$, $2$, $3$, $4$, $6$ ou $8$ (puisqu'on a exclu $\infty$), 
si l'on suppose de plus toute composante connexe du diagramme de Coxeter de
cardinal au moins 3. Cela ne diff\`ere donc du cas Kac-Moody que par la possibilit\'e
du nombre $8$.

\par On peut am\'eliorer l\'eg\`erement ce r\'esultat. Il est clair que deux r\'esidus
oppos\'es de $\r I^{\infty}$ de type $J\subset I$ (\cf [M\H uhlherr-99; 6 p 135])
forment un immeuble jumel\'e de type $W^v(J)$. Ainsi tous les sous-diagrammes (du
diagramme de Coxeter de $W^v$) connexes, \`a 3 sommets et sans coefficient $\infty$
ne comportent que les coefficients $3$, $4$, $6$ ou $8$.

\par 3) Si le diagramme de Coxeter de $W^v$ est connexe, sans coefficient
$\infty$ et si $\vert I\vert\ge 3$, on sait donc que tous les $2-$r\'esidus de $\r
I^{\pm \infty}$ sont de Moufang. L'immeuble jumel\'e $\r I^{\infty}$ est alors lui-m\^eme
de Moufang [Abramenko-Brown-08; prop. 8.27]. D'apr\`es [Tits-92; 4.3 prop. 7] il
existe alors un groupe $G$ d'automorphismes de $\r I^{\infty}$ qui est muni 
d'une donn\'ee radicielle "jumel\'ee" de type $W^v$ (qui g\'en\'eralise la notion de
donn\'ee radicielle de type $\qF$ de [Rousseau-06; 1.4] d\'ej\`a signal\'ee en
2.4.2). De plus cette donn\'ee radicielle d\'etermine enti\`erement l'immeuble jumel\'e
$\r I^\infty$ [Tits-89; 8.1].

\par La classification des immeubles jumel\'es \'epais de Moufang est assez largement
connue dans le cas $2-$sph\'erique (\ie quand tous les coefficients de la matrice
de Coxeter de $W^v$ sont finis): \cf [Tits-92], [M\H uhlherr-Ronan-95] et [M\H
uhlherr-99-02].

\bigskip
\parni {\bf 3.9. Appartements jumel\'es}

\par Les appartements {\it admissibles} de $(\r I^{\pm \infty},d_\pm )$ ou les
appartements {\it jumel\'es} de $(\r I^{\infty},d_\pm ,d_*)$ sont associ\'es aux paires
de chambres \`a l'infini (\ie de germes de secteurs) oppos\'ees, \cf [R\'emy-02;
2.5.2]. Si on note $A$ l'appartement (unique) contenant deux germes de quartier
oppos\'es $\g Q^+$ et $\g Q^-$, le lemme 3.10 suivant, compar\'e \`a la d\'efinition de
\lc, dit que les appartements admissibles (resp. l'appartement jumel\'e) associ\'es \`a
$\g Q^+$ et $\g Q^-$ sont $A^{\pm \infty}$ (resp.
$A^{\infty}=A^{+\infty}\cup A^{-\infty}$).

\par Inversement tout appartement $A\in\r A$ est l'enveloppe convexe de deux
germes de quartier oppos\'es $\g Q^+$ et $\g Q^-$; donc $A^{\pm \infty}$ est un
appartement admissible de $\r I^{\pm \infty}$ et $A^{\infty}$ un appartement jumel\'e 
de $\r I^{\infty}$. D'apr\`es [Abramenko-96; lemme 2iv p.24] l'appartement
admissible $A^{+\infty}$ (resp. $A^{-\infty}$) d\'etermine enti\`erement son jumeau
$A^{-\infty}$ (resp. $A^{+\infty}$). On a donc des bijections canoniques entre
l'ensemble \r A des appartements de \r I, les ensembles $\r A^{\pm \infty}$ 
d'appartements de $\r I^{\pm \infty}$ et l'ensemble $\r A^*$ des appartements 
jumel\'es de $\r I^{\infty}$.

\bigskip
\parni {\bf Lemme 3.10.} {\sl Soient $\g Q^+\in\r I^{+\infty}$, $\g Q^-\in\r
I^{-\infty}$ deux germes de quartier oppos\'es et $A$ un appartement les contenant.
Pour tout germe de quartier positif $\g Q'$, on a:
\par $d_+(\g Q^+,\g Q')=d_*(\g Q^-,\g Q')$ si et seulement si $\g Q'$ est contenu
dans $A$.
    }

\parni {\bf  D\'emonstration.} Si $\g Q'\subset A$, la relation est \'evidente par
construction. Supposons inversement que $d_+(\g Q^+,\g Q')=d_*(\g Q^-,\g Q')=w$.
Soit $\g Q_0=\g Q^+,\g Q_1,\dots,\g Q_n=\g Q'$ une galerie minimale de germes de
quartiers dans $\r I^{+\infty}$ (donc $n=\ell(w)$). Si $\g Q'\not\subset A$, soit
$i$ le plus petit entier tel que $\g Q_{i}\subset A$ et $\g Q_{i+1}\not\subset
A$. La r\'etraction $\qr=\qr_{A,\sg Q^-}$ conserve les codistances \`a $\g Q^-$ donc
$d_*(\g Q^-,\qr(\g Q'))=w$. Par ailleurs, si $\{s\}$ est le type du germe de
cloison de quartier \g F domin\'e par $\g Q_{i}$ et $\g Q_{i+1}$ et si $\g
Q'_{i+1}$ est le germe de quartier de $A$ adjacent \`a $\g Q_{i}$ le long de \g F,
on a $d_+(\g Q_{i},\g Q_{i+1})=d_+(\g Q_{i},\g Q'_{i+1})=d_+(\g Q'_{i+1},\g
Q_{i+1})=s$; $d_*(\g Q^-,\g Q_{i})=d_+(\g Q^+,\g Q_{i})=w_i\in W^v$; 
$d_*(\g Q^-,\g Q'_{i+1})=d_+(\g Q^+,\g Q'_{i+1})=d_+(\g Q^+,\g Q_{i+1})=w_is$ et
$\ell(w_is)>\ell(w_i)=\ell(w_iss)$. Le second axiome des jumelages nous dit alors
que $d_*(\g Q^-,\g Q_{i+1})=w_i$; donc $d_*(\g Q^-,\qr\g Q_{i+1})=w_i=d_*(\g
Q^-,\g Q_{i})$ et $\qr\g Q_{i+1}=\g Q_i$. La galerie $\qr\g Q_0=\g Q^+,\qr\g
Q_1,\dots,\qr\g Q_n=\qr\g Q'$ b\'egaye donc et ainsi $d_+(\g Q^+,\qr\g Q')$ a une
longueur strictement inf\'erieure \`a $n$, longueur de $d_*(\g Q^-,\qr\g Q')$; ceci
est impossible dans l'appartement $A$.
\qed

\bigskip
\parni {\S$\,${\bf 4.}$\quad${\bftwelve Immeubles microaffines et arbres \`a
l'infini }}
\bigskip
\par On consid\`ere toujours une masure affine \r I

\bigskip \noindent
{\bf 4.1. L'appartement t\'emoin $\A_J$} 

\par Soit $J\subset I$. Le quadruplet $(V,W^v(J),(\alpha_i)_{i\in
J},(\tch{\alpha}{_i})_{i\in J})$ satisfait aux conditions de 1.1. Son syst\`eme de
racines est $\qF_J=\Phi\cap(\bigoplus_{i\in J}\;
\R\alpha_i)=\{\qa\in\qF\mid\qa(F^v(J))=\{0\}\;\}=W^v(J).\{\qa_i\mid i\in J\;\}$.
On note $V_J=\bigcap_{i\in J}\; Ker(\qa_i)$ le support de la facette $F^v(J)$.

\par On consid\`ere l'ensemble $\r M^J$ (resp. $\r M^J_i$) des hyperplans de \r M dont
la direction est de la forme $Ker(\qa)$ avec $\qa\in\qF_J$ (resp. $\qa\in\qD_{im}$
et $\qa(F^v(J))=\{0\}$); ils d\'efinissent une structure d'appartement affine
(mod\'er\'ement imaginaire) sur $\A$ associ\'ee \`a
$(V,W^v(J),(\alpha_i)_{i\in J},(\tch{\alpha}{_i})_{i\in J})$. Tous ces hyperplans
sont stables par translation par $V_J$. Ainsi $\A_J=\A/V_J$ est muni de syst\`emes
$\r M_J$ et $\r M_J^i$ d'hyperplans (en bijection avec $\r M^J$ et $\r M^J_i$) qui
en font un appartement affine essentiel (mod\'er\'ement imaginaire) de groupe de Weyl
vectoriel $W^v(J)$.

\par Le groupe $W^J=W^v(J)\ltimes Q \,\check{ }\subset W$ stabilise $\r M^J$ et
commute \`a $V_J$. Il induit donc un groupe de transformations affines de $\A_J$ 
qui contient le groupe de Weyl $W_J$ de cet appartement (et est contenu
dans le groupe $(W_J)_{P_J}$ correspondant).

\par L'ensemble $\A(F^v(J))$ des germes de faces de quartier de $\A$ de direction
$F^v(J)$ s'identifie de mani\`ere $W^J-$\'equivariante \`a $\A_J$. Un germe de chemin\'ee
\g R de $\A$ de direction dominant $F^v(J)$ d\'efinit un filtre de parties de
$\A_J$ dont une base est obtenue comme suit: \`a un \'el\'ement du filtre \g R on
associe l'ensemble des $\g F\in\A(F^v(J))=\A_J$ contenus dans cet \'el\'ement. Ce
filtre est un germe de chemin\'ee de $\A_J$; si \g R a pour direction $F^v(J)$
c'est une facette-ferm\'ee de $\A_J$ et toutes les facettes-ferm\'ees sont ainsi
obtenues. De m\^eme un germe de face de quartier \g F de direction dominant
$F^v(J)$ d\'efinit un germe de face de quartier (ou un point) de $\A_J$.

\par On s'int\'eressera essentiellement au cas $J$ sph\'erique. Alors $W^v(J)$ et
$\qF_J$ sont finis et il n'y a pas de racine imaginaire pour $\A_J$ ou $\A^J$;
de plus $\A_J$ poss\'ede un produit scalaire
$W^v(J)-$invariant qui en fait un appartement au sens de [Rousseau-08]. Si
l'appartement $\A$ est semi-discret, l'appartement $\A_J$ est discret.

\par Plus g\'en\'eralement on d\'efinit les ensembles $\A(F^v)$ et $\A_{F^v}$ (en
bijection canonique avec $\A_J$) pour une facette vectorielle $F^v$ de type $J$ \`a
la place de $F^v(J)$. On dit que $\A(F^v)$ est la {\it fa{\c c}ade} de $\A$ {\it dans la
direction} $F^v$, \cf [Charignon].

\bigskip \noindent
{\bf 4.2. L'ensemble $\r I(\g F^\infty)$} 

\par Soit $\g F^\infty$ une direction de face de quartier sph\'erique (par exemple
positive) de type $J$. On note $\r I(\g F^\infty)$ l'ensemble des germes de faces
de quartier \g F dans cette direction.

\par Soit $\g f=x+F^v\subset A$ un repr\'esentant d'un $\g F\in \r I(\g F^\infty)$.
On choisit un isomorphisme $\qf\in Isom(\A,A)$ tel que $\qf^{-1}(F^v)=F^v(J)$.
L'ensemble $A(\g F^\infty)$ des \'el\'ements de $\r I(\g F^\infty)$ contenus dans $A$
est identifi\'e par $\qf$ avec $\A(F^v(J))$ et donc avec $\A_J$.

\par Soit \g R un germe de chemin\'ee de \r I de direction dominant $\g F^\infty$
(donc \'evas\'ee); consid\'erons un appartement $A$ contenant \g R. Comme en 4.1
ci-dessus, on peut associer \`a \g R un filtre de parties de $\r I(\g F^\infty)$
(not\'e $\g R_{\sg F^\infty}$) contenu dans $A(\g F^\infty)$ mais ind\'ependant du
choix de $A$. En choisissant $\qf$ comme ci-dessus, $\g R_{\sg F^\infty}$
s'identifie \`a une chemin\'ee de l'appartement $\A_J$; si la direction de \g R est
$\g F^\infty$ alors $\g R_{\sg F^\infty}$ s'identifie \`a une facette de $\A_J$.

\parni{\bf Remarque.} Si $J$ n'est pas sph\'erique, la direction d'une face de
quartier  \g f de type $J$ n'a pas \'et\'e d\'efinie; il est donc a priori plus
difficile de d\'efinir $\r I(\g F^\infty)$ qui devrait \^etre une masure affine.

\bigskip
\parni {\bf Proposition 4.3.} {\sl  L'ensemble $\r I(\g F^\infty)$ muni de ses
{\it facettes} $\g R_{\sg F^\infty}$ (pour les germes de chemin\'ees (\'evas\'ees) \g R
de direction $\g F^\infty$) et de ses {\it appartements} $A(\g F^\infty)$ (pour
$A\in\r A$ contenant un $\g F\in\r I(\g F^\infty))$ est un immeuble affine de
type $\A_J$ (au sens de [Rousseau-08], \cf 2.2.6). Cet immeuble est la {\it fa{\c c}ade}
de \r I {\it dans la direction} $\r F^\infty$.
 }

\parni{\bf Remarque.} Dans le cas fini et si $\g F^\infty$ est la direction de
face de quartier minimale (\ie de type $J=I$), alors $\r I(\g F^\infty)$ est
l'essentialis\'e $\r I^e$ de l'immeuble \r I, c'est \`a dire son quotient par $V_0$.

\parni {\bf  D\'emonstration.} On a vu ci-dessus que les appartements avec leurs
facettes sont isomorphes \`a $\A_J$; d'o\`u (I0). D'apr\`es l'axiome (MA3) deux
facettes $\g R_{\sg F^\infty}$ et $\g R'_{\sg F^\infty}$ sont dans un m\^eme
appartement $A(\g F^\infty)$, d'o\`u (I1); et l'axiome (MA4) dit que cet
appartement est unique \`a un isomorphisme fixant ces deux facettes pr\`es. De plus
l'intersection de deux appartements est une union de facettes d'apr\`es l'axiome
(MA2); d'o\`u (I2).
\qed

 \noindent
{\bf 4.4. L'ach\`evement de Satake de la masure \r I} 

\par On note $\r I^{\qm+}$ (resp. $\r I^{\qm-}$) la r\'eunion (disjointe) des $\r
I(\g F^\infty)$ pour les directions de faces de quartier sph\'eriques positives
(resp. n\'egatives) $\g F^\infty$.Ainsi $\r I^{\qm+}$ (resp. $\r I^{\qm-}$) est
l'ensemble des germes de faces de quartier sph\'eriques positives (resp. n\'egatives)
de \r I.

\par L'{\it ach\`evement de Satake} $\r I^{sat}$ de \r I est la r\'eunion de \r I, $\r
I^{\qm+}$ et $\r I^{\qm-}$. Dans le cas fini $\r I^{\qm+}=\r I^{\qm-}$ contient
$\r I^e$, on se limite au cas essentiel et alors $\r I^{sat}=\r I^{\qm+}=\r
I^{\qm-}$. En dehors du cas fini la r\'eunion est disjointe mais cet ach\`evement
semble inachev\'e, on aimerait lui adjoindre des masures $\r I(\g F^\infty)$ pour
$\g F^\infty$ non sph\'erique.

\par $\r I^{sat}$ est r\'eunion d'appartements $A^{sat}=A\cup A^{\qm+}\cup
A^{\qm-}$ (pour $A\in\r A$) o\`u $A^{\qm\pm }$ est l'ensemble des germes de faces de
quartier sph\'eriques de $A$, positives ou n\'egatives. L'appartement $A^{\qm+}$ est
isomorphe \`a la r\'eunion disjointe $\A^{\qm+}$ des $\A_{F^v}$ pour $F^v$ facette
vectorielle sph\'erique positive de $V$. En fait $\A^{\qm+}$ est l'appartement
t\'emoin de l'immeuble microaffine de [Rousseau-06] dans sa r\'ealisation de Satake
(\ie $\A^{\qm+}=\A^{s}$, \cf [\lc; 4.2.1]).

\par \`A un germe de chemin\'ee \'evas\'ee \g R on associe une facette $\g R^\qm$ situ\'ee
dans $\r I(\g F^\infty)$ o\`u $\g F^\infty$ est la direction de \g R. Si \g R est
dans un appartement $A$, un isomorphisme de $\A$ sur $A$ identifie $\g R^\qm$ \`a
une facette $F^s$ de $\A^s$ (\cf [\lc; 4.2.1]).

\par On peut donc dire que $\r I^{\qm+}$ muni de ses appartements $A^{\qm+}$
(pour $A\in\r A$) et de ses facettes $\g R^\qm$ (pour \g R germe de chemin\'ee
\'evas\'ee positive) est un {\it immeuble microaffine} de type $\A^{\qm+}$ (\ie au
sens de la r\'ealisation de Satake). Aucune d\'efinition abstraite n'a \'et\'e donn\'ee
dans \lc, mais on a les propri\'et\'es habituelles suivantes: les appartements sont
isomorphes \`a $\A^{\qm+}$, deux facettes $\g R_1^\qm$ et $\g R_2^\qm$ sont
contenues dans un m\^eme appartement (\cf (MA3)) et si deux appartements
contiennent $\g R_1^\qm$ et $\g R_2^\qm$, ils sont isomorphes par un isomorphisme
fixant ces deux facettes (\cf (MA4)).

\par Dans le cas fini, $\r I^{\qm+}$ est l'ach\`evement de Satake ou ach\`evement
poly\'edral de l'immeuble (essentiel) \r I [Rousseau-06; 4.3]. On va le munir d'une
topologie qui en fait une compactification de \r I, si ce dernier est
localement fini.

\bigskip \noindent
{\bf 4.5. Topologies sur les immeubles microaffines} 

\par Pour $A\in\r A$, l'appartement microaffine $A^{\qm+}$ est muni d'une
topologie [Rousseau-06; 4.2.4] qui est compacte dans le cas fini. On peut m\^eme
d\'efinir cette topologie sur $A^{\qm-}\cup A\cup A^{\qm+}$ (quand $W^v$ n'a pas de
facteur direct fini). 

\par On va "recoller" ces topologies sur les appartements en deux topologies 
distinctes sur $\r I^{\qm+}$. Pour $\g F\in \r I^{\qm+}$, on doit d\'efinir une
base de voisinages de \g F dans $\r I^{\qm+}$; on va m\^eme d\'efinir des voisinages
dans la r\'eunion disjointe de \r I et $\r I^{\qm+}$.

\par Soit $\g f=x+F^v$ un repr\'esentant de \g F contenu dans un appartement $A$.
Le sous-groupe $W^v_{F^v}$ de $W^v(A)$ fixant $F^v$ est fini, il existe donc une
base de voisinages ouverts de $0$ dans $V_A$ stables par $W^v_{F^v}$. Pour un tel
voisinage ouvert $U$ et un $\qx\in F^v$, l'ensemble $x+\qx+U+F^v$ est un
voisinage ouvert dans $A$ de l'adh\'erence du raccourci $x+\qx+F^v$ de \g f; on note
$\qO_{A,U,\qx}$ la r\'eunion de cet ensemble et des $\g F'\in A^{\qm+}$ qui sont
contenus dedans. Si maintenant $A'$ est un appartement contenant
$germ_{x+\qx}(x+\qx+F^v)$, il existe un isomorphisme $\qf$ de $A$ sur $A'$ fixant
$germ_{x+\qx}(x+\qx+F^v)$ qui est unique \`a $W^v_{F^v}$ pr\`es (au maximum); ainsi
$\qf(\qO_{A,U,\qx})$ ne d\'epend pas du choix de $\qf$, on le note $\qO_{A',U,\qx}$.

\par La base de voisinages de \g F pour la {\it topologie forte} (resp. {\it
faible}) est form\'ee d'ensembles index\'es par $\qx$, $U$ comme ci-dessus, chaque
ensemble \'etant la r\'eunion des $\qO_{A',U,\qx}$ pour $A'$
un appartement contenant $x+\qx+F^v$ (resp. $germ_{x+\qx}(x+\qx+F^v)$). On note
$\r I_s^{\qm+}$ (resp. $\r I_w^{\qm+}$) l'espace $\r I^{\qm+}$ muni de cette
topologie.

\parni {\bf Propri\'et\'es attendues.} Les r\'esultats suivants sont vraisemblables
mais non d\'emontr\'es ici car non utilis\'es dans la suite.

\par 1) Pour les deux topologies les appartements sont ferm\'es et munis de la 
topologie de [Rousseau-06; 4.2.4].

\par 2) Pour $\g F^\infty\in \r I^{+\infty}$ et pour les deux topologies, $\r
I(\g F^\infty)$ est muni de sa topologie naturelle d'immeuble affine euclidien,
il est ouvert dans son adh\'erence qui est la r\'eunion des $\r I(\g F_1^\infty)$
pour $\g F_1^\infty\ge \g F^\infty$.

\par 3) L'ensemble des germes de quartier positifs est discret dans $\r
I_s^{\qm+}$ mais pas dans $\r I_w^{\qm+}$.

\par 4) Dans le cas classique et si \r I est localement fini, $\r I_w^{\qm+}$ est
la compactification de Satake de l'immeuble \r I, \cf [Landvogt-96],
[Charignon].

\bigskip \noindent
{\bf 4.6. L'arbre associ\'e \`a une cloison \`a l'infini $\g F^\infty$}

\par Dans la suite de ce paragraphe on va s'attacher \`a g\'en\'eraliser un certain
nombre de r\'esultats de [Tits-86]. L'un des ingr\'edients essentiels est l'arbre \`a
l'infini associ\'e \`a une cloison de $\r I^{\pm \infty}$: si $\g F^\infty$ est une
cloison \`a l'infini de $\r I^{+\infty}$, c'est \`a dire une direction de cloison de
quartier positive, l'immeuble affine $\r I(\g F^\infty)$ est de dimension $1$ (la
codimension de $\g F^\infty$) donc un arbre.

\par On a sugg\'er\'e en 4.5 que la fronti\`ere de $\r I(\g F^\infty)$ (dans $\r
I_s^{\qm+}$ ou
$\r I_w^{\qm+}$) est form\'ee des $\r I(\g Q)=\{\g Q\}$ pour $\g Q=\g Q^\infty$ un
germe de quartier dominant $\g F^\infty$. Si $\g F^\infty\in A^\infty$, la
fronti\`ere de $A(\g F^\infty)$ est r\'eduite \`a $A(\g Q_1^\infty)\cup A(\g
Q_2^\infty)=\{\g Q_1,\g Q_2\}$, o\`u $\g Q_1$ et $\g Q_2$ sont les deux germes de
quartier de $A$ dominant $\g F^\infty$.

\par Ind\'ependamment de ces r\'esultats, il est facile de voir que $\g Q_1$ et $\g
Q_2$ d\'eterminent les deux bouts de l'appartement $A(\g F^\infty)$ de l'arbre
$\r I(\g F^\infty)$: si $\g q_i=x+C_i^v\subset A$ est un quartier de germe $\g
Q_i$, l'ensemble des germes de cloisons de quartier de $A$ de direction $\g
F^\infty$ et contenus dans
$\g q_i$ est une demi-droite $\qd_i^x$ de la droite $A(\g F^\infty)$ et les
demi-droites
$\qd_1^x$, $\qd_2^x$ sont oppos\'ees.

\par Ainsi $\r I(\g F^\infty)$ est un arbre (discret ou non) muni d'un syst\`eme
d'appartements dont les bouts correspondent bijectivement aux germes de quartiers
de \r I dominant $\g F^\infty$. Ce syst\`eme d'appartements satisfait aux
conditions de [Tits-86]: deux bouts distincts sont toujours dans un m\^eme
appartement (qu'ils d\'eterminent enti\`erement).

\bigskip
\parni {\bf Proposition 4.7.} {\sl  1) Soient $x\in\r I$ et $\g F^\infty\in\r
I^\infty$, il existe une et une seule face de quartier $\g f$ dans \r I de sommet
$x$ et direction $\g F^\infty$.

\par 2) Soient $\g f_1$ et $\g f_2$ deux faces de quartier sph\'eriques dans \r I
de m\^eme direction $\g F^\infty$. Leurs adh\'erences sont disjointes ou elles ont
m\^eme germe. Si de plus elles ont m\^eme sommet, elles co\"{\i}ncident.
 }

\parni{\bf N.B.} On g\'en\'eralise ainsi [Tits-86; prop. 5 et cor. 17.4].

\parni {\bf  D\'emonstration.} 1) Soient $\g F'$ un germe de face de quartier de
direction $\g F^\infty$ et $A$ un appartement contenant $x$ et $\g F'$. Dans $A$
la face de quartier \g f de sommet $x$, parall\`ele \`a $\g F'$ convient: en effet
$\g f^\infty=\g F^\infty$. Si on a deux solutions $\g f'$ et $\g f''$, la
proposition 2.7 nous fournit une suite $\g F_0=\g F',\g F_1,\dots,\g F_n=\g F''$
de germes de faces de quartier parall\`eles et des appartements $A_1,\dots,A_n$
tels que $A_i$ contienne $x$, $\g F_{i-1}$ et $\g F_i$. Tout appartement $A'_i$
contenant $x$ et $\g F_i$ contient une unique face de quartier de sommet $x$ et
direction $\g F^\infty=\g F^\infty_i$; celle-ci est dans l'enclos de x et $\g
F_i$, donc ne d\'epend pas de $A'_i$, on la note $\g f(\g F_i)$. Ainsi on a $\g
f'=\g f(\g F_0)=\g f(\g F_1)=\dots=\g f(\g F_n)=\g f''$.

\par 2) Soit $x\in{\overline{\g f_1}}\cap{\overline{\g f_2}}$, les raccourcis de
$\g f_1$ et $\g f_2$ de sommet $x$ ont m\^eme direction et m\^eme sommet, elles sont
donc \'egales d'apr\`es 1). Ainsi $\g F_1=\g F_2$.
\qed

\parni {\bf Proposition 4.8.} \qquad\cf [Tits-86; 17.3]
{\sl  \par 1) Deux murs (\'eventuellement fant\^omes) de m\^eme direction $M^\infty$
(\cf 3.3.4) sont des hyperplans parall\`eles d'un m\^eme appartement.

\par 2) Soient $M^\infty$ un mur de $\r I^\infty$ et $\g F^\infty\subset
M^\infty$ une cloison. Alors, pour tout germe de cloison de quartier \g F de
direction $\g F^\infty$, il existe un unique mur (\'eventuellement fant\^ome) $M_{\sg
F}$ de \r I de direction $M^\infty$ contenant \g F. Ce mur est un vrai mur si
l'enclos de \g F est \'egal \`a son adh\'erence.
 }

\parni{\bf Remarque.}  Soient $M^\infty$ un mur de $\r I^\infty$ et  $\g
F^\infty\subset M^\infty$ une cloison. Alors, pour tout germe de chemin\'ee \g R de
direction $\g F^\infty$ et pour toute chemin\'ee assez petite \g r de germe \g R,
tous les murs $M_{\sg F}$ (associ\'es aux germes de cloisons de quartier \g F de
direction $\g F^\infty$ contenus dans \g r) sont contenus dans un m\^eme
appartement (dans la d\'emonstration ci-dessous de la partie 2 il suffit de prendre
$A$ contenant
\g r et un $\g F^-$ de direction $\g F^{-\infty}$). On note $\g r+M^\infty$ la
r\'eunion de ces murs. Quand \g r varie, les $\g r+M^\infty$ engendrent un filtre
$\g R+M^\infty$ de parties de \r I (contenu dans l'appartement $A$).

\parni {\bf  D\'emonstration.} 1) Choisissons deux directions de cloison de
quartier oppos\'ees $\g F^{+\infty}$ et $\g F^{-\infty}$ dans $M^{\infty}$. Les
murs $M_1$ et $M_2$ (de direction $M^\infty$) contiennent donc des cloisons de
quartier $\g f^\pm _i$ repr\'esentant ces directions. Soit $A$ un appartement 
contenant les germes
$\g F^-_1$ et $\g F^+_2$; il suffit de montrer qu'il contient $M_1$ et $M_2$.
Faisons le pour $M_1$. La proposition 2.7 nous fournit une suite $\g F_0=\g
F^+_1,\g F_1,\dots,\g F_n=\g F_2^+$ de germes de cloisons de quartier parall\`eles
et des appartements $A_1,\dots,A_n$ tels que $A_i$ contienne $\g F^-_1$, $\g
F_{i-1}$ et
$\g F_{i}$. Pour tout appartement $A'$ contenant $\g F^-_1$ et $\g F_{i}$, le mur
engendr\'e par $\g F^-_1$ ne d\'epend pas de l'appartement $A'$ car il est dans
l'enclos de $\g F^-_1$ et $\g F_{i}$; on le note $M(\g F_i)$. Gr\^ace aux
appartements $A_i$ on a donc $M_1=M(\g F_1^+=\g F_0)=M(\g F_1)=\dots=M(\g
F_n=\g F_2^+)\subset A$.

\par 2) On peut supposer $\g F^\infty=\g F^{+\infty}\in\r I^{+\infty}$. Notons
$\g F^{-\infty}$ sa cloison oppos\'ee dans $M^{\infty}$ et $\g F^-$ un germe de
cloison de quartier de direction $\g F^{-\infty}$. Il existe un appartement $A$
contenant \g F et $\g F^-$. Le support $M$ de \g F dans $A$ est un mur
\'eventuellement fant\^ome; c'est un vrai mur si $cl(\g F)=\overline{\g F}$. Il est
clair dans $A$ que $M^{+\infty}\ni\g F^{+\infty}$, $M^{-\infty}\ni\g
F^{-\infty}$ et donc que le mur de $\r I^\infty$ associ\'e \`a $M$ est bien
$M^\infty$. Soit $M'$ un autre mur avec la m\^eme propri\'et\'e. D'apr\`es la partie 1)
$M$ et $M'$ sont parall\`eles dans un m\^eme appartement; comme ils contiennent le
m\^eme germe de cloison, ils sont \'egaux.
\qed

\parni{\bf 4.9. L'arbre \'etendu associ\'e \`a un mur \`a l'infini}

\par 1) Soit $M^\infty$ un mur de $\r I^\infty$, on note $\r I(M^\infty)$ la
r\'eunion des appartements $A$ de \r I contenant un mur $M$ de direction $M^\infty$
(qui sont appel\'es {\it appartements} de $\r I(M^\infty)$). D'apr\`es la proposition
pr\'ec\'edente, $\r I(M^\infty)$ est r\'eunion disjointe des murs (vrais ou fant\^omes) de
direction $M^\infty$. L'ensemble $\r I^e(M^\infty)$ de ces murs est donc un
quotient de $\r I(M^\infty)$.

\par Choisissons une cloison $\g F^\infty\in M^\infty$ (\eg positive) et notons 
$\g F^{-\infty}$ la cloison oppos\'ee dans $M^\infty$; en fait $M^\infty$ est
d\'etermin\'e par $\g F^{\infty}$ et $\g F^{-\infty}$. La seconde partie de la
proposition pr\'ec\'edente d\'ecrit une application bijective de $\r I(\g F^\infty)$
sur $\r I^e(M^\infty)$. Ainsi $\r I(M^\infty)$ est un immeuble affine inessentiel
dont le quotient essentiel $\r I^e(M^\infty)$ est l'arbre $\r I(\g F^\infty)$. Si
\g R est un germe de chemin\'ee de direction $\g F^\infty$, on a vu que $\g
R_{\sg F^\infty}$ est une facette de $\r I(\g F^\infty)$ (\cf 4.3); l'ensemble
correspondant de $\r I(M^\infty)$ est $\g R+M^\infty$ (\cf remarque 4.8). Quand
\r R varie, ces ensembles $\g R+M^\infty$ d\'ecrivent donc les facettes de
l'immeuble $\r I(M^\infty)$.

\par 2) On va mieux identifier les appartements de $\r I^e(M^\infty)$. Il est
d'abord clair que tout appartement de $\r I(M^\infty)$ induit un appartement de
$\r I(\g F^\infty)$, on a donc la propri\'et\'e suivante (axiome (A5) de [Ronan-89;
appendix 3]):

\par{\bf Propri\'et\'e du Y}: Si $A_1$, $A_2$ et $A_3$ sont trois appartements de \r
I tels que chacune des intersections $A_1\cap A_2$, $A_1\cap A_3$ et $A_3\cap A_2$
 soit un demi-appartement, alors $A_1\cap A_2\cap A_3$ est non vide.

\par En effet, si les murs que sont les bords de deux de ces demi-appartements ne
sont pas parall\`eles dans l'appartement $A_i$ les contenant, la conclusion est
\'evidente. Le cas parall\`ele d\'ecoule de ce que $\r I(\g F^\infty)$ est un arbre si
\g F est un germe de cloison de quartier dans un de ces murs.

\par 3) Un appartement $A$ contenant un germe de cloison de quartier \g F de
direction $\g F^\infty$ d\'etermine deux germes de quartiers $\g Q_1$, $\g Q_2$
dont $\g F^\infty$ est une cloison. L'appartement $A(\g F^\infty)$ de l'arbre $\r
I(\g F^\infty)$ est enti\`erement d\'etermin\'e par ses deux bouts $\g Q_1$ et $\g
Q_2$. Il y a donc (en g\'en\'eral) beaucoup d'appartements $A'$ de \r I d\'eterminant
le m\^eme appartement \`a l'infini $A'(\g F^\infty)=A(\g F^\infty)$. La
proposition ci-dessous montre que parmi eux un et un seul est un appartement de
$\r I(M^\infty)$. On d\'ecrit ainsi la bijection entre les appartements de $\r I(\g
F^\infty)$ et de $\r I(M^\infty)$.

\bigskip
\parni {\bf Proposition 4.10.} {\sl  Dans les conditions pr\'ec\'edentes, il existe un
unique appartement $A'$ contenant \g F, $\g Q_1$ et $\g Q_2$, tel que le mur $M$ 
de $A'$ engendr\'e par \g F ait pour direction $M^\infty$. }
 
\parni {\bf  D\'emonstration.} Si $A'$ est tel qu'indiqu\'e, le mur $M$ est l'unique
mur (\'eventuellement fant\^ome) de 4.8.2. Il est clair que $A'$ est r\'eunion des
enclos $cl(M\cup\g Q_1)$ et $cl(M\cup\g Q_2)$ donc unique.

\par D'apr\`es 4.8.2 il existe un appartement $A$ tel que, dans $A$, le mur $M$
engendr\'e par \g F ait pour direction $M^\infty$. Soit $\g F^-$ le germe de
cloison de quartier oppos\'e \`a \g F dans $M$. Soit $A_i$ un appartement contenant
$\g F^-$ et $\g Q_i$; le mur $M_i$ de $A_i$ engendr\'e par $\g F^-$ a donc comme
direction $M^\infty$ (puisque $\g F^\infty\le \g Q_i$). D'apr\`es l'unicit\'e de
4.8.2 on a $M_i=M$. L'enveloppe convexe ferm\'ee de $\g F^-$ et $\g Q_i$ (\ie
l'enveloppe convexe ferm\'ee de $M$ et $\g Q_i$) dans $A_i$ est un demi-appartement
$D_i$ (\'eventuellement fant\^ome) de bord $M$. Dans $A_1$ consid\'erons l'image $\g
Q_3$ par la r\'eflexion $r_{M}$ de $-\g Q_1$. On a donc $\g Q_3\subset D_1$ et $(\g
F^-)^\infty\le \g Q_3$. Soit $A'$ un appartement contenant $\g Q_2$ et $\g
Q_3$, on a $\g F^\infty,(\g F^-)^\infty\subset A'^\infty$, donc $M^\infty\subset
A'^\infty$ et $A'$ est r\'eunion de murs de direction $M^\infty$. Par unicit\'e des
murs de direction $M^\infty$ contenant un germe de cloison (4.8.2) $A'$ contient
tout le demi-appartement de $D_1\subset A_1$ engendr\'e par $\g F^\infty$ et un
quartier $\g q_3\subset D_1$ de germe $\g Q_3$; en particulier $A'\supset\g Q_1$.
Comme il contient aussi $\g Q_2$ il contient \g F et donc (par unicit\'e) $M$ qui
est de direction $M^\infty$.
\qed

\parni {\bf Proposition 4.11.} {\sl Soient $A_1$ et $A_2$ deux appartements de \r
I tels que $A_1\cap A_2$ soit un demi-appartement $D$, alors il existe une
r\'etraction $\qr$ (de centre un germe de quartier) de \r I sur $A_1$ telle que
$\qr(A_2)=D$. 
 }

\parni{\bf Remarque.}  A part la condition de diminution des distances (qui n'a
un sens que pour les immeubles affines) ce r\'esultat est l'axiome (A5) de
[Tits-86].
\parni {\bf  D\'emonstration.} Soient $M$ le bord de $D$ et \g F un germe de
cloison de quartier contenu dans le mur $M$. Tout se passe dans $\r I(M^\infty)$
et donc principalement dans l'arbre $\r I(\g F^\infty)$. Il suffit de choisir le
germe de quartier \g Q de $A_1\setminus D$ dont une cloison est $\g F^\infty$ et
de consid\'erer la r\'etraction $\qr_{A_1,\sg Q}$.
\qed

\parni{\bf 4.12. Le probl\`eme de classification}

\par Le th\'eor\`eme 2 de [Tits-86; p 166] dit qu'un immeuble affine \r I est
d\'etermin\'e \`a isomorphisme pr\`es par son immeuble sph\'erique \`a l'infini $\r I^\infty$
(suppos\'e \'epais) et les arbres $\r I(\g F^\infty)\simeq\r I(M^\infty)$ pour $\g
F^\infty$ cloison de $\r I^\infty$ contenue dans un mur $M^\infty$. Si $\r
I^\infty$ est de Moufang, ces arbres d\'eterminent une valuation de la donn\'ee
radicielle $G$ associ\'ee \`a $\r I^\infty$ [\lc; th 3 p 170] et cette valuation
permet de reconstruire l'immeuble affine \r I;
on obtient ainsi la classification des immeubles affines irr\'eductibles localement
finis de dimension $\ge 3$ [\lc; 14, 15].

\par On peut envisager une strat\'egie semblable pour classifier les masures
affines. On a vu en 3.8.3 qu'\`a une masure affine \'epaisse $2-$sph\'erique est
souvent associ\'ee une donn\'ee radicielle "jumel\'ee" et que ces donn\'ees radicielles
sont classifi\'ees dans de nombreux cas. Les r\'esultats pr\'ec\'edents sur les arbres $\r
I(\g F^\infty)$ ou $\r I(M^\infty)$ doivent permettre de d\'efinir une valuation de
la donn\'ee radicielle. Malheureusement il n'y a pas encore de construction d'une
masure affine associ\'ee \`a une donn\'ee radicielle valu\'ee (hors le cas du paragraphe
6). Ainsi ce processus de classification s'arr\^ete l\`a (pour l'instant en tout cas).

\bigskip
\parni {\S$\,${\bf 5.}$\quad${\bftwelve Propri\'et\'es \`a distance finie des masures }}
\bigskip
\par On consid\`ere une masure affine \r I que l'on supposera ordonn\'ee \`a partir de
5.4.

\bigskip \noindent
{\bf Proposition 5.1.} {\sl Soient $A$ un appartement de la masure \r I, $x_1$ et
$x_2$ deux points de $A$ tels que $x_1\le _Ax_2$. On consid\`ere des facettes locales
de \r I, $F_1$ et $F_2$ de sommets respectifs $x_1$ et $x_2$. Si $F_1$ est
positive et $F_2$ n\'egative, on suppose de plus $x_1=x_2$ ou $x_2-x_1\in\r
T^\circ(A)$ (c\^one de Tits ouvert de $A$) \ie $x_1\bc\le _A x_2$. Alors il existe un
appartement $A'$ contenant $F_1$ et $F_2$.
 }

\parni {\bf  D\'emonstration.} 

\par Quitte \`a \'echanger $x_1$, $x_2$ et changer les signes, on peut supposer
v\'erifi\'ee l'une des trois conditions suivantes: $F_2$ positive, $x_1=x_2$ ou 
$x_2-x_1\in\r T^\circ(A)$. Soit \g q un quartier (n\'egatif si $x_1\ne x_2$) de sommet
$x_2$ dans
$A$ tel que $x_1\in{\overline{\g q}}$; alors, par convexit\'e, il existe un 
appartement $A_1$ contenant \g q et $F_2$ et cet appartement est isomorphe \`a $A$
par un isomorphisme fixant \g q. Ainsi $x_1\le _{A_1}x_2$, et m\^eme $x_2-x_1\in\r
T^\circ(A_1)$, si c'est vrai dans $A$. On s'est donc ramen\'e au cas o\`u $F_2\subset
A$.

\par On voit facilement qu'il existe dans $A_1$ un quartier (positif si
$x_1\ne x_2$) $\g q_1$ de sommet $x_1$ tel que $F_2\subset{\overline{\g q_1}}$. Par
convexit\'e il existe un appartement $A'$ contenant $F_1$ et $\g q_1$, donc aussi
$F_2$.
\qed

 \noindent
{\bf Proposition 5.2.} {\sl Sous les hypoth\`eses de la proposition 5.1, supposons
de plus $x_1=x_2$ et les deux facettes locales $F_1$, $F_2$ de m\^eme signe. Alors
deux appartements contenant $F_1$ et $F_2$ sont isomorphes par un isomorphisme
fixant $F_1$ et $F_2$.
 }
\parni{\bf N.B.} Il r\'esulte de la d\'emonstration qui suit que, si deux
appartements contiennent une m\^eme chambre locale $C_0$ en $x$, leur intersection
est (localement)  convexe par galeries: toute galerie tendue d'une chambre locale
en $x$ \`a une facette locale en $x$ de cette intersection est enti\`erement contenue
dans cette intersection. De plus tout isomorphisme fixant $C_0$ entre les deux
appartements fixe toutes les chambres de ces galeries.

\parni {\bf  D\'emonstration.} 
\par On se ram\`ene par le proc\'ed\'e habituel [Brown-89; IV 1] au cas o\`u l'une des
deux facettes (\eg $F_2$) est une chambre. Soient $A$ un appartement contenant
$F_1$, $F_2$ et $\g q_2$ le quartier de sommet $x=x_1=x_2$ dans $A$ engendr\'e par
$F_2$. Si $A'$ est un autre appartement contenant $F_1$ et $F_2$, on consid\`ere
dans $A'$ une galerie de chambres locales en $x$ : $F_2=C_0,C_1,\dots,C_n\supset 
F_1$ qui est minimale (de $C_0$ \`a $F_1$). En particulier le type de cette galerie
donne un mot r\'eduit de $W^v$, donc dans tout appartement $A''$ contenant
$C_0,C_1,\dots,C_{k-1}$, ces chambres locales sont toutes du m\^eme c\^ot\'e du mur
(\'eventuellement fant\^ome) de $A''$ contenant la cloison
${\overline{C_{k-1}}}\cap{\overline{C_{k}}}$.

\par En utilisant 2.9.1 successivement pour les murs s\'eparant $C_{k-1}$ et
$C_{k}$, on trouve un appartement $A_2$ de \r I contenant $\g q_2$,
$C_0,C_1,\dots,C_n$ et donc $F_1$; de plus la galerie $C_0,C_1,\dots,C_n$ est
minimale de $F_2$ \`a $F_1$ dans cet appartement $A_2$. D'apr\`es (MA4) $A\supset
cl_{A_2}(F_1,\g Q_2)=cl_{A_2}(F_1,\g q_2)\supset
cl_{A_2}(F_1,C_0)\supset{\overline{C_0}}\cup{\overline{C_1}}\cup\dots\cup
{\overline{C_n}}$ et il existe un isomorphisme $\qf:A_2\rightarrow A$ fixant cet
enclos. D'apr\`es (MA2) appliqu\'e \`a $C_0=F_2$, il existe un isomorphisme
$\psi:A'\rightarrow A_2$ qui fixe $C_0$. Mais les chambres locales
$C_0,C_1,\dots,C_n$ sont dans $A'$ et $A_2$, donc, par r\'ecurrence, $\psi$ fixe
chaque $C_j$. Finalement l'isomorphisme $\qf\circ\psi$ fixe $F_1$ et $F_2$, c'est
l'isomorphisme cherch\'e.
\qed

 \noindent
{\bf 5.3. Immeubles r\'esiduels} \qquad Soit $x$ un point de la masure \r I.

\par 1) On note $\r I^+_x$ (resp. $\r I^-_x$) l'ensemble des germes de segments
positifs (resp. n\'egatifs) $[x,y)$ de \r I et $\r I_x=\r I^+_x\cup\r I^-_x$. Si
$A$ est un appartement de \r I contenant $x$, on note $A^+_x$ (resp. $A^-_x$)
l'ensemble de ces germes de segments positifs (resp. n\'egatifs) contenus dans $A$
et $A_x=A^+_x\cup A^-_x$. Dans le cas fini $A_x=A^+_x= A^-_x$ (resp. $\r I_x=\r
I^+_x=\r I^-_x$) est la "sph\`ere unit\'e tangente" en $x$ \`a $A$ (resp. \r I). Si $F$
est une facette ou une facette locale de \r I positive (resp. n\'egative) contenant
$x$ dans son adh\'erence, on note $F_x$ l'ensemble des germes de segments positifs
(resp. n\'egatifs) $[x,y)$ pour $]x,y)\subset F$.

\par 2) Soient $\qe=\pm $ et $A$ un appartement contenant $x$. Le choix de $x$ comme
origine identifie $A$ \`a $V(A)$. Alors $A^\qe_x$ est identifi\'e \`a l'ensemble
$V(A)^{s\qe}$ des demi-droites d'origine $0$ de $\qe\r T(A)\subset V(A)$; il
h\'erite des deux structures de complexes de Coxeter que l'on va d\'efinir sur $\qe\r
T(A)$:

\par - La structure non restreinte est celle d\'ecrite en 1.3 (pour $V$), son
groupe de Weyl est $W^v(A)\simeq W^v$. Les facettes correspondantes dans $A_x^\qe$
sont les $F^\ell_x$ o\`u $F^\ell$ est une facette locale en $x$ de signe
$\qe$ contenue dans $A$. Les murs correspondants sont \'eventuellement fant\^omes.

\par - La structure restreinte a pour groupe de Weyl le sous-groupe $W^v(A)_x$ de
$W^v(A)$ engendr\'e par les r\'eflexions vectorielles par rapport aux directions des
(vrais) murs de $A$ contenant $x$. Ce groupe est un groupe de Coxeter pour un
ensemble de g\'en\'erateurs \'eventuellement infini [D\'eodhar-89]; des r\'esultats
compl\'ementaires sont prouv\'es dans [Tits-88] ainsi que dans [Bardy-96; 5.1] ou
[Moody-Pianzola-95; 5.7] (o\`u la restriction \`a un contexte l\'eg\`erement diff\'erent
peut \^etre lev\'ee), la meilleure r\'ef\'erence ici semble \^etre [H\'ee-90],
malheureusement peu disponible. Il r\'esulte de ces raisonnements que
$A_x^\qe$ est muni d'une structure de  complexe de Coxeter dont les murs sont les
(vrais) murs de $A$ contenant $x$ et dont les facettes sont les $F_x$ o\`u $F$ est
une facette de $A$ de signe $\qe$ dont l'adh\'erence contient $x$. Ce complexe de
Coxeter peut cependant \^etre tronqu\'e,  car l'ensemble des facettes $F_x$ dans
l'adh\'erence d'une chambre $C_x$ peut ne pas  \^etre en bijection (d\'ecroissante)
avec l'ensemble des parties de l'ensemble des cloisons de $C_x$. Par contre le
groupe $W^v(A)_x$ (agissant sur $A$ en fixant $x$) est bien simplement transitif
sur les chambres $C_x$ de m\^eme signe et engendr\'e par l'ensemble $S(C^0_x)$
(\'eventuellement infini) des r\'eflexions par rapport aux murs de l'une de ces
chambres $C^0_x$, voir [Tits-88; prop.3]. L'axiome (MA2) et cette simple 
transitivit\'e montrent que le syst\`eme de Coxeter
$(W^v(A)_x,S(C^0_x))$ ne d\'epend pas, \`a isomorphisme unique pr\`es, des choix de $A$
et $C^0_x$, on le note $(W^v_x,S_x)$; bien sur $W^v_x$ est un sous-groupe de
$W^v$.

\par 3) On munit les appartements $A^\qe_x$ de $\r I^\qe_x$  de l'une des deux 
structures de complexe de Coxeter d\'ecrites ci-dessus. Alors les propositions 5.1
et 5.2 montrent que $\r I^\qe_x$ avec ses appartements est un immeuble au sens
classique, \cf \eg [Brown-89; IV 1].

\par On note $Ch_r(\r I^\qe_x)$ l'ensemble des chambres de $\r I^\qe_x$ pour la
structure restreinte. Il est muni d'une $W^v_x-$distance $d_\qe^r$ (v\'erifiant les
axiomes de [Tits-92; 2.1] ou [R\'emy-02; 2.3.1] puisqu'elle traduit la structure
d'immeuble): on choisit un appartement $A_0$ contenant $x$ et une chambre
positive $C^0$ de $A_0$ contenant $x$ dans son adh\'erence. Si
$C_x,C'_x\in Ch_r(\r I^\qe_x)$, on choisit un appartement $A$ contenant les
chambres correspondantes $C$, $C'$ et on identifie $A$ \`a $V(A)$ par le choix de
$x$ pour origine. D'apr\`es (MA2)  il existe un isomorphisme $\qf$ de $(A_0,x)$ sur
$(A,x)$. Alors, si $\qf^{-1}(C)=\qe wC^0$ et $\qf^{-1}(C')=\qe w'C^0$ pour
$w,w'\in W^v(A_0)_x\simeq W^v_x$, on d\'efinit
\qquad$d^r_\qe(C_x,C'_x)=w^{-1}w'$\qquad\quad; ceci ne d\'epend pas des choix
effectu\'es.

\par On note $Ch_{nr}(\r I^\qe_x)$ l'ensemble des chambres de $\r I^\qe_x$ pour la
structure non restreinte. Il est muni d'une $W^v-$distance $d_\qe^{nr}$, d\'efinie
de mani\`ere analogue \`a la pr\'ec\'edente et qui traduit la structure non restreinte
d'immeuble sur $\r I^\qe_x$, \cf [Gaussent-Rousseau-08; sect. 4.5].
\medskip
\par On suppose dans le reste de ce \S\kern 2pt 5 la masure \r I {\bf ordonn\'ee}.
\bigskip
 \noindent
{\bf Proposition 5.4.} {\sl Dans la masure ordonn\'ee \r I, on consid\`ere deux
appartements $A_1$, $A_2$ et deux points distincts $x,y$ de $A_1\cap A_2$
v\'erifiant $x\le _{A_1}y$ (resp. $x\bc\le _{A_1}y$). Alors $x\le _{A_2}y$ (resp.
$x\bc\le _{A_2}y$) et il existe un isomorphisme de
$A_1$ sur $A_2$ qui fixe $cl_{A_1}([x,y])$.
 }
\parni{\bf Cons\'equences.} 
\par 1) On peut donc d\'efinir des relations $\le $ et $\bc\le $ dans la masure
ordonn\'ee \r I:

\par Pour tous $x,y\in \r I$, on dit que  $x\le y$ (resp. $x\bc\le y$) si et
seulement si il existe un appartement $A$ de \r I contenant $x,y$ tel que
$x\le _Ay$ (resp. $x\bc\le _{A}y$) (ceci ne d\'epend pas du choix de $A$).

\par On va voir en 5.9 que $\le $ ou $\bc\le $ est un pr\'eordre. Bien sur $x\bc\le y
\Rightarrow x\le y$.

\par 2) Pour $x\le y$ dans \r I, le segment $[x,y]$ et m\^eme son enclos est
ind\'ependant de l'appartement choisi contenant $x$ et $y$. On peut en particulier
parler de {\it convexit\'e pr\'eordonn\'ee} dans \r I.

\parni {\bf  D\'emonstration.} 
\par 1) D'apr\`es l'axiome (MAO), on a $[x,y)_{A_1}=[x,y)_{A_2}\subset A_1\cap
A_2$; l'axiome (MA2) dit donc qu'il existe un isomorphisme $\qf$ de $A_1$ sur
$A_2$ qui fixe $[x,y)$. Ainsi, dans $A_2$, $[x,y)$ est un germe de segment positif
(resp. et g\'en\'erique) et donc $x\le _{A_2}y$ (resp. $x\bc\le _{A_2}y$).

\par 2) Notons $\qD_1$ (resp. $\qD_2$) la demi-droite de $A_1$ (resp. $A_2$)
d'origine $x$ (resp. $y$) et contenant $y$ (resp. $x$). Ces deux demi-droites
sont pr\'eordonn\'ees et ont en commun le segment $[x,y]$. On va prouver l'existence
d'un appartement les contenant toutes les deux et dans lequel leur r\'eunion
$\qD$ est une droite.

\par Soit $\g Q_2$ un germe de quartier de $A_2$ contenant la direction de
$\qD_2$ dans son adh\'erence. D'apr\`es (MA3) et par compacit\'e de
$\qD_1\cup\{\infty\}$ (\cf la preuve de 2.7) il existe des points distincts 
$y_1=x,y_2=y,y_3,\dots,y_n$ dans cet ordre sur $\qD_1$ et des appartements
$A_2,A_3,\dots,A_n,A_{n+1}$ contenant tous $\g Q_2$ et respectivement
$[y_1,y_2],[y_2,y_3],\dots,[y_{n-1},y_n]$, $[y_n,\infty]=\qD_1\setminus[x,y_n[$.
L'appartement $A_2$ contient $\qD_2$. L'appartement $A_3$ contient $y_2$ et $\g
Q_2$ et donc aussi $\qD_2$ (par convexit\'e, \cf (MA4)). Mais
$y=y_2\in[x,y_3]_{A_1}$ et $x\le _{A_1}y_3$; d'apr\`es (MAO),
$[x,y_3]_{A_1}=[x,y_3]_{A_3}$ et $y_2\in[x,y_3]$. Ainsi
$\qD_3=\qD_2\cup[y_2,y_3]$ est une demi-droite de $A_3$ d'origine $y_3$ dans
l'enveloppe convexe de $\g Q_2$ et $y_3$. On est dans la situation du d\'epart avec
$A_2$ devenu inutile. On a donc montr\'e, par r\'ecurrence sur $n$, que $A_{n+1}$
contient $\qD=\qD_1\cup\qD_2$  et que $\qD$ est une droite de $A_{n+1}$.

\par 3)  Un meilleur choix de l'appartement:

\par Soit $C_1^v$ une chambre vectorielle de $A_1$ et $F_1^v$ la facette
vectorielle engendr\'ee par la direction de $\qD_1$. On consid\`ere dans $A_1$ la
chemin\'ee pleine $\g r_1=cl(germ_y(y+C_1^v)+F_1^v)$ et son germe $\g R_1$. Il est
clair que $cl(x,\g R_1)\supset \g r_1\supset \qD_1$. Dans 2) ci-dessus on peut
supposer que $A'=A_{n+1}$ contient $\g Q_2$ et $\g R_1$; on a vu qu'il contient
$\qD_1$ et $\qD_2$ en particulier $x$ et $y$.  Ainsi $A'$ contient le quartier
$\g q_2$ de $A_2$ de sommet $y$ et direction $\g Q_2$ ainsi que la chemin\'ee
pleine $\g r_1$. D'apr\`es l'axiome (MA4) on a un isomorphisme $\qf$ de $A_1$ sur
$A'$ fixant $cl(x,\g R_1)\supset cl([x,y])$ et un isomorphisme $\psi$ de $A'$ sur
$A_2$ fixant $cl(\g q_2)\supset([x,y])$. Ainsi $\psi\circ\qf$ est l'isomorphisme
cherch\'e.
\qed

 \par\noindent
{\bf Proposition 5.5.} {\sl Soient $A$ un appartement de la masure ordonn\'ee \r
I,  
$x_1$ et $x_2$ deux points de $A$ tels que $x_1\le _A x_2$. On consid\`ere des facettes 
locales de \r I, $F_1$ et $F_2$ de sommets respectifs $x_1$ et $x_2$. On suppose
$x_2-x_1\in\r T^\circ(A)$ ou $F_2$ positive et $F_1$ n\'egative. Alors, si des
appartements $A'$ et $A''$ contiennent $F_1$ et $F_2$, ils contiennent leur
enveloppe convexe et il existe un isomorphisme de $A'$ sur $A''$ qui fixe cette
enveloppe convexe.
 }
\parni{\bf Remarque.} On a montr\'e l'existence d'appartements $A'$ et $A''$ comme
dans l'\'enonc\'e en 5.1 (sous une hypoth\`ese moins restrictive). Avec l'hypoth\`ese
ci-dessus, si $\qO_2$ (resp. $\qO_1$) est un \'el\'ement assez petit du filtre $F_2$
(resp. $F_1$) contenu dans $A'$, alors pour tous $y_2\in\qO_2$ et $y_1\in\qO_1$
on a $y_1\le y_2$.

\parni {\bf  D\'emonstration.} On suppose $\qO_2$ et $\qO_1$ comme ci-dessus,
convexes et contenus dans $A'$ et $A''$. D'apr\`es 5.4 $A'\cap A''$ contient
$[y_1,y_2]_{A'}=[y_1,y_2]_{A''}$ pour tous $y_2\in\qO_2$ et $y_1\in\qO_1$. Ainsi
$A'\cap A''$ contient l'enveloppe convexe de $\qO_2\cup \qO_1$ et donc de
$F_2\cup F_1$. Pour construire l'isomorphisme de $A'$ sur $A''$ on se ram\`ene (par
le proc\'ed\'e habituel [Brown-89; IV 1]) au cas o\`u l'une des facettes locales (\eg
$F_1$) est une chambre.

\par Soient $y_1\in\qO_1$ tel que $F_1\subset cl([x_1,y_1))$ et $z_1$ le milieu
de $[x_1,y_1]$. Alors $z_1$ est dans un ouvert de $A'$ ou $A''$ contenu dans 
$cl([x_1,y_1])\cap\qO_1$. D'apr\`es 5.4 il existe un isomorphisme $\qf$ de $A'$ sur
$A''$ qui fixe
$cl([x_1,y_1])$. Soit
$z_2\in\qO_2$, alors $[z_1,z_2]\subset A'\cap A''$ rencontre $cl([x_1,y_1])$
selon un voisinage ouvert de $z_1$ dans $[z_1,z_2]$. Donc $\qf$ est une bijection
affine de $[z_1,z_2]$ sur son image dans $A''$ qui est l'identit\'e sur un
voisinage de $z_1$. Or 5.4 appliqu\'e \`a $[z_1,z_2]$ dit que l'identit\'e de
$[z_1,z_2]$ est une bijection affine (pour les structures affines h\'erit\'ees
respectivement de $A'$ et $A''$), donc $\qf$ co\"{\i}ncide avec cette bijection
affine; en particulier $\qf(z_2)=z_2$. Ainsi $\qf$ fixe $\qO_2$ et
$cl([x_1,y_1])$, donc $\qO_1$ (quitte \`a le remplacer par un \'el\'ement plus petit de
$F_1$). Pour $t_1\in\qO_1$, $t_2\in\qO_2$ le m\^eme raisonnement que ci-dessus
montre que $\qf$ fixe $[t_1,t_2]$. Mais l'enveloppe convexe de $\qO_1\cup\qO_2$
est la r\'eunion de ces segments, d'o\`u la proposition.
\qed

\par \noindent
{\bf 5.6. Jumelage des immeubles r\'esiduels} 

\par Pour un point $x$ de la masure ordonn\'ee \r I, on a d\'efini (en 5.3) deux
ensembles $\r I^+_x$, $\r I^-_x$ et deux structures d'immeubles sur chacun d'eux,
la non restreinte de groupe de Weyl $W^v$ et la restreinte de groupe de Weyl
$W^v_x$. Pour chacune de ces structures on peut d\'efinir un jumelage de $\r
I^+_x$ et $\r I^-_x$. Le cas non restreint \'etant trait\'e en d\'etail dans
[Gaussent-Rousseau-08], on se concentre ici sur le cas restreint.

\par On d\'efinit ci-dessous une application {\it codistance}:
\par\qquad\qquad $d_*^r\; :\; Ch_r(\r I^+_x)\times Ch_r(\r I^-_x)\cup Ch_r(\r
I^-_x)\times Ch_r(\r I^+_x)\rightarrow W^v_x$

\par On choisit un appartement $A_0$ contenant $x$ et une chambre positive $C^0$
de $A_0$ contenant $x$ dans son adh\'erence. Si $\qe=\pm $,  $C_x\in Ch_r(\r
I^{\qe}_x)$ et $C'_x\in Ch_r(\r I^{-\qe}_x)$, on choisit un appartement $A$ 
contenant les chambres correspondantes $C$, $C'$ et on identifie $A_0$ \`a $V(A_0)$
par le choix de $x$ comme origine. D'apr\`es (MA2) il existe un isomorphisme $\qf$
de 
$(A_0,x)$ sur $(A,x)$. Alors, si $\qf^{-1}(C)=\qe wC^0$ et 
 $\qf^{-1}(C')=-\qe w'C^0$ pour $w,w'\in W^v(A_0)_x\simeq W^v_x$, on d\'efinit
$d_*^r(C_x,C'_x)=d_{-\qe}^r(-C_x,C'_x)=d_\qe^r(C_x,-C'_x)=w^{-1}w'$; ceci ne
d\'epend pas des choix effectu\'es d'apr\`es 5.5.

{\bf Proposition.} {\sl $d_*^r$ est un jumelage des immeubles $\r I^+_x$ et 
$\r I^-_x$ munis de leurs structures restreintes.
 }
\parni{\bf Remarque.} On a de m\^eme un jumelage $d_*^{nr}$ (codistance \`a valeurs
dans $W^v$) de ces immeubles munis de leurs structures non restreintes.

\parni {\bf  D\'emonstration.} 
\par Il faut v\'erifier les axiomes de [Tits-92; 2.2], voir
aussi [R\'emy-02; 2.5.1].

\par On a $d_*^r(C'_x,C_x)=w'^{-1}w=d_*^r(C_x,C'_x)^{-1}$, d'o\`u (Tw1). Soient
maintenant $C_x\in Ch_r(\r I^{\qe}_x)$ et $C'_x,C''_x\in Ch_r(\r I^{-\qe}_x)$ des
chambres telles que $d_*^r(C_x,C'_x)=w\in W^v_x$ et $d_{-\qe}^r(C'_x,C''_x)=s\in
S_x$, avec $\ell(ws)=\ell(w)-1$ (longueurs calcul\'ees dans le syst\`eme de Coxeter
$(W^v_x,S_x)$ ). Soit $A$ un appartement contenant $C$ et $C'$ (et donc $x$); on
choisit de r\'ealiser $(W^v_x,S_x)$ dans $A$ par le choix de $\qe C$ comme chambre
fondamentale. Comme $\ell(ws)=\ell(w)-1$, le (vrai) mur $M$ engendr\'e par la
cloison de $-\qe C'$ de type $\{s\}$ s\'epare $-\qe C'$ de $\qe C$; ainsi $C$ et
$C'$ sont du m\^eme c\^ot\'e de $M$. D'apr\`es 2.9.1 il existe un appartement $A'$
contenant $C$, $C'$ et $C''$. Dans cet appartement $C'=w.(-C)$ et
$C''=(wsw^{-1}).C'$ donc $C''=ws.(-C)$ et $d_*^r(C_x,C''_x)=ws$, c'est l'axiome
(Tw2).

\par Pour v\'erifier le troisi\`eme axiome (Tw3), soient $C_x\in Ch_r(\r I^{\qe}_x)$,
$C'_x\in Ch_r(\r I^{-\qe}_x)$, $w=d_*^r(C_x,C'_x)\in W^v_x$ et $s\in S_x$. Dans
un appartement $A$ contenant $C$ et $C'$, consid\'erons la chambre $C''\ne C'$
adjacente \`a $C'$ le long de la cloison de type $\{s\}$. On a bien
$d_{-\qe}^r(C'_x,C''_x)=s$ et $d_*^r(C_x,C''_x)=ws$ comme demand\'e.
\qed

\par \noindent
{\bf Proposition 5.7.} {\sl Soient $C^-$, $C$ deux chambres n\'egatives d'un
immeuble jumel\'e (de groupe de Weyl $W$) et $C^+$ une chambre positive oppos\'ee \`a
$C^-$. On consid\`ere un appartement jumel\'e $A=A^-\cup A^+$ contenant $C$ et la
r\'etraction
$\qr=\qr_{A,C}$ de l'immeuble sur $A$ de centre $C$. On note $w^-$ et $w^+$
les \'el\'ements de $W=W(A)$ tels que $\qr(C^-)=w^-C$ et $\qr(C^+)=-w^+C$ (la chambre
de $A^+$ oppos\'ee \`a $w^+C\in A^-$). Alors $w^+\le w^-$ pour l'ordre de
Bruhat-Chevalley sur $W(A)$ correspondant au choix du syst\`eme de g\'en\'erateurs
form\'e des r\'eflexions par rapport aux murs de $C$ (ou $-C$).
 }

\parni {\bf  D\'emonstration.} On a $w^+=d_*(C,C^+)$ et $w^-=d_-(C,C^-)$. La
relation $d_*(C,C^+)\le d_-(C,C^-)$ ou plut\^ot $d_*(C^+,C)\le d_-(C^-,C)$ est exactement
ce qu'obtient Peter Abramenko au cours de la d\'emonstration de
$\ell(d_*(C^+,C))\le \ell(d_-(C^-,C))$ dans sa Remark 3 page 24 de [Abramenko-96].
\qed

\par \noindent
{\bf Corollaire 5.8.} {\sl Soient $\g Q=germ_\infty(x_0-C^v)$ un germe de
quartier n\'egatif dans un appartement $A$ de la masure ordonn\'ee \r I,
$\qr=\qr_{A,\sg Q}$ la r\'etraction sur $A$ de centre \g Q (\cf 2.6) et $x\le y$ deux
points de \r I. Alors l'image $\qr([x,y])$ du segment $[x,y]$ dans $A$ est une
ligne bris\'ee qui est "pli\'ee positivement", c'est \`a dire:

\par pour tout $z_1=\qr(z)\in\qr(]x,y[)$, on note $\qp_+=\qr([z,y))$ (resp.
$\qp_-=\qr([z,x))$ ) et $w_+$ (resp. $w_-$) l'\'el\'ement minimal de $W^v$ tel que
$\qp_+\subset w_+(z_1+{\overline{C^v}})$ (resp. $\qp_-\subset
w_-(z_1-{\overline{C^v}})$), on a alors
$w_+\le w_-$ pour l'ordre de Bruhat-Chevalley.
 }
\parni{\bf Remarques.} 1) L'ordre de Bruhat-Chevalley choisi sur $W^v=W^v(A)$
correspond au choix du syst\`eme de g\'en\'erateurs form\'e des r\'eflexions par rapport aux 
murs de la chambre vectorielle $C^v$ de $A$.

\par 2) La relation $\le $ et le segment $[x,y]$ sont bien d\'efinis dans \r I (5.4).
On a d\'emontr\'e en 2.8 que $\qr([x,y])$ est une ligne bris\'ee "croissante", en
particulier $\qp_+$ (resp. $\qp_-$) est un germe de segment positif (resp.
n\'egatif) en $z_1$.

\parni {\bf  D\'emonstration.} D'apr\`es (MA2) on peut supposer que $A$ contient \g Q
et $[z,x)$, donc $z=z_1$ et $\qp_-=[z,x)$. Consid\'erons la chambre locale (en $z$)
$C=germ_z(z-C^v)$. Par d\'efinition de $w_-$ on a $\qp_-\subset
{\overline{C^-}}=w_-{\overline{C}}$. Dans un appartement contenant $C^-$ et $[z,y)$
(5.1), on note $C^+$ la chambre locale en
$z$ oppos\'ee \`a $C^-$ et $w'_+=d_{*z}^{nr}(C,C^+)$. Il est clair que
$\qp_+\subset\qr({\overline{C^+}})=w'_+.germ_z(z+{\overline{C^v}})\subset
w'_+.(z+{\overline{C^v}})$, donc
$w_+\le w'_+$ [Humphreys-90; 5.12 p. 123]. Il est clair que $\qr$ induit dans
l'immeuble jumel\'e
$\r I_z$ (avec sa structure non restreinte) la r\'etraction de centre $C_z$. De
plus $C_z^+$ et $C_z^-$ sont oppos\'ees et se r\'etractent respectivement sur
$-w'_+C_z$ et $w_-C_z$. Enfin les choix de syst\`emes de g\'en\'erateurs de $W^v=w^v(A)$
effectu\'es en 5.7 et 5.8 sont les m\^emes. D'apr\`es 5.7 on a $w'_+\le w_-$, donc
$w_+\le w_-$.
\qed

\par \noindent
{\bf Th\'eor\`eme 5.9.} {\sl Dans la masure ordonn\'ee \r I la relation $\le $ ou $\bc\le $
(d\'efinie en 5.4.1) est un pr\'eordre. Plus pr\'ecis\'ement si trois points $x$, $y$ et
$z$ dans \r I sont tels que $x\le y$ et $y\le z$, alors $x\le z$ et en particulier les points $x$,
$z$ sont dans un m\^eme appartement; de m\^eme pour $\bc\le $.
 }

\parni {\bf  Remarque.} D'apr\`es 1.5.4 la relation $\bc\le $  est une relation
d'ordre en dehors du cas fini (o\`u elle est triviale). On a $x\le y$, $y\bc\le z$, 
$y\ne z$ (ou $x\bc\le y$, $y\le z$, $x\ne y$) $\Rightarrow x\bc\le z$.

\parni {\bf  D\'emonstration.} On suppose $x\ne y$ et $y\ne z$. Dans le cas semi-discret, il
suffit de traduire mot \`a mot la d\'emonstration du Theorem 6.9 de
[Gaussent-Rousseau-08] en rajoutant au dictionnaire les traductions suivantes:
[\lc; sect. 4.3.4 ou lemma 6.11]
$\mapsto$ 5.1 et 5.5, [\lc; sect. 4.4] $\mapsto$ 2.7.1 et [\lc; prop. 6.1]
$\mapsto$ 5.8. Pour le cas g\'en\'eral il faut r\'eorganiser les ingr\'edients de cette
d\'emonstration:

\par 1) Pour $z'\in[y,z[$ tel que $x\leq z'$, on choisit un appartement $A$ 
contenant $[z',x)$ et $[z',z)$ (5.1); cet appartement a un syst\`eme de racines
r\'eelles associ\'e $\qF(A)$ et on d\'efinit l'ensemble fini $\qF(z')$ des
racines $\qa\in\qF(A)$ telles que $\qa(z')>\qa(x_1)$ et $\qa(z')>\qa(z_1)$ pour 
certains $x_1\in[x,z']\cap A$ et $z_1\in[z,z']\cap A$. Comme $[z',x)$ et 
$[z',z)$ sont pr\'eordonn\'es, 5.5 montre que $\qF(z')$  d\'epend, \`a isomorphisme pr\`es,
seulement de $[z',x)$ et $[z',z)$ mais pas de $A$.

\par On raisonne par r\'ecurrence sur $\vert\qF(y)\vert$; s'il vaut $-1$ le
th\'eor\`eme est d\'emontr\'e.

\par 2) D'apr\`es 5.1, il existe un appartement $A_1$ contenant $x$ et $[y,z)$. 
On choisit une chambre vectorielle $C^v$ dans $A_1$ telle que son syst\`eme de
racines positives associ\'e $\qF^+(C^v)$ contienne les racines $\qa\in\qF(A_1)$
telles que $\qa(y)>\qa(x)$ ou $\qa(y)=\qa(x)$ et $\qa(z_1)>\qa(y)$ (pour un  
$z_1\in[y,z]\cap A_1$); en particulier $[x,y]\subset (y-\overline{C^v})\cap
(x+\overline{C^v})$ et $[y,z)\subset x+\overline{C^v}$. Maintenant si
$\qa\in\qF^+(C^v)$ est tel que 
$\qa(z_1)<\qa(y)$  (pour un $z_1\in[y,z]\cap A_1$) alors $\qa(y)>\qa(x)$; donc
$\qF(y)$ (calcul\'e dans $A_1$) est l'ensemble des racines $\qa\in\qF^+(C^v)$
telles que $\qa(z_1)<\qa(y)$  (pour un $z_1\in[y,z]\cap A_1$).

\par Soient \g Q le germe de quartier associ\'e \`a $-C^v$ dans $A_1$ et $\qr$ la 
r\'etraction de centre \g Q sur $A_1$.

\par 3) On note $\qS_1=[y,Z]$ le segment de $A_1$ d'origine $y$, contenant
$[y,z)$ et affinement isomorphe \`a $[y,z]$. On note $Z_1$ le point le plus \'eloign\'e
de $y$ dans $\qS_1\cap(x+\overline{C^v})$, donc $\forall Z_2\in[y,Z_1]$,
$x\in Z_2-\overline{C^v}$. On note enfin $z_1\in]y,z]$ le point correspondant \`a
$Z_1$ par l'isomorphisme affine de $[y,Z]$ sur $[y,z]$.

\par Comme en 2.7.1 on obtient une suite $y_0=y,\,y_1,\dots,\,y_n=z_1\in[y,z_1]$
et des  appartements $A_1,A_2,\dots,A_n$ tels que $A_i$ contienne \g Q et
$[y_{i-1},y_i]$. En particulier $x\in y_1-\overline{C^v}$ et $x\le y_1$ (resp.
$x\bc\le y_1$): le th\'eor\`eme est d\'emontr\'e si $y_1=z$. Si $y_1\ne z$ on raisonne par
r\'ecurrence sur n. On consid\`ere l'appartement $A_2$ ci-dessus si $n\ge 2$ et un
appartement $A_2$ contenant $[y_1,z)$ et $y_1-\overline{C^v}$ (donc aussi $x$ et
\g Q) si $n=1$ \ie $y_1=z_1\ne z$. On va maintenant comparer $\qF(y)$ et $\qF(y_1)$.

\par 4) La r\'etraction $\qr$ envoie isomorphiquement $A_2$ sur $A_1$.  Ceci permet
d'identifier $\qF(y_1)$ avec l'ensemble $\qF'(y_1)$ des racines $\qa\in\qF(A_1)$
telles que $\qa(y_1)>\qa(x)$ (donc $\qa\in\qF^+(C^v)$) et 
$\qa(y_1)>\qa(\qr z_2)$ (pour un $z_2\in [y_1,z]\cap A_2$). On a
$\qr([y_1,z))=y_1+[0,1)w^+\ql\,$, $[y,y_1)=y+[0,1)w^-\ql$ pour un certain 
$\ql\in{\overline{C^v}}$  et certains $w^+,w^-\in W^v$, avec $w^+\leq w^-$ d'apr\`es
5.8. En particulier, pour $\qa\in\qF^+(C^v)$, $\qa(y_1)>\qa(\qr z_2)$ signifie
$\qa(w^+\ql)<0$, donc 
$\qF'(y_1)\subset\{\qa\in\qF^+(C^v)\mid\qa(w^+\ql)<0\}$ et (comme $w^+$ est
choisi minimal) cet ensemble a pour cardinal $\ell(w^+)$. Mais nous avons vu en
2) que $\qF(y)=\{\qa\in\qF^+(C^v)\mid\qa(w^-\ql)<0\}$. Donc, comme $w^+\leq w^-$, 
$\vert\qF'(y_1)\vert\leq\ell(w^+)\leq\ell(w^-)\leq\vert \qF(y)\vert$ .

\par Si $\vert\qF'(y_1)\vert<\vert\qF(y)\vert$ le th\'eor\`eme est d\'emontr\'e par
r\'ecurrence. Sinon les 4 nombres ci-dessus sont \'egaux; en particulier, comme  
$w^+\leq w^-$, on a $w^+=w^-$ et $\qF'(y_1)=\qF(y)$.

\par 5) Si $y_1=z_1\ne z$, alors, par d\'efinition, il existe $\qa\in\qF^+(C^v)$ tel
que $\qa(Z)<\qa(y_1)=\qa(x)$, donc $\qa(Z)<\qa(y)$. Ainsi $\qa\in\qF(y)$ et
$\qa\notin\qF'(y_1)$; c'est absurde.

\par Donc maintenant $y_1\ne z_1$. Comme $w^+=w^-$, $y$ et $\qr([y_1,z))$ sont
align\'es dans $A_1$. L'isomorphisme $\qr$ de $A_2$ sur $A_1$ identifie la chambre
vectorielle $C^v$ de $A_1$ d\'efinie en 2) \`a celle que l'on peut d\'efinir de la m\^eme
mani\`ere dans $A_2$ avec $y_1$ (en effet, comme $y_1\ne z_1$,
$\qa(y)>\qa(x)\Rightarrow\qa(y_1)>\qa(x)$).

\par Si on d\'efinit $\qS'_1=[y_1,Z']$, $Z'_1\in\qS'_1$ et $z'_1\in[y_1,z]$ comme
en 3) avec $(A_2,y_1,C^v)$ \`a la place de $(A_1,y,C^v)$, l'isomorphisme $\qr$ de
$A_2$ sur $A_1$ envoie $\qS'_1$ sur $[y_1,Z]$ et $Z'_1$ sur $Z_1$; donc
$z'_1=z_1$. Le raisonnement de l'\'etape 3) (pour $y_1$ et $A_2$) fournit alors les
m\^emes points $y_1,\dots,y_n=z_1\in[y_1,z_1]$. On est donc pass\'e de $n$ \`a $n-1$ et
le th\'eor\`eme est prouv\'e par r\'ecurrence.  \qed

\par \noindent
{\bf Remarque 5.10.} 

\par Ce th\'eor\`eme est le meilleur r\'esultat de cet article pour caract\'eriser des
paires de points de la masure ordonn\'ee \r I qui sont situ\'es dans un m\^eme
appartement.

\par Dans le cas fini, on a $x\le y$ (et m\^eme $x\bc\le y$) pour tous $x$,
$y$ dans un m\^eme appartement $A$. Pour tous $x,y\in\r I$ il existe un point z et
deux appartements $A_1,A_2$ tels que $x,z\in A_1$ et $y,z\in A_2$; on a donc
$x\le z\le y$ et $x,y$ sont dans un m\^eme appartement d'apr\`es 5.9. D'apr\`es 5.1 et 5.5
deux facettes de \r I sont donc toujours contenues dans un m\^eme appartement $A$,
qui est unique \`a un isomorphisme fixant les deux facettes pr\`es. Ainsi 5.9 peut 
\^etre vu comme une g\'en\'eralisation de [Bruhat-Tits-72; 7.3.4 et 7.3.6].

\par Si le groupe de Weyl $W^v$ de la masure ordonn\'ee \r I est affine, il existe
une forme lin\'eaire $\qd$ sur $V$ (la "plus petite racine imaginaire positive")
qui est combinaison lin\'eaire \`a coefficients positifs des $\qa_i$, invariante par
$W^v$ et telle que $\r T^\circ=\{v\in V\mid\qd(v)>0\}$, $\r T=\r
T^\circ\cup V_0$. Chaque appartement $A$ est muni d'une forme lin\'eaire $\qd_A$
telle que $x\le _Ay\Leftrightarrow \qd_A(y-x)>0$ ou $x-y\in V_0$; de plus $\qd_A(y-x)$
est ind\'ependant du choix de l'appartement $A$ contenant $x$ et $y$. Choisissons un
germe de quartier positif \g Q, on peut donc d\'efinir $\qd_{\sg Q}(y-x)$ pour tous
$x$, $y$ dans \r I par $\qd_{\sg Q}(y-x)=\qd_A(\qr y-\qr x)$ pour tout
appartement $A$ contenant \g Q et pour $\qr=\qr_{A,\sg Q}$. D'apr\`es 2.7.1 on a
$\qd_{\sg Q}(y-x)=\qd_A(y-x)$ d\`es que $x,y\in A$ et $\qd_A(y-x)\ne 0$. Donc $y\ge x$ et
$y\ne x\Rightarrow\qd_{\sg Q}(y-x)>0$.

\par Howard Garland [1995] consid\`ere la situation d'un groupe de lacets sur un
corps discr\`etement valu\'e $K$ (on supposera ici que son corps r\'esiduel contient
$\C$). Plus pr\'ecis\'ement, pour un groupe semi-simple simplement lac\'e $G$, il
\'etudie un sous-groupe ${\widehat G}_b$ de $G(K((t))\,)$ qui est contenu dans le
groupe engendr\'e par $U_{\sg Q}^{max+}$ et $U_{\sg Q}^-$ (notations de
[Gaussent-Rousseau-08; sect. 3.2 et 3.3] en consid\'erant le groupe de Kac-Moody
$G_1=G(K[t,t^{-1}])\rtimes K^*$). Il semble que ${\widehat G}_b$ n'agisse que sur
un espace $\r I'$ plus gros que la masure \r I construite dans \lc pour $G_1$,
mais il fixe \g Q. Dans $\r I'$ la d\'ecomposition de Cartan tordue que Garland
prouve semble montrer que, si $x,y\in\r I'$ et $\qd_{\sg Q}(y-x)>0$, il existe un
appartement contenant $x,y$ et \g Q. On peut donc se demander si la r\'eciproque de
la derni\`ere phrase de l'alin\'ea pr\'ec\'edent est vraie dans toute masure ordonn\'ee de
groupe de Weyl $W^v$ affine. Un tel r\'esultat pr\'eciserait beaucoup le th\'eor\`eme
5.9: presque toutes les paires de points seraient contenues dans un m\^eme
appartement.

\bigskip
\parni {\S$\,${\bf 6.}$\quad${\bftwelve Masures associ\'ees \`a certains groupes de
Kac-Moody }}
\bigskip
\par Dans ce dernier paragraphe on se propose de montrer que les masures (hovels)
construites dans [Gaussent-Rousseau-08] (ou en tout cas la plupart
d'entre elles) forment des masures affines, ordonn\'ees, \'epaisses et
semi-discr\`etes au sens adopt\'e ici. On ne rappellera pas toutes les d\'efinitions
introduites dans
\LC.

\bigskip \noindent
{\bf  6.1. La situation}

\par On consid\`ere donc un groupe de Kac-Moody $G$ d\'efini sur un corps $K$ muni
d'une valuation discr\`ete $\qo$ pour laquelle le corps r\'esiduel $\qk$ contient
$\C$.

\par On suppose que $G$ v\'erifie les hypoth\`eses techniques de \LC; en particulier
il est d\'eploy\'e et sym\'etrisable. On suppose de plus qu'il v\'erifie l'hypoth\`ese (SC)
de [\lc; sect. 2.1.5]: 
\qquad\qquad(SC)\qquad $\tch{Q}{\;}=\sum_{i\in I}\;\Z\tch{\alpha}{_i}$ est
sans cotorsion dans $Y$.

\parni En fait, d'apr\`es le raisonnement de \LC, on doit pouvoir toujours se
ramener \`a ce cas sans beaucoup changer l'espace \r I.

\par On a construit dans \LC un ensemble \r I, muni d'un recouvrement par un
ensemble \r A de sous-ensembles appel\'es appartements. Tous ces appartements sont
permut\'es par $G(K)$ et donc isomorphes \`a un appartement fondamental $A_f=\A$. Le
stabilisateur de $A_f$ est un groupe $N(K)$ (normalisateur du tore fondamental
$T_f$ de groupe de caract\`eres (resp. cocaract\`eres) $X$ (resp. $Y$)) qui agit sur
$\A$ par un groupe
$W_Y$.

\par L'espace $\A$ est un appartement semi-discret au sens de 1.4 ci-dessus. On le
consid\`ere comme mod\'er\'ement imaginaire, car muni des murs d\'ecrits en 1.2.1 et 1.4.d.
Le groupe
$W_Y$ est un groupe de transformations affines stabilisant l'ensemble des murs
r\'eels ou imaginaires et contenant le groupe de Weyl affine $W$ (engendr\'e par les
r\'eflexions par rapport aux murs r\'eels). En particulier $W_Y$ normalise $W$. Ainsi
tout appartement de
\r A est muni d'une structure d'appartement de type $\A$ et l'axiome (MA1) est
v\'erifi\'e.

\bigskip \noindent
{\bf  6.2. Le sous-groupe} $G(K)_1$

\par Le groupe $Y_1=Y/\tch{Q}{\;}$ est ab\'elien libre. C'est le groupe des
cocaract\`eres d'un tore $T_1$ (de groupe de caract\`eres $X_1=\{\chi\in X \mid
\chi(\tch{Q}{\;})=0\} ) $ et l'homomorphisme naturel $\qd$ de $T_f$ sur $T_1$ se
prolonge en un homomorphisme $\qd$ de $G$ sur $T_1$ trivial sur tous les
sous-groupes radiciels: il suffit, pour tout corps $K'$ contenant $K$, de
consid\'erer la pr\'esentation de $G(K')$ donn\'ee par J. Tits, \cf [R\'emy-02; 8.3.3 et
8.4.2].

\par Si \r O est l'anneau des entiers de $(K,\qo)$, on a $T_1(K)/T_1(\r O)\simeq
Y_1$ et on note $G(K)_1$ le sous-groupe distingu\'e $G(K)_1=\qd^{-1}(T_1(\r O))$ de
$G(K)$, donc $G(K)/G(K)_1\simeq Y_1$. Si $T_f(K)_1=T_f(K)\cap G(K)_1$, on a 
$T_f(K)/T_f(K)_1\simeq Y_1$.

\par On sait que l'image $\qn(N(K))$ de $N(K)$ dans le groupe des automorphismes
affines de $A_f=\A$ est $W_Y=W^v\ltimes Y$, que $\qn(T_f(K))=Y$ et que
Ker$\qn=H=T_f(\r O)$. En fait 
$N(K)\cap G(K)_1$ a pour image (par $\qn$) $W=W^v\ltimes \tch{Q}{\;}$; en effet,
comme $\qd$ est trivial sur les sous-groupes radiciels de $G$, $G(K)_1$ contient
des repr\'esentants dans $N(K)$ de g\'en\'erateurs de $W^v$.

\par On sait qu'un point ou un germe de quartier a un bon fixateur
[Gaussent-Rousseau-08; sect 4.1 et 4.2.4] et qu'un point et un germe de quartier
sont toujours contenus dans un m\^eme appartement [\lc; sect 4.3.3], le lemme suivant
montre donc que
$G(K)_1$ est transitif sur \r A. Comme le stabilisateur dans 
$G(K)_1$ de $A_f$ est
$N(K)\cap G(K)_1$ qui agit sur $A_f$ via $W$, il est clair que chaque appartement
de \r I est muni d'une unique structure d'appartement de type $\A$ (au sens de
1.13) et que tout \'el\'ement de $G(K)_1$ induit un isomorphisme d'un appartement
quelconque de \r I sur son image.

\bigskip\parni
{\bf Lemme 6.3.} {\sl Soit $\qO\subset\A$ un filtre qui a un assez bon fixateur,
alors $G_\qO={\widehat P}_\qO$ est contenu dans $G(K)_1$, en particulier il induit 
un isomorphisme entre un appartement et son image. Bien sur $G_\qO$ permute
transitivement les appartements contenant $\qO$.
 }
\parni {\bf  D\'emonstration.} On a $G_\qO=(G_\qO\cap U^+).(G_\qO\cap
U^-).{\widehat N}_\qO$ (\`a l'\'echange de $+$ et $-$ pr\`es). On sait que $\qd$ est
trivial sur
$U^+$ et
$U^-$, il suffit donc de voir que ${\widehat N}_\qO\subset G(K)_1$. Mais
$\qn({\widehat N}_\qO)={\widehat N}_\qO/H$ est le fixateur de $\qO$ dans $W_Y$, on
doit montrer qu'il est dans $W$.  Soit $w=w^v.y\in W_Y=W^v\ltimes Y$ (avec $y\in Y$
et $w^v\in W^v$) qui fixe un point $x$ de $\A$; comme $w^v$ fixe $x$ \`a
$\tch{Q}{\;}\otimes\R$ pr\`es, $y$ doit \^etre dans
$Y\cap(\tch{Q}{\;}\otimes\R)=\tch{Q}{\;}$: on a bien $w\in W=W^v\ltimes
\tch{Q}{\;}$.  \qed

\parni{\bf 6.4. Groupe parabolique associ\'e \`a une facette vectorielle $F^v$ de
$\A$}

\par On associe \`a $F^v$ un ensemble parabolique de racines
$\qF(F^v)=\qF^m(F^v)\cup\qF^u(F^v)$ avec $\qF(F^v)=\{\qa\in\qF\mid\qa(F^v)\ge 0\}$,
$\qF^m(F^v)=\{\qa\in\qF\mid\qa(F^v)=0\}=\{\qa\in\qF\mid F^v\subset M^v(\qa)\}$ et
$\qF^u(F^v)=\qF(F^v)\setminus \qF^m(F^v)$. On en d\'eduit des sous-groupes de
$G(K)$:

\par- le sous-groupe parabolique $P(F^v)$ est engendr\'e par $T_f(K)$ et les
$U_\qa(K)$ pour $\qa\in \qF(F^v)$,
\par- son facteur de L\'evi $M(F^v)$ est engendr\'e par $T_f(K)$ et les
$U_\qa(K)$ pour $\qa\in \qF^m(F^v)$,
\par- son radical unipotent $U(F^v)$ est le plus petit sous-groupe normal de
$P(F^v)$ contenant les $U_\qa(K)$ pour $\qa\in \qF^u(F^v)$,
\par et on a la d\'ecomposition en produit semi-direct $P(F^v)=M(F^v)\ltimes
U(F^v)$, \cf [Rousseau-06; 1.7] ou [R\'emy-02; 6.2].

\par Si on suppose $F^v$ dans l'adh\'erence de la chambre fondamentale positive
$C^v_f$, $U(F^v)$ est aussi le plus petit sous-groupe normal de $U^+(K)$
contenant les $U_\qa(K)$ pour $\qa\in \qF^u(F^v)$ [R\'emy-02; 6.2]. On a donc 
$U(F^v)\subset U^+(K)\cap U^{max}(\qF^u(F^v))=G(K)\cap
U^{max}(\qF^u(F^v))=U^{pmax}(\qF^u(F^v))$ , \cf [Gaussent-Rousseau-08; sect.
3.3.3].

\par D'apr\`es la d\'efinition par g\'en\'erateurs et relations des groupes de Kac-Moody
[R\'emy-02; 8.3.3] $M(F^v)$ est le (ou au moins un quotient du) groupe de Kac-Moody
$G_{F^v}(K)$ de syst\`eme g\'en\'erateur de racines $(A(J),X,Y,(\alpha_i)_{i\in
J}, (\tch{\alpha}{_i})_{i\in J})$ si $F^v=F^v(J)$. Si $F^v=V_0$, on a
$M(F^v)=G(K)$ et $U(F^v)=\{1\}$. Si $F^v$ est sph\'erique, $\qF^m(F^v)$ est fini et
$M(F^v)$ est un groupe r\'eductif, \cf [R\'emy-02; 12.5.2] ou [Rousseau-06; 1.7].

\bigskip
\parni{\bf 6.5. Facette microaffine et groupe parahorique associ\'es \`a une chemin\'ee}

\par Soit $\g r=\g r(F,F^v)$ une chemin\'ee de $\A$ associ\'ee \`a une facette
$F=F(x,F^v_0)$ et une facette vectorielle $F^v$. 

\par \`A $F^v$, on associe l'ensemble $\r M^{\qF^m(F^v)}$ des murs de \r M de
direction Ker$(\qa)$ pour $\qa\in\qF^m(F^v)$ et l'appartement affine 
$\A^{F^v}$ de $G_{F^v}(K)$ correspondant au choix dans $\A$ des murs de $\r
M^{\qF^m(F^v)}$ et des murs imaginaires de direction contenant $F^v$. \`A la
chemin\'ee \g r, on associe la facette microaffine $F^\qm=F^\qm(\g r)=(\r F,F^v)$ o\`u 
$\r F=F_{\qF^m(F^v)}(x,F^v_0)$ est la facette engendr\'ee par $F$ dans $\A^{F^v}$.
Ainsi $\overline{\r F}$ est le filtre des parties de $\A$ contenant une intersection
de demi-espaces $D(\qa,k)$ de $\A^{F^v}$ (\ie pour $\qa\in\qD$, $\qa(F^v)=0$ et
$k\in\qG_\qa$) contenant $F$ ; en particulier $\overline{\r F}\supset \g r$.

\par Le sous-groupe parahorique $P^\qm(\g r)=P(F^\qm)$ est le produit semi-direct 
 $P(F^\qm)=M(F,F^v)\ltimes U(F^v)$ o\`u $M(F,F^v)$ est le (ou l'image dans $M(F^v)$
du) sous-groupe ${\widehat P}_{\r F}$ de $G_{F^v}(K)$ associ\'e \`a la facette \r F
dans [Gaussent-Rousseau-08].

\parni{\bf Remarques.} 1) Si $F^v=V_0$, alors $\r F=F$ et
$P(F^\qm)=M(F,F^v)={\widehat P}_F$.

\par 2) $F^\qm$ n'est une vraie facette microaffine (au sens de [Rousseau-06;
2.5]) que si $F^v$ est sph\'erique \ie \g r \'evas\'ee. Dans ce cas \g F est une
facette de l'appartement de Bruhat-Tits $\A^{F^v}$ du groupe r\'eductif $M(F^v)$;
comme nous avons suppos\'e la valuation discr\`ete, c'est m\^eme un sous-ensemble de
$\A^{F^v}$ (un polysimplexe).

\par 3)  En fait \r F, $F^\qm$, $M(F,F^v)$ et $P^\qm(\g r)$ ne d\'ependent que du
germe \g R de  \g r.

\par 4) Le groupe $P^\qm(\g r)$ n'est parahorique qu'en un sens g\'en\'eralis\'e, celui
de [Rousseau-06; 2.10] si \g r est \'evas\'ee. En fait $M(F,F^v)$ est un sous-groupe
parahorique de $M(F^v)$ (\cf [Gaussent-Rousseau-08; sect. 3.7]), au sens
classique si \g r est \'evas\'ee.

\bigskip\parni
{\bf Proposition 6.6.} {\sl Le groupe $P^\qm(\g r)$ fixe (point par point) le
germe de chemin\'ee \g R.
 }
\parni{\bf N.B.} Il n'est pas sur que $U(F^v)=U^+(K)\cap U^{max}(\qF^u(F^v))$ et
pas plus \'evident que $U^+(K)\cap U^{max}(\qF^u(F^v))$ fixe \g R.

\parni {\bf  D\'emonstration.}
\par Supposons $F^v$ dans l'adh\'erence de la chambre fondamentale $C^v_f$; alors un
\'el\'ement $u_0$ de $U(F^v)$ est un produit fini de conjugu\'es
$v_iu_iv_i^{-1}$ avec $v_i\in U^+(K)=U(C^v_f)$, $u_i\in U_{\qa_i}(K)$,
$\qa_i\in\qF^u(F^v)\subset\qF^+$; de m\^eme les $v_i$ sont des produits finis
d'\'el\'ements $u_{i,j}\in U_{\qb_{i,j}}(K)$ avec $\qb_{i,j}\in\qF^+$. Donc il existe
$y\in \A$ tel que $u_i\in U_{\qa_i,y}$ et $u_{i,j}\in U_{\qb_{i,j},y}$, $\forall
i,j$. Ainsi $u_0\in U^{pm}_y(\qF^u(F^v))$ d'apr\`es les relations de commutation
dans le groupe compl\'et\'e $U^{max}(\qF^+)$. Il est alors clair que $u_0$ fixe
(point par point) $\g r(F+\qx,F^v)$ pour un bon choix de $\qx\in F^v$. Par
ailleurs $M(F,F^v)$ fixe la facette-ferm\'ee $\overline{\r F}$ qui contient 
$\g r(F,F^v)$ et donc $\g r(F+\qx,F^v)$.  \qed

\parni
{\bf Proposition 6.7.} {\sl Soient $\g R_1=\g R(F_1,F_1^v)$ et $\g R_2=\g
R(F_2,F_2^v)$ deux germes de chemin\'ees de $\A$ avec $\g R_1$ \'evas\'e, alors \qquad
$G(K)=P^\qm(\g R_1).N(K).P^\qm(\g R_2)$.
 }
\parni{\bf Remarques.} 1) Ceci constitue une d\'ecomposition de Bruhat ou
Birkhoff mixte, \cf 2.1 et [Rousseau-06; 3.5] ou m\^eme une d\'ecomposition d'Iwasawa
mixte quand $\g R_2$ est d\'eg\'en\'er\'e en particulier une facette ferm\'ee.

\par 2) D'apr\`es la d\'emonstration ci-dessous, il suffit en fait de supposer que,
$\forall w\in W^v$ $\qF^m(F^v_1)\cap w\qF^m(F^v_2)$ est fini (au lieu de $\g R_1$
\'evas\'e \ie $\qF^m(F_1^v)$ fini).

\par 3) On d\'eveloppe ici la remarque 4.6 de [Gaussent-Rousseau-08].

\parni {\bf  D\'emonstration.}

\par D'apr\`es la d\'ecomposition de Bruhat ou de Birkhoff de $G(K)$, pour tout $g\in
G(K)$ il existe $n\in N(K)$ tel que $g\in P(F_1^v)nM(F_2^v)U(F_2^v)$. Mais
$\qF(n^{-1}F_1^v)\cap\qF^m(F_2^v)$ est un syst\`eme parabolique de racines dans
$\qF^m(F_2^v)$ (associ\'e \`a la facette $cl_{\qF^m(F_2^v)}(n^{-1}F_1^v)$ ). La
d\'ecomposition d'Iwasawa de $M(F_2^v)$ [Gaussent-Rousseau-08; prop. 3.6] donne
donc: $M(F_2^v)=(M(F_2^v)\cap P(n^{-1}F_1^v)).(N(K)\cap M(F_2^v)).M(F_2,F_2^v)$.
Ainsi $g\in P(F_1^v)N(K)M(F_2,F_2^v)U(F_2^v)$ et il existe $m\in M(F_1^v)$,
$n_1\in N(K)$ tels que $g\in U(F_1^v)mn_1P^\qm(\g R_2)=U(F_1^v)mP^\qm(n_1\g
R_2)n_1$.

\par Mais $\qF(n_1F_2^v)\cap\qF^m(F_1^v)$ est un syst\`eme parabolique de racines
dans $\qF^m(F_1^v)$ (associ\'e \`a la facette $cl_{\qF^m(F_1^v)}(n_1F_2^v)$ ).  La
d\'ecomposition d'Iwasawa de $M(F_1^v)$ donne: $M(F_1^v)=M(F_1,F_1^v).(N(K)\cap
M(F_1^v)).M(n_1F_2^v,F_1^v).(U(n_1F_2^v)\cap M(F_1^v))$, o\`u on note
$M(n_1F_2^v,F_1^v)$  le sous-groupe de $M(F_1^v)$ associ\'e au
sous-syst\`eme de racines
$\qF^m(n_1F_2^v)\cap \qF^m(F_1^v)$ (par hypoth\`ese ce  syst\`eme de racines
est fini, donc $M(n_1F_2^v,F_1^v)$ est un groupe r\'eductif). Il existe donc
$n_2\in N(K)\cap M(F_1^v)$ tel que: 
\par\qquad $m\in n_2M(n_2^{-1}F_1,n_2^{-1}F_1^v=F_1^v).M(n_1F_2^v,F_1^v).
(U(n_1F_2^v)\cap M(F_1^v))$.

\par La facette $n_2^{-1}F_1$ engendre une facette $F_1'$ de $\A$ pour
$\qF^m(n_1F_2^v)\cap \qF^m(F_1^v)$ \ie pour $M(n_1F_2^v,F_1^v)$. Le sous-groupe
$M(F'_1)$ correspondant de $M(n_1F_2^v,F_1^v)$ v\'erifie donc $M(F'_1)\subset
M(n_2^{-1}F_1,F_1^v)$. De m\^eme la facette $n_1F_2$ engendre une facette $F'_2$ de
$\A$ pour $\qF^m(n_1F_2^v)\cap \qF^m(F_1^v)$ \ie pour $M(n_1F_2^v,F_1^v)$. Le 
sous-groupe $M(F'_2)$ correspondant de $M(n_1F_2^v,F_1^v)$ v\'erifie donc
$M(F'_2)\subset M(n_1F_2,n_1F_2^v)$.

\par La d\'ecomposition de Bruhat-Iwahori classique dans le groupe r\'eductif
$M(n_1F_2^v,F_1^v)$ (contenu dans $M(F_1^v)$) donne donc:

\par  $M(n_1F_2^v,F_1^v)\subset M(F'_1).(N(K)\cap
M(n_1F_2^v,F_1^v)).M(F'_2)$
\par\qquad\qquad\qquad $\subset M(n_2^{-1}F_1,F_1^v).N(K).M(n_1F_2,n_1F_2^v)$.

\parni Ainsi $m\in
n_2M(n_2^{-1}F_1,F_1^v).N(K).M(n_1F_2,n_1F_2^v).U(n_1F_2^v)$ 
\par\qquad\qquad$=M(F_1,F_1^v).N(K).P^\qm(n_1\g R_2)$,

\parni et $g\in U(F_1^v).M(F_1,F_1^v).N(K).P^\qm(n_1\g R_2)n_1=
P^\qm(\g R_1).N(K).P^\qm(\g R_2)$.   \qed

\parni
{\bf Proposition 6.8.} {\sl Soit $\g r=\g r(F,F^v)$ une chemin\'ee solide de $\A$.
Alors \g r et son germe \g R ont de bons fixateurs. Plus pr\'ecis\'ement si \g r est
\'evas\'ee et $F^v$ est dans l'adh\'erence de la chambre vectorielle fondamentale, ces
bons fixateurs sont $G_{\sg r}={\widehat N}_{\sg
r}.U^{nm-}_{F+F^v}.U^{pm+}_F=M(F,F^v).U^{pm}_F(\qF^u(F^v))$ et $G_{\sg
R}=M(F,F^v).U(F^v)=P^\qm(\g R)$
 }
\parni{\bf Remarque.} Une face de quartier sph\'erique $\g f=x+F^v$ et son germe (\`a
l'infini) \g F ont aussi de bons fixateurs. Il suffit, dans la d\'emonstration
ci-dessous, de remplacer $F$ par $\{x\}$ et $cl(F\cup (F+n\qx))$ par
$cl(\{x,x+n\qx\})\cap supp(x+F^v)$.

\parni {\bf  D\'emonstration.}

\par On choisit $\qx\in F^v$ et une chambre vectorielle ferm\'ee ${\overline{C^v}}$
contenant $-F^v$. Pour $n\in\N$ les facettes $F=F(x,F_0^v)$ et
$F_n=F(x+n\qx,F_0^v)$ ont un bon fixateur [Gaussent-Rousseau-08; sect. 4.2.2]; de
plus $F\subset F_n+{\overline{C^v}}$ et $F_n\subset F-{\overline{C^v}}$. D'apr\`es
[\lc; rem. 4.4 et prop. 4.3 1)], $F\cup F_n$ et $cl(F\cup F_n)$ ont de bons
fixateurs. Mais pour $n>0$, le fixateur dans $W^e_{\sR}$ de $cl(F\cup F_n)$ est
fini (puisque \g r est solide), donc $\g r = \cup_n\;cl(F\cup F_n)$ et $\g
R=germ_\infty(\g r)$ ont de bons fixateurs [\lc; prop. 4.3 3) et 2)].

\par Il reste \`a identifier ces bons fixateurs quand \g r est \'evas\'ee. On a  $G_{\sg
r}={\widehat N}_{\sg r}.U^{nm-}_{F+F^v}.U^{pm+}_F$. Mais 
$U_F^{max}(\qF^+)=U_F^{max}(\qF^{m+}(F^v))\ltimes U_F^{max}(\qF^u(F^v))$ et, par
hypoth\`ese, $\qF^{m}(F^v)$ est fini, donc $U^{pm+}_F=V^+_F\ltimes
U_F^{pm}(\qF^u(F^v))$ o\`u $V^+_F$ est le groupe $U^+_F$ (d\'efini dans
[\lc; sect. 3.2]) relatif au groupe r\'eductif $M(F^v)$. De plus
$U^{nm-}_{F+F^v}$ est le groupe $V^-_F$ (\ie le $U^-_F$ relatif \`a ce groupe
$M(F^v)$). Enfin ${\widehat N}_{\sg r}$ est le ${\widehat N}_{F}$ relatif \`a ce 
groupe $M(F^v)$. Ainsi $G_{\sg r}={\widehat N}_{\sg
r}.V^-_F.V^+_F.U_F^{pm}(\qF^u(F^v))=M(F,F^v).U_F^{pm}(\qF^u(F^v))$. Le r\'esultat
pour $G_{\sg R}$ est maintenant clair, d'apr\`es 6.6.  \qed

\parni{\bf 6.9. Cons\'equences}

\par 1) L'axiome (MA2) est v\'erifi\'e: d'apr\`es [Gaussent-Rousseau-08; sect 4.1 et
4.2] un point, un germe d'intervalle pr\'eordonn\'e (\cf sect. 4.4 et prop. 4.3 1)
et 2) de \lc) et une demi-droite g\'en\'erique  ont de bons fixateurs. On vient de voir
en 6.8 que c'est aussi vrai pour une chemin\'ee solide. D'apr\`es [\lc; remark
4.2 et prop. 4.3.1] un sous-groupe de ce bon fixateur fixe l'enclos de cette
partie $F$ et est transitif sur les appartements contenant $F$. D'apr\`es le lemme
6.3, $G_F$ induit des isomorphismes entre ces appartements.

\par 2) L'axiome (MA3) est v\'erifi\'e: on a $G(K)=G_{\sg R}.N(K).G_F$ d'apr\`es 6.6 et
6.7. Vue la transitivit\'e de $G_{\sg R}$ ou $G_F$ sur les appartements contenant
\g R ou $F$ (si $F$ est une facette ou une chemin\'ee solide) un raisonnement
classique prouve l'axiome (MA3).

\medskip
\par Il reste \`a v\'erifier l'axiome (MA4).

\medskip\parni
{\bf Proposition 6.10.} {\sl Soient $\g R_1=\g R(F_1,F^v_1)$ et $\g R_2=\g
R(F_2,F^v_2)$ deux germes de chemin\'ees de $\A$ et $\qO=\g R_1\cup\g R_2$. On
suppose que $\g R_1$ est \'evas\'e et que $\g R_2$ a un bon fixateur (ce qui est en
particulier v\'erifi\'e si $\g R_2$ est solide ou une facette). Alors $\qO$ a un assez
bon fixateur, le groupe
$G_\qO$ fixe l'enclos de $\qO$ et agit transitivement sur les appartements
contenant $\qO$.
\par Si de plus $\g R_2$ est \'evas\'e, alors $\qO$ a un bon fixateur.
 }

\parni {\bf  D\'emonstration.}

\par On sait qu'un germe de chemin\'ee solide ou une facette a un bon fixateur (6.8
et [Gaussent-Rousseau-08; sect. 4.2.2]).

\par On note ${\overline{C_1^v}},\dots,{\overline{C_n^v}}$ les chambres
vectorielles ferm\'ees contenant $F_1^v$ (en nombre fini car $F_1^v$ est
sph\'erique). La r\'eunion
$\Theta_1={\overline{C_1^v}}\cup\dots\cup{\overline{C_n^v}}$ contient un c\^one
ouvert contenant $F_1^v$. Mais ${\overline{F_1}}$ est born\'e (\`a $V_0$ pr\`es) donc,
pour $x$ quelconque et $\qx\in F_1^v$ assez grand, on a: 
${\overline{F_1}}+\qx\subset x+\Theta_1$ et $\g r(F_1,F_1^v)+\qx\subset
x+\Theta_1$. Ainsi $\g R_1\subset \g R_2+\Theta_1$. D'apr\`es la proposition 4.3 4) de
\LC, $\qO$ a un assez bon fixateur et  le groupe
$G_\qO$ a les propri\'et\'es annonc\'ees.

\par Si de plus $\g R_2$ est \'evas\'e, on note $\Theta_2$ l'ensemble analogue \`a
$\Theta_1$ associ\'e \`a $F_2^v$. On a $\g R_1\subset \g R_2+\Theta_1$ et $\g
R_2\subset \g R_1+\Theta_2$. Mais $\g r(F_2,F_2^v)-\Theta_2=\A$, donc $\g
R_1\subset \g R_2-\Theta_2$. Ainsi la remarque 4.4 de \LC montre (en
diff\'erenciant les cas o\`u $F_1^v$ et $F_2^v$ sont de signes contraires ou oppos\'es)
que $\qO$ a un bon fixateur. \qed

\parni
{\bf Th\'eor\`eme 6.11.} {\sl Sous la condition (SC) de 6.1, les masures construites
dans \LC sont des masures affines, ordonn\'ees, \'epaisses et semi-discr\`etes (au sens 
de cet article).
 }

\parni {\bf  D\'emonstration.}

\par On a d\'ej\`a vu l'axiome (MA1) en 6.2 et les axiomes (MA2), (MA3) en 6.9. Pour
(MA4), si deux appartements $A$ et $A_f=\A$ contiennent \g R et $F$, alors $A$
est transform\'e de $\A$ par un \'el\'ement $g$ de $G_{\sg R\cup F}$ qui fixe
$cl{_{\sA}}(\g R\cup F)$. Ainsi $A\cap\A$ contient $cl{_{\sA}}(\g R\cup F)$. De
plus
$g\in G_{\sg R}$; d'apr\`es le lemme 6.3, $g$ induit donc un isomorphisme de $\A$
sur $A$ fixant $cl{_{\sA}}(\g R\cup F)$.

\par Pour l'axiome (MAO), si $x\le y$ dans $\A$, on sait que $\{x,y\}$ a un bon
fixateur et $G_{\{x,y\}}$ fixe $[x,y]{_{\sA}}$ [\lc; sect. 4.2.1]. Donc, pour
tout appartement $A$ contenant $x$ et $y$, on a: $[x,y]_A=[x,y]{_{\sA}}$.

\par On a indiqu\'e en 6.1 que $\A$ est semi-discret.

\par Si $C$ est une chambre de $\A$ de mur $M(\qa,k)$ avec $C\subset D(\qa,k)$,
on a vu en [\lc; sect. 4.3.4] que $U_{\qa,k}$ agit transitivement sur les
chambres de \r I adjacentes \`a la chambre $C$ le long de ${\overline C}\cap
M(\qa,k)$. De plus $U_{\qa,k+1}$ est le fixateur dans $U_{\qa,k}$ de celle de ces
chambres qui est contenue dans $\A$. Enfin par construction
$U_{\qa,k}/U_{\qa,k+1}$ est un groupe isomorphe \`a $(\qk,+)$. Ainsi la cloison 
${\overline C}\cap M(\qa,k)$ est domin\'ee par une infinit\'e de chambres. \qed

\bigskip\parni
{\bf Remarque 6.12.} Les masures construites dans \LC ont des propri\'et\'es que ne
poss\`edent pas forc\'ement toutes les masures affines d\'efinies ici. Par exemple on a
montr\'e en [\lc; sect 4.3.3], comme cons\'equence de la d\'ecomposition d'Iwasawa, qu'un
germe de quartier et un germe de segment sont toujours dans un m\^eme appartement,
m\^eme si ce germe de segment n'est pas pr\'eordonn\'e.

\vskip 1.8cm \centerline{\bftwelve R\'ef\'erences.}

\bigskip

{\parindent=-1cm \leftskip= 1cm

\par  Peter {\petcap Abramenko}
\ref{96}{\it Twin buildings and applications to S-arithmetic groups}, 
 Lecture notes in Math. {\bf1641 }(Springer, Berlin, 1996). 

\par  Peter {\petcap Abramenko} et Kenneth S. {\petcap Brown }
\ref{08}{\it Buildings: theory and applications}, graduate texts in Math. {\bf ?}
(Springer, Berlin, 2008).


\par Nicole {\petcap  Bardy}
\ref{96}{\it Syst\`emes de racine infinis}, M\'emoire {\bf65} (Soc. Math.
France, 1996).

\par Nicolas {\petcap Bourbaki}
\ref{68}{\it Groupes et alg\`ebres de Lie}, Chapitres IV, V, VI (Hermann,
Paris, 1968).



\par Paul {\petcap Broussous} 
\ref{04}Simplicial complexes lying equivariantly over the affine building of $GL(N)$,
{\it Math. Annalen} {\bf 329} (2004), 495-511.
 
\par Kenneth S. {\petcap Brown } 
\ref{89}{\it Buildings} (Springer, Berlin, 1989).  

\par  Fran{\c c}ois {\petcap Bruhat} et Jacques {\petcap Tits }

\ref{72}Groupes r\'eductifs sur un corps local I , Donn\'ees radicielles
 valu\'ees, {\it Publ. Math. I.H.\'E.S.} {\bf41} (1972), 5-251.


\par Cyril {\petcap Charignon } 
\ref{?}{\it Compactifications d'immeubles affines}, en pr\'eparation. 


\par Vinay V. {\petcap Deodhar} 
\ref{89}A note on subgroups generated by reflections in Coxeter groups, {\it Arch.
Math.} {\bf 53} (1989), 543-546.

\par Howard {\petcap Garland} 
\ref{95}A Cartan decomposition for p-adic loop groups, {\it Math. Ann.}
 {\bf 302} (1995), 151-175.

\par Paul {\petcap Garrett  } 
\ref{97}{\it Buildings and classical groups}, Chapman and Hall (1997).

\par St\'ephane {\petcap Gaussent} et Guy {\petcap Rousseau}
\ref{08}Kac-Moody groups, hovels and Littelmann paths, {\it Annales Inst.
Fourier} {\bf58} (2008), ? .

\par Jean-Yves {\petcap H\'ee}
\ref{90}Sur une courte note de V.V. Deodhar, 5 pages
 (1990) in Th\`ese d'\'etat, Universit\'e Paris Sud Orsay (1993).

\par James E. {\petcap  Humphreys}
\ref{90}{\it Reflection groups and Coxeter groups}, 
Cambridge studies in adv. Math. {\bf29} (Cambridge U. Press, 
Cambridge,1990).

\par Victor G. {\petcap  Kac}
\ref{90}{\it Infinite dimensional Lie algebras}, Troisi\`eme 
\'edition, (Cambridge U. Press, Cambridge,1990).


\par Bruce {\petcap Kleiner} et Bernhard {\petcap Leeb}
\ref{97}Rigidity of quasi-isometries for symmetric spaces
 and euclidean buildings, {\it Publ. Math. I.H.\'E.S.}
 {\bf 86} (1997), 115-197.


\par Erasmus {\petcap Landvogt  } 
\ref{96}{\it A compactification of the Bruhat-Tits building},
 Springer lecture note in Math. {\bf1619} (1996).

\par Robert {\petcap Moody } et Arturo {\petcap Pianzola}
\ref{95}{\it Lie algebras with triangular decompositions},
 Wiley-Interscience (1995).

\par Bernhard {\petcap M\H uhlherr}
\ref{99}Locally split and locally finite twin buildings of 2-spherical
 type,  {\it J. reine angew. Math.} {\bf 511} (1999), 119-143.

\ref{02}Twin buildings, in {\it Tits buildings and the model theory of 
groups, W\H urzburg (2000)}, K. Tent \'editrice, London Math. Soc. 
lecture note {\bf 291} Cambridge U. Press (2002), 103-117.

\par Bernhard {\petcap M\H uhlherr} et Mark {\petcap Ronan}
\ref{95}Local to global structure in twin buildings, {\it Inventiones 
Math.} {\bf 122} (1995), 71-81.

\par Anne {\petcap Parreau  } 
\ref{00}Immeubles affines: construction par les normes et \'etude 
des isom\'etries, {\it Contemporary Math.} {\bf 262 } (2000),
 263-302.


\goodbreak
\par Bertrand {\petcap R\'emy  }\nobreak
\ref{02}{\it Groupes de Kac-Moody d\'eploy\'es et presque 
d\'eploy\'es}, Ast\'erisque {\bf277} (2002). 

\par Mark A. {\petcap Ronan}
\ref{89}{\it Lectures on buildings}, Perspectives in Math. {\bf 7} 
Academic Press (1989). 

\goodbreak
\par Guy {\petcap Rousseau  }\nobreak 
\ref{77}{\it Immeubles des groupes r\'eductifs sur les corps 
locaux},  Th\`ese d'\'etat, Universit\'e Paris-Sud Orsay (1977).





\ref{01}Exercices m\'etriques immobiliers, {\it Indag. Math.}
 {\bf 12} (2001), 383-405. 
 
\ref{06}Groupes de Kac-Moody d\'eploy\'es sur un corps local, 
immeubles microaffines {\it Compositio Mathematica} {\bf 142} (2006), 
501-528.
 
\ref{08}Euclidean buildings,  In "G\'eom\'etries 
\`a courbure n\'egative ou nulle, groupes discrets et rigidit\'e",  \'ecole d'\'et\'e,
Grenoble, 14 juin-2 juillet 2004, {\it S\'eminaires et
congr\`es} {\bf 18}, Soc. Math. France (2008), 77-116.

\ref{?}Groupes de Kac-Moody d\'eploy\'es sur un corps local 2, 
masures ordonn\'ees, en pr\'eparation. 




\par Rudolf {\petcap Scharlau  } 
\ref{95}Buildings, in {\it Handbook of incidence geometry},
 \'Editeur F. Buekenhout, Elsevier (1995), 1085-1114.

\par Jacques {\petcap Tits} 

\ref{74}{\it Buildings of spherical type and finite BN pairs},
 Springer lecture note in Math. {\bf386} (1974).





\ref{86}Immeubles de type affine, in {\it Buildings and the
 geometry of diagrams, Como (1984)}, Springer lecture note in
 Math.  {\bf 1181 } (1986), 159-190. 

\ref{87}Uniqueness and presentation of Kac-Moody groups over
 fields, {\it J. of Algebra} {\bf105} (1987), 542-573.

\ref{88}Sur le groupe des automorphismes de certains groupes de 
Coxeter, {\it J. of Algebra} {\bf 113} (1988), 346-357.

\ref{89}Immeubles jumel\'es I, R\'esum\'e de cours, 
{\it Annuaire du Coll\`ege de France} (1989), 81-96.

\ref{92}Twin buildings and groups of Kac-Moody type, in {\it 
Groups combinatorics and geometry, Durham (1990)}, Liebeck et
 Saxl \'editeurs, London Math. Soc. lecture note {\bf165}, 
Cambridge U. Press (1992), 249-286.

\par Jacques {\petcap Tits} et Richard M. {\petcap Weiss}
\ref{03}{\it Moufang Polygons}, Springer monographs in Math. 
 (2003).

\par Richard M. {\petcap Weiss}
\ref{03}{\it The structure of spherical buildings}, 
(Princeton U. Press, Princeton, 2003).
\ref{08}{\it The structure of affine buildings}, 
(Princeton U. Press, Princeton, 2008).

}

\vskip 10mm
\parni {\petcap \qquad Guy Rousseau}
\parni 
\parni Institut \'Elie Cartan  
\parni Unit\'e mixte de Recherche 7502 
\parni Nancy-Universit\'e, CNRS, INRIA 
\parni B.P 239
\parni 54506  Vand\oe uvre l\`es Nancy Cedex
\parni {\petcap France}

rousseau@iecn.u-nancy.fr

\end